\documentclass{amsart}

\usepackage{mathtools}
\usepackage{bm}
\usepackage{amsthm}
\usepackage{mathrsfs}
\usepackage{cleveref}

\newcommand{\rme}{\mathrm{e}}
\newcommand{\rmf}{\mathrm{f}}
\newcommand{\rmi}{\mathrm{i}}
\renewcommand{\Re}{\operatorname{Re}}

\newcommand{\argmin}{\operatorname*{argmin}}

\numberwithin{equation}{section}

\newtheorem{theorem}{Theorem}[section]
\newtheorem{conj}[theorem]{Conjecture}

\allowdisplaybreaks

\makeatletter
\newcommand{\definetitlefootnote}[1]{%
  \newcommand\addtitlefootnote{%
    \makebox[0pt][l]{$^{*}$}%
    \footnote{\protect\@titlefootnotetext}
  }%
  \newcommand\@titlefootnotetext{\spaceskip=\z@skip $^{*}$#1}%
}
\makeatother

\begin{document}
    \definetitlefootnote{This work was supported by JSPS KAKENHI Grant Numbers 17K18732 (HM), 18K13455 (KS), 19K03556 (MN), NSF DMS-1907583 (YM) and the Math+X grant from the Simons Foundation (YM).}
    \title[Fronts and pulses of the bidomain model]{Asymptotic behavior of fronts and pulses of the bidomain model\addtitlefootnote}
    \author[H.~Matano]{Hiroshi Matano}
    \address[H.~Matano]{Meiji Institute for Advanced Study of Mathematical Sciences, Meiji University, 4-21-1 Nakano, Nakano-ku, Tokyo 164-8525, Japan}
    \email{matano@meiji.ac.jp}
    \author[Y.~Mori]{Yoichiro Mori}
    \address[Y.~Mori]{Department of Mathematics, University of Pennsylvania, Philadelphia 19104, USA}
    \email{y1mori@sas.upenn.edu}
    \author[M.~Nara]{Mitsunori Nara}
    \address[M.~Nara]{Faculty of Science and Engineering, Iwate University, 3-18-34 Ueda, Morioka-shi, Iwate 020-8550, Japan}
    \email{nara@iwate-u.ac.jp}
    \author[K.~Sakakibara]{Koya Sakakibara}
    \address[K.~Sakakibara]{Department of Applied Mathematics, Faculty of Science, Okayama University of Science, 1-1 Ridaicho, Kita-ku, Oakayama-shi, Okayama 700-0005, Japan;
    RIKEN iTHEMS, 2-1 Hirosawa, Wako-shi, Saitama 351-0198, Japan}
    \email{ksakaki@xmath.ous.ac.jp}
    \begin{abstract}
        The bidomain model is the standard model for cardiac electrophysiology. In this paper, we investigate the instability and asymptotic behavior of planar fronts and planar pulses of the bidomain Allen--Cahn equation and the bidomain FitzHugh--Nagumo equation in two spatial dimension. In previous work, it was shown that planar fronts of the bidomain Allen--Cahn equation can become unstable in contrast to the classical Allen--Cahn equation. We find that, after the planar front is destabilized, a rotating zigzag front develops whose shape can be explained by simple geometric arguments using a suitable Frank diagram. We also show that the Hopf bifurcation through which the front becomes unstable can be either supercritical or subcritical, by demonstrating a parameter regime in which a stable planar front and zigzag front can coexist. In our computational studies of the bidomain FitzHugh--Nagumo pulse solution, we show that the pulses can also become unstable much like the bidomain Allen--Cahn fronts. However, unlike the bidomain Allen--Cahn case, the destabilized pulse does not necessarily develop into a zigzag pulse. For certain choice of parameters, the destabilized pulse can disintegrate entirely. These studies are made possible by the development of a numerical scheme that allows for the accurate computation of the bidomain equation in a two dimensional strip domain of infinite extent.
    \end{abstract}
    \maketitle

    \section{Introduction}

    The cardiac bidomain model is the standard mathematical model for cardiac electrophysiology:
    \begin{align*}
        \begin{dcases}
            \frac{\partial u}{\partial t}-f(u,s)=\nabla\cdot(A_\rmi\nabla u_\rmi),\quad u=u_\rmi-u_\rme,\\
            \nabla\cdot(A_\rmi\nabla u_\rmi)+\nabla\cdot(A_\rme\nabla u_\rme)=0,\\
            \frac{\partial s}{\partial t}=g(u,s),
            \quad
            s=(s_1,\ldots,s_G),
        \end{dcases}
    \end{align*}
    where $u_{\rmi,\rme}$ are intracellular/extracellular voltages, $u$ is the transmembrane voltage, $s_1,\ldots,s_G$ are gating variables, and $A_{\rmi,\rme}$ are conductivity tensors.
    Nonlinear terms $f(u,s)$ and $g(u,s)$ are of Hodgkin--Huxley (or FitzHugh--Nagumo) type.
    
    It is challenging to study this model mathematically, which has led many to study the monodomain reduction.
    More precisely, let us assume that the condition $A_\rme=\beta A_\rmi$ holds; intracellular and extracellular anisotropies are proportional.
    Then, the bidomain system reduces to the monodomain system:
    \begin{align}
        \begin{dcases}
            \frac{\partial u}{\partial t}=\nabla\cdot(A_{\mathrm{mono}}\nabla u)+f(u,s),
            \quad
            A_{\mathrm{mono}}=\frac{\beta}{1+\beta}A_{\rmi},\\
            \frac{\partial s}{\partial t}=g(u,s).
        \end{dcases}
        \label{eq:monodomain_FHN}
    \end{align}
    The monodomain reduction is an equation of reaction-diffusion type with Hodgkin--Huxley or FitzHugh--Nagumo type nonlinearities, and thus support traveling pulse solutions and other patterns characteristic of excitable systems \cite{KS98,CPS14}.
    These traveling pulse solutions describe the propagating electrical signal in the heart. Extensive simulations indicate that the bidomain model support qualitatively similar solutions.
    However, the bidomain model is a \textit{quantitatively} better model than the monodomain model, especially under extreme conditions, defibrillation, for instance \cite{keener1996biophysical,KS98,GP13}.
    A natural question arises as to how different the behavior of the bidomain model is qualitatively from that of the monodomain model.
    
    The bidomain model, initially introduced in \cite{EJ70,T78,MRE79}, is
    the standard tissue-level model of cardiac electrophysiology widely used in simulations (see for instance \cite{clayton2011models,H93,KB98,CPT05,PDRVG06}).
    Well-posedness is studied in  \cite{CS02,BCP09,V09,CPS14,GK18,GKK19}. It is possible to derive the bidomain model from an underlying microscropic model through homogenization. Such a calculation was first performed formally in 
    \cite{NK93,KS98} and was given analytical justification in \cite{PSC05}.
    
    Very little is known mathematically of the qualitative properties of the bidomain equation. 
    As discussed earlier, from both mathematical and physiological points of view, it is important to study the traveling front and pulse solutions of the bidomain equations. In \cite{MM16}, the authors study the bidomain Allen--Cahn equation in $\mathbb{R}^2$ (to be introduced shortly), which should be seen as the bidomain analogue of the classical Allen--Cahn equation. The bidomain Allen--Cahen equation supports traveling planar front solutions in every direction much like the classical Allen-Cahn equation. However, in sharp contrast to the classical case, the planar front solutions of the bidomain Allen--Cahn equation was found to be unstable under certain parametric conditions.
    
    The study in \cite{MM16} was of a perturbative character, confined to spectral computations of the linearized operator around the traveling front solution. Nothing is known beyond this perturbative regime. Tentative numerical simulations have shown that a zigzag front appears when the planar front is unstable \cite[Section 6]{MM16}, but the mechanism that determines the shape of zigzag fronts is unknown.
    Furthermore, when we turn our attention to the bidomain FitzHugh--Nagumo equation, nothing is known about the pulse's stability.
    
    In this paper, we perform a computational study of the asymptotic behavior of fronts and pulses in bidomain equations in $\mathbb{R}^2$.
    Here, the bidomain equation refers to the bidomain Allen--Cahn equation
    \begin{align}
        \frac{\partial u}{\partial t}-f(u)
        =
        \nabla\cdot(A_\rmi\nabla u_\rmi)
        =
        -\nabla\cdot(A_\rme\nabla u_\rme),
        \quad
        u=u_\rmi-u_\rme
        \label{eq:bidomain_AC}
    \end{align}
    or the bidomain FitzHugh--Nagumo equation
    \begin{align}
        \begin{dcases}
            \frac{\partial u}{\partial t}-f(u,v)
            =
            \nabla\cdot(A_\rmi\nabla u_\rmi)
            =
            -\nabla\cdot(A_\rme\nabla u_\rme),
            \quad
            u=u_\rmi-u_\rme,\\
            \frac{\partial v}{\partial t}
            =
            g(u,v),
        \end{dcases}
        \label{eq:bidomain_FHN}
    \end{align}
    where $A_\rmi$ and $A_\rme$ are $2\times 2$ positive definite symmetric matrices called the conductivity matrices.
    In this paper, we focus on the case where the nonlinearities are given by
    \begin{alignat*}{2}
        &f(u)=u(1-u)(u-\alpha)
        &\quad
        &\text{in \eqref{eq:bidomain_AC}}\\
        &f(u,v)=u(1-u)(u-\alpha)-v,
        \quad
        g(u,v)=\epsilon(u-\gamma v)
        &\quad
        &\text{in \eqref{eq:bidomain_FHN}},
    \end{alignat*}
    where $\alpha\in(0,1)$, $\epsilon>0$, and $\gamma>0$ are constants.
    
    In \cref{sec:preliminaries}, we briefly summarize the results of \cite{MM16}, which we will use later. In particular, it was shown that the convexity of a suitably defined Frank diagram is closely tied to the stability of planar fronts (see \cref{fig:Frank}). Let us consider a planar front solution propagating in a certain direction. If the Frank diagram is convex in this direction, the front is stable to long wavelength perturbations and is unstable otherwise.
    We will also quote an explicit expression on the asymptotic behavior of the principal eigenvalue to be used as a benchmark for the algorithm we propose for computing the principal eigenvalues.
    We will use the characterization of planar fronts' stability by the Frank diagram to guide our study on the shape of the zigzag fronts.

    \Cref{sec:numerical_scheme} summarizes the different numerical methods used in this paper. In \cite{MM16}, the bidomain equation was simulated on a bounded rectangular domain with periodic boundary conditions. A similar method is described in \cref{subsec:spreading}, which we will use to compute spreading front and pulse solutions. However, such a method is not suitable for a detailed computational study of traveling front or pulse solutions, since these solutions reside in regions of infinite extent.
    
    The planar fronts and the planar pulses are originally defined in the whole plane, but it is not easy to compute them numerically in the whole plane.
    Therefore, we consider the bidomain equations in the strip region, which is infinite in the direction of propagation $\xi$ but periodic in the orthogonal direction $\eta$.
    
    We first apply a coordinate transformation, 
    in the direction of propagation $\xi$, centered at the appropriately defined front or pulse location, mapping the infinite strip into a bounded rectangle. We solve the time-discretized equation in the resulting rectangular domain using finite differences in the modified $\xi$ direction and the Fourier transform in the $\eta$ direction. We use a splitting method and alternate between the evolution governed by the bidomain operator (see \eqref{eq:bidomain_op}) and by the nonlinearities. To obtain second order accuracy in time, Strang splitting is used and the evolution substeps corresponding to the bidomain and nonlinear terms are both solved with second order methods. At each time step, we re-center the coordinate transformation so that the front position is fully numerically resolved. This step requires an interpolation operation from the old to the new grid, for which we employ Lagrange interpolation to minimize the error incurred through this step. We perform a numerical convergence study in \cref{sect:accuracy} to confirm that our numerical scheme is second order accurate in space and time. The dynamics of planar fronts and pulses and their instabilities can now be accurately captured. Based on the above, we can also compute the principal eigenvalues and corresponding eigenfunctions of the linearization around the planar fronts, the numerical results of which are tested against analytical calculations in \cref{sect:accuracy}. Furthermore, we develop algorithms to compute the rotationary fronts to which the planar fronts asymptotically approach by devising a suitable iterative algorithm.

    \Cref{sec:AC} deals with the bidomain Allen--Cahn equation.
    By performing numerical computations of the principal eigenvalues and the planar front, we investigate the relationship between the sign of the real part of the principal eigenvalue, the strip region's width, and the planar front's stability. Stability criteria based on eigenvalue calculations are shown to be consistent with the onset of planar front instabilities as exhibited by the numerical computation of the full bidomain model. 
    
    Our numerical experiments strongly suggest that unstable planar fronts asymptotically approach a rotating zigzag front with a constant translational speed $c_{\xi}^\theta$ and a rotational speed $c_\eta^\theta$. The shape of the zigzag front is characterized by the two angles $\theta_{\mathrm{m}}$ and $\theta_{\mathrm{p}}$ between the $\eta$ axis and the level sets of the zigzag front (see \cref{fig:coarsened_zigzag_front}). The shape and speed of the eventual zigzag front can be predicted from elementary geometric arguments using the Frank diagram. In particular, the angles $\theta_{\mathrm{m}}$ and $\theta_{\mathrm{p}}$ correspond to the contact points between the Frank diagram and its convex hull.
    We numerically verify this geometric prediction by comparing the predicted values of $c_{\xi}^\theta, c_\eta^\theta, \theta_{\mathrm{m}}, \theta_{\mathrm{p}}$ against the values obtained by a full numerical computation. The asymptotic zigzag front shape and its speeds are computed by the iterative algorithm mentioned previously. 
    
    Then, we investigate the relationship between the region where the zigzag front exists, the convex hull of the Frank diagram, and the Frank diagram's curvature.
    We observe the supercritical Hopf bifurcation and the subcritical Hopf bifurcation depending on parameters, and Hopf--Hopf bifurcation at the boundary of these bifurcations.
    The zigzag front caused by the planar front's instability has several peaks at the beginning and finally converges to a shape with a single peak after repeated coarsening.
    We examine the relationship between the number of peaks and their duration and discuss which number of peaks is, in some sense, stable.
    Finally, we discuss the spreading front.
    In the paper \cite{MM16}, the authors predicted that the spreading front would converge to the Wulff shape, and we confirm numerically that this is indeed the case.
    
    \Cref{sec:FHN} deals with the bidomain FitzHugh--Nagumo equation.
    In the bidomain Allen--Cahn equation, the planar front exists in all directions, and we study its qualitative behavior, such as the asymptotic behavior and the stability from several points of view in \cref{sec:AC}.
    Therefore, it is natural to conduct a similar study on the bidomain FitzHugh--Nagumo equation. We first recall that, for given $0<\alpha<1/2$, $\epsilon>0$ must be made small enough for a stable planar pulse solution to exist for the classical FitzHugh-Nagumo equation. We verify that, the planar pluse solution for the classical FitzHugh--Nagumo equations are also planar pulse solutions for the bidomain FitzHugh--Nagumo equations. We are thus led to the study of the stability of the planar pulse solutions and their asymptotic evolution.
    
    Much like the bidomain Allen-Cahn equation, we verify through numerical experiments that the planar pulse solutions are indeed unstable in directions in which the Frank diagram is non-convex. In contrast to the bidomain Allen-Cahn case, however, we find that the destabilized planar pulse solution may not necessarily approach a zigzag pulse for large time. Depending on the parameter values of $\epsilon$ and $\alpha$, the planar pulse either develops into a zigzag pulse or disintegrates completely. We find the pulse disintegrates when $\epsilon$ and $\alpha$ are close to the edge of the parametric region in which planar pulse solutions exist. 
    
    Finally, we examine the spreading pulse.
    We numerically explore the existence and asymptotic shape of the spreading pulse, corresponding to the region where the pulse exists.
    
    In \cref{sec:discussion}, we summarize our results and discuss the possible significance for our results in the study of cardiac arrhythmias.

    \section{Preliminaries}
    \label{sec:preliminaries}

    Before going into the main issue, in this section, we briefly summarize the results obtained in Mori--Matano~\cite{MM16}.
    We apply their results throughout this paper. 

    \subsection{Expressions using pseudo-differential operators}
    \label{subsec:pseudo}

    We can represent bidomain equations as a closed-form for $u$ by using pseudo-differential operators.
    These are useful in computing spreading fronts and spreading pulses in a rectangular region with periodic boundary conditions.

    Denote by $\mathcal{F}$ the two-dimensional Fourier transform; $\mathcal{F}v$ is defined for a function $v$ on $\mathbb{R}^2$ by
    \begin{align*}
        (\mathcal{F}v)(\bm{k})
        =
        \hat{v}(\bm{k})
        \coloneqq
        \int_{\mathbb{R}^2}v(\bm{x})\exp(-\rmi\bm{k}\cdot\bm{x})\,\mathrm{d}\bm{x},
        \quad
        \bm{k}
        =
        (k,l)^\top
        \in
        \mathbb{R}^2.
    \end{align*}
    Applying the Fourier transform $\mathcal{F}$ to the relation $\nabla\cdot(A_\rmi\nabla u_\rmi)=-\nabla\cdot(A_\rme\nabla u_\rme)$, we can represent the term $\nabla\cdot(A_\rmi\nabla u_\rmi)$ using a pseudo-differential operator $\mathcal{L}$ as follows:
    \begin{align}
        \nabla\cdot(A_\rmi\nabla u_\rmi)
        =
        -\mathcal{F}^{-1}Q\mathcal{F}u
        \equiv
        -\mathcal{L}u.
        \label{eq:bidomain_op}
    \end{align}
    Namely, the bidomain operator $\mathcal{L}$ is a Fourier multiplier operator with symbol $Q$ given by
    \begin{align*}
        Q(\bm{k})
        =
        \frac{Q_\rmi(\bm{k})Q_\rme(\bm{k})}{Q_\rmi(\bm{k})+Q_\rme(\bm{k})},
        \quad
        Q_{\rmi,\rme}(\bm{k})
        =
        \bm{k}^\top A_{\rmi,\rme}\bm{k}.
    \end{align*}
    Thus the bidomain Allen--Cahn equation \eqref{eq:bidomain_AC} and the bidomain FitzHugh--Nagumo equation \eqref{eq:bidomain_FHN} are rewritten as
    \begin{align}
        \frac{\partial u}{\partial t}=-\mathcal{L}u+f(u)
        \label{eq:bidomain_AC_L}
    \end{align}
    and
    \begin{align*}
        \begin{dcases}
            \frac{\partial u}{\partial t}=-\mathcal{L}u+f(u,v),\\
            \frac{\partial v}{\partial t}=g(u,v).
        \end{dcases}
    \end{align*}
    A suitable linear transformation on the coordinate system makes it possible to transform the conductivity matrices $A_\rmi$ and $A_\rme$ into the following standard forms:
    \begin{align}
        \begin{split}
            A_\rmi
            &=
            \begin{pmatrix}
                1+b+a&0\\
                0&1+b-a
            \end{pmatrix},\\
            A_\rme
            &=
            \begin{pmatrix}
                1-b-a&0\\
                0&1-b+a
            \end{pmatrix},
        \end{split}
        \quad
        |a\pm b|<1.
        \label{eq:conductivity_matrices}
    \end{align}
    Accordingly, we always assume that the conductivity matrices $A_\rmi$ and $A_\rme$ are given by \eqref{eq:conductivity_matrices}.

    \subsection{Directional bidomain equations}

    We will scrutinize the asymptotic behavior of fronts and pulses in the bidomain equations propagating in several directions. 
    To this end, it is necessary to write down them in a new coordinate $\bm{\xi}=(\xi,\eta)^\top$, which we obtain by rotating the original coordinate $\bm{x}$ by angle $\theta$.
    Consequently, the bidomain operator $\mathcal{L}$ reduces to the form $\mathcal{L}^\theta=\mathcal{F}^{-1}Q^\theta\mathcal{F}$, where
    \begin{align*}
        &Q^\theta(\bm{k})
        =
        \frac{Q_\rmi^\theta(\bm{k})Q_\rme^\theta(\bm{k})}{Q_\rmi^\theta(\bm{k})+Q_\rme^\theta(\bm{k})},
        \quad
        Q_{\rmi,\rme}^\theta(\bm{k})
        =
        \bm{k}^\top A_{\rmi,\rme}^\theta\bm{k},\\
        &A_{\rmi,\rme}^\theta
        =
        R^\theta A_{\rmi,\rme}R^{-\theta},
        \quad
        R^\theta
        =
        \begin{pmatrix}
            \cos\theta&-\sin\theta\\
            \sin\theta&\cos\theta
        \end{pmatrix}.
    \end{align*}
    Here, $\mathcal{F}$ denotes the two-dimensional Fourier transform in $\bm{\xi}$. 
    The $\bm{n}^\theta$-directional bidomain Allen--Cahn equation is given by
    \begin{align*}
        \frac{\partial u}{\partial t}
        =
        -\mathcal{L}^\theta u+f(u),
    \end{align*}
    and the $\bm{n}^\theta$-directional bidomain FitzHugh--Nagumo equation is given by
    \begin{align*}
        \begin{dcases}
            \frac{\partial u}{\partial t}=-\mathcal{L}^\theta u+f(u,v),\\
            \frac{\partial v}{\partial t}=g(u,v),
        \end{dcases}
    \end{align*}
    where $\bm{n}^\theta=(\cos\theta,\sin\theta)^\top$.
    Besides, we define $a_{\rmi,\rme}^\theta$, $b_{\rmi,\rme}^\theta$, and $c_{\rmi,\rme}^\theta$ as follows:
    \begin{align*}
        \begin{pmatrix}
            a_{\rmi,\rme}^\theta&b_{\rmi,\rme}^\theta\\[1ex]
            b_{\rmi,\rme}^\theta&c_{\rmi,\rme}^\theta
        \end{pmatrix}
        \coloneqq
        A_{\rmi,\rme}^\theta.
    \end{align*}

    \subsection{Planar fronts in the bidomain Allen--Cahn equation}

    Let us consider a planar front $u_\rmf^\theta(\xi)=u_\rmf^\theta(\bm{n}^\theta\cdot\bm{x}-c_\rmf^\theta t)$, which propagates in the direction of $\bm{n}^\theta$.
    We impose boundary conditions at infinity:
    \begin{align*}
        \lim_{\xi\to-\infty}u_\rmf^\theta(\xi)=1,
        \quad
        \lim_{\xi\to\infty}u_\rmf^\theta(\xi)=0.
    \end{align*}
    Substituting $u_\rmf^\theta$ into the bidomain Allen-Cahn equation \eqref{eq:bidomain_AC}, we obtain the following ordinary differential equation:
    \begin{align*}
        c_\rmf^\theta\frac{\mathrm{d}u_\rmf^\theta}{\mathrm{d}\xi}+Q(\bm{n}^\theta)\frac{\mathrm{d}^2u_\rmf^\theta}{\mathrm{d}\xi^2}+f(u_\rmf^\theta)=0.
    \end{align*}
    Let $u_\rmf^*$ be the normalized planar front, and $c_\rmf^*$ be its speed; that is, $(u_\rmf^*,c_\rmf^*)$ is a solution to the following boundary value problem:
    \begin{align*}
        \begin{dcases}
            c_\rmf^*\frac{\mathrm{d}u_\rmf^*}{\mathrm{d}\xi}+\frac{\mathrm{d}^2u_\rmf^*}{\mathrm{d}\xi^2}+f(u_\rmf^*)=0,\\
            u_\rmf^*(-\infty)=1,\quad u_\rmf^*(\infty)=0.
        \end{dcases}
    \end{align*}
    This problem is that the planar front in the Allen-Cahn equation satisfies. 
    $c_{\rmf}^*$ is unique, and $u_{\rmf}^*$ is uniquely determined up to translation.
    We can explicitly represent them as
    \begin{align*}
        u_{\rmf}^*(\xi)=\frac{1}{1+\exp(\xi/\sqrt{2})},
        \quad
        c_\rmf^*=\sqrt{2}\left(\frac{1}{2}-\alpha\right),
    \end{align*}
    and as a result, we may express $u_\rmf^\theta$ and $c_\rmf^\theta$ as follows:
    \begin{align*}
        u_\rmf^\theta(\xi)=u_{\rmf}^*\left(\xi/\sqrt{Q(\bm{n}^\theta)}\right),
        \quad
        c_\rmf^\theta=\sqrt{Q(\bm{n}^\theta)}c_\rmf^*.
    \end{align*}

    \subsection{Principal eigenvalues in the bidomain Allen--Cahn equation}

    We will examine planar fronts' stability in the bidomain Allen-Cahn equation associated with the corresponding linearized operator's principal eigenvalue.
    To this end, we introduce a moving coordinate $\bm{\xi}$ that travels with the planar front so that the $\xi$-axis aligns with the direction of propagation $\bm{n}^\theta$ and the $\eta$-axis is parallel to the planar front.
    The bidomain Allen-Cahn equation in this new coordinate becomes
    \begin{align*}
        \frac{\partial u}{\partial t}
        =
        c_\rmf^\theta\frac{\partial u}{\partial\xi}-\mathcal{L}^\theta u+f(u).
    \end{align*}
    We linearize this equation around the traveling front $(u_\rmf^\theta,c_\rmf^\theta)$ to obtain
    \begin{align*}
        \frac{\partial v^\theta}{\partial t}
        =
        c_\rmf^\theta\frac{\partial v^\theta}{\partial\xi}-\mathcal{L}^\theta v^\theta+f'(u_\rmf^\theta)v^\theta
        \equiv
        \mathcal{P}^\theta v^\theta.
    \end{align*}
    Applying the Fourier transform in the $\eta$ direction to this equation, we obtain the $l$th mode's eigenvalue problem for each $l\in\mathbb{R}$:
    \begin{align}
        \mathcal{P}_l^\theta v_l^\theta=\lambda_l^\theta v_l^\theta,
        \quad
        \lambda_0^\theta=0,
        \quad
        v_0^\theta=-\frac{\partial u_{\mathrm{f}}^\theta}{\partial\xi},
        \label{eq:Eig1}
    \end{align}
    where 
    \begin{align*}
        \mathcal{P}_l^\theta
        =
        c_\rmf^\theta\frac{\partial}{\partial\xi}-\mathcal{L}_l^\theta+f'(u_\rmf^\theta),
        \quad
        \mathcal{L}_l^\theta
        =
        \mathcal{F}_\xi^{-1}Q^\theta(k,l)\mathcal{F}_\xi,
    \end{align*}
    and $\mathcal{F}_\xi$ is the one-dimensional Fourier transform in $\xi$.
    We impose the normalizing condition
    \begin{align}
        v_l^\theta(0)=v_0^\theta(0).
        \label{eq:Eig2}
    \end{align}

    The following theorem describes the detailed asymptotic behavior of the principal eigenvalue $\lambda_l^\theta$ as $l$ tends to $0$.
    See \cite[Theorem 4.2]{MM16} for details.

    \begin{theorem}
        \label{thm:2.1}
        There is a $\delta>0$ such that for $|l|<\delta$ there is an eigenvector-eigenvalue pair $(v_l^\theta,\lambda_l^\theta)\in H^2(\mathbb{R})\times\mathbb{C}$ satisfying \eqref{eq:Eig1} and \eqref{eq:Eig2} with the following properties:
        \begin{enumerate}
            \item $\lambda_l^\theta$ is a simple principal eigenvalue of $\mathcal{P}_l^\theta$.
            Namely, there exists a constant $\nu_\delta>0$ independent of $l$ such that
            \begin{align*}
                \Sigma(\mathcal{P}_l^\theta)\setminus\{\lambda_l^\theta\}\subset\{z\in\mathbb{C}\mid\Re z<-\nu_\delta\},
            \end{align*}
            where $\Sigma(\cdot)$ refers to the spectrum of the operator.
            \item $\lambda_l^\theta$ is a $C^2$ function for $l$ and satisfies the following asymptotic expansion:
            \begin{align*}
                \lambda_l^\theta
                =
                \rmi\alpha_1c_\rmf^\theta l-\alpha_0l^2+O(l^3)
                \quad
                \text{as}\ l\to0,
            \end{align*}
            where
            \begin{align*}
                \alpha_0
                &=
                \frac{1}{2}+\frac{1}{2}\left(3a^2\cos^2(2\theta)+2ab\cos(2\theta)-4a^2\sin^2(2\theta)-b^2\right)\\
                &\hspace{150pt}
                -\frac{2a^2\sin^2(2\theta)(b+a\cos(2\theta))^2}{1-(b+a\cos(2\theta))^2},\\
                \alpha_1
                &=
                \frac{2a\sin(2\theta)(b+a\cos(2\theta))}{1-(b+a\cos(2\theta))^2}.
            \end{align*}
        \end{enumerate}
    \end{theorem}

    \subsection{Stability of planar fronts}

    \begin{figure}[tb]
        \centering
        \includegraphics[width=.4\hsize]{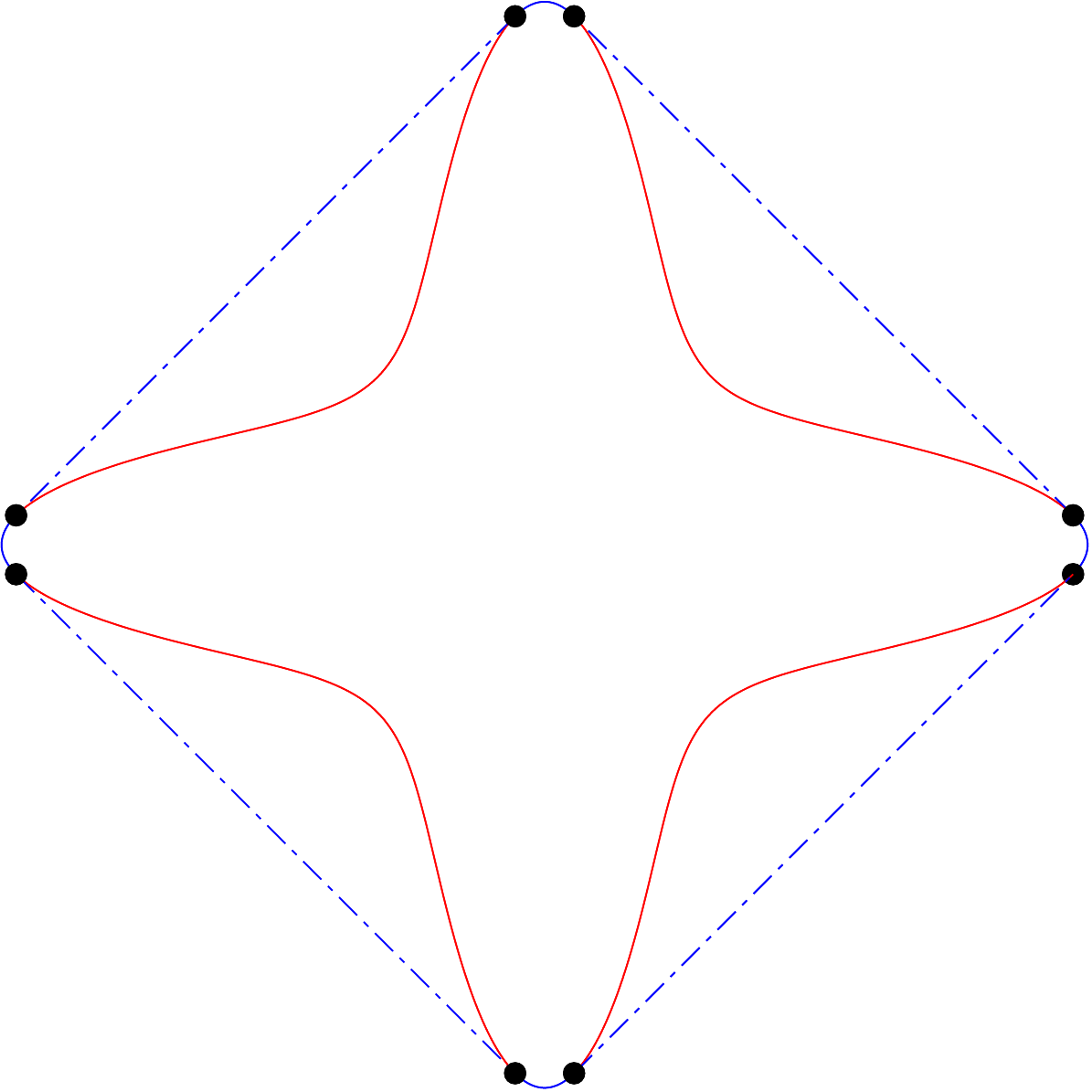}

        \caption{The Frank diagram when $a=0.9$ and $b=0$.
        The blue solid and broken lines represent the Frank diagram's convex hull, and the black dots represent the tangent points of the Frank diagram and the convex hull.}
        \label{fig:Frank}
    \end{figure}

    We will investigate the stability of planar fronts.
    As analytical result, in \cite{MM16}, the authors investigated the stability of planar fronts in the bidomain Allen--Cahn equation in terms of the Frank diagram.

    A \textit{Frank plot} $\mathscr{F}$ is a curve defined as follows~\cite{BCP97,HMR05}:
    \begin{align*}
        \mathscr{F}=\{(\cos\theta,\sin\theta)^\top/\sqrt{Q(\bm{n}^\theta)}\mid\theta\in[0,2\pi]\}.
    \end{align*}
    The figure surrounded by the Frank plot is a \textit{Frank diagram}.
    \Cref{fig:Frank} shows an example of the Frank diagram.
    In particular, when $a>1/2$, the Frank diagram is nonconvex, and there is always an angle at which the planar front becomes unstable~\cite[Eq.~(4.98)]{MM16}.
    More precisely, in \cite{MM16}, the following theorem about the relationship between the planar fronts' stability and the Frank diagram's convexity is presented.

    \begin{theorem}
        \label{thm:stability}
        The planar front propagating in the direction of $\bm{n}^\theta$ where the Frank plot is nonconvex at $(\cos\theta,\sin\theta)^\top/\sqrt{Q(\bm{n}^\theta)}$ is spectrally unstable.
    \end{theorem}

    \section{Numerical scheme}
    \label{sec:numerical_scheme}

    In this section, we develop several numerical methods for solving bidomain equations. 
    Numerical methods for the bidomain FitzHugh-Nagumo equation can be obtained by naturally modifying ones for the bidomain Allen-Cahn equation. 
    Therefore, in the following, we explain our numerical methods for the bidomain Allen-Cahn equation in detail and only make brief comments for the bidomain FitzHugh-Nagumo equation.

    \subsection{Spreading fronts and spreading pulses}
    \label{subsec:spreading}

    We consider the bidomain Allen-Cahn equation \eqref{eq:bidomain_AC_L} in a periodic rectangular region $\Omega=\mathbb{S}_{d_1}^1\times\mathbb{S}_{d_2}^1$, where $\mathbb{S}_d^1\coloneqq\mathbb{R}/d\mathbb{Z}$.
    Namely, we solve the following problem:
    \begin{align}
        \begin{dcases*}
            \frac{\partial u}{\partial t}=-\mathcal{L}u+f(u)&in $\Omega$, $t>0$,\\
            u(\cdot,\cdot,0)=u_0&in $\Omega$,
        \end{dcases*}
        \label{eq:bidomain_AC_spreading}
    \end{align}
    where $u_0$ is the initial value with compact support.
    We adopt the operator splitting method to discretize in time and the Fourier transform to discretize in space to solve this problem.

    The essence of the operator splitting method is to decompose the problem \eqref{eq:bidomain_AC_spreading} into two parts:
    \begin{align}
            \frac{\partial v}{\partial t}=-\mathcal{L}v
            \quad
            \text{in}\ \Omega,\ t>0,
        \label{eq:splitting_1}
    \end{align}
    and
    \begin{align}
        \frac{\partial w}{\partial t}=f(w)
        \quad
        \text{in}\ \Omega,\ t>0.
        \label{eq:splitting_2}
    \end{align}
    Denote time evolution maps, which advance solutions by time $\Delta t$, for \eqref{eq:splitting_1} and \eqref{eq:splitting_2} by $\psi^{\Delta t}$ and $\varphi_f^{\Delta t}$, respectively.
    Here we write the subscript $f$ in $\varphi_f^{\Delta t}$ in order to emphasize that the nonlinearity is $f$.
    According to the operator splitting method of Strang type, a second-order discretization in time, we compute the time evolution as follows:
    \begin{align}
        u^{n+1}=(\varphi_f^{\Delta t/2}\circ\psi^{\Delta t}\circ\varphi_f^{\Delta t/2})u^n,
        \label{eq:Strang_splitting}
    \end{align}
    where $u^n(\cdot)$ is an abbreviation of $u(\cdot,t^n)$, and $t^n=n\Delta t$ denotes the $n$th time step with uniform time increment $\Delta t$.

    \subsubsection*{Construction of time evolution maps}

    We construct $\varphi_f^{\Delta t}$ based on the second-order explicit Runge--Kutta method.
    Namely, define $\varphi_f^{\Delta t}$ as
    \begin{align}
        \varphi_f^{\Delta t}(\mathsf{F})
        =
        \mathsf{F}+\Delta tf(\mathsf{F}_*),
        \quad
        \mathsf{F}_*
        =
        \mathsf{F}+\frac{\Delta t}{2}f(\mathsf{F})
        \label{eq:second-order_RK}
    \end{align}
    for a general function $\mathsf{F}$. 
    Concerning $\psi^{\Delta t}$, we adopt the Fourier transform.
    Namely, define $\psi^{\Delta t}$ as
    \begin{align*}
        \psi^{\Delta t}(\mathsf{F})
        =
        \mathcal{F}_h^{-1}\exp(-Q_h\Delta t)\mathcal{F}_h\mathsf{F}
    \end{align*}
    for a general function $\mathsf{F}$, where $\mathcal{F}_h$ and $\mathcal{F}_h^{-1}$ are the discrete Fourier transform and its inverse, respectively.
    $Q_h$ is a restriction of the Fourier multiplier to the space of discrete wavenumbers.
    Utilizing the above constructed time evolution maps $\varphi_f^{\Delta t}$ and $\psi^{\Delta t}$, we compute the time evolution of the bidomain Allen--Cahn equation \eqref{eq:bidomain_AC_spreading} by \eqref{eq:Strang_splitting}.

    We similarly compute the bidomain FitzHugh--Nagumo equation in a rectangular region with periodic boundary conditions.
    From solution $(u^n,v^n)$ at the $n$th time step $t^n$, we compute the solution $(u^{n+1},v^{n+1})$ at the next time step $t^{n+1}$ using the following procedure:
    \begin{align}
        u^{n+1}=(\varphi_f^{\Delta t/2}\circ\psi^{\Delta t}\circ\varphi_f^{\Delta t/2})(u^n,v^n),
        \quad
        v^{n+1}=\varphi_g^{\Delta t}(u^{n+1},v^n),
        \label{eq:Strang_splitting_FHN}
    \end{align}
    where the operator $\varphi_f^{\Delta t/2}\circ\psi^{\Delta t}\circ\varphi_f^{\Delta t/2}$ acts on the first variable $u^n$, while the operator $\varphi_g^{\Delta t}$ acts on the second variable $v^n$.

    \subsection{Planar/zigzag fronts and planar/zigzag pulses}
    \label{subsec:planar_fronts_pulses}

    The fronts and pulses are defined in the whole plane; therefore, it is natural to compute them in the whole plane with boundary conditions at infinity.
    However, because of computational difficulty, in most previous studies, the bidomain equations are solved in a bounded region, which does not correctly reflect the boundary conditions at infinity, raising questions about numerical solutions' reliability.
    To study the asymptotic behavior of fronts and pulses in the direction of $\bm{n}^\theta$, we solve the $\bm{n}^\theta$-directional bidomain equations in the strip region that is not bounded in $\xi$ direction while periodic in $\eta$ direction.
    To be precise, define the strip region $\mathcal{S}_d$ by
    \begin{align*}
        \mathcal{S}_d
        =
        \mathbb{R}\times\mathbb{S}_d^1,
        \quad
        \mathbb{S}_d^1
        =
        \mathbb{R}/d\mathbb{Z}.
    \end{align*}
    As in the previous section, we adopt the operator splitting method in time discretization.
    We adopt the Fourier transform in the $\eta$ direction and the finite difference method in the $\xi$ direction concerning spatial discretization.

    By eliminating $u_\rme$ from the $\bm{n}^\theta$-directional bidomain Allen--Cahn equation using relation $u=u_\rmi-u_\rme$, we obtain the following partial differential equations for $u$ and $u_\rmi$:
    \begin{align*}
        \begin{dcases*}
            \frac{\partial u}{\partial t}-f(u)=\nabla\cdot(A_\rmi^\theta\nabla u_\rmi)&in $\mathcal{S}_d$, $t>0$,\\
            \nabla\cdot((A_\rmi^\theta+A_\rme^\theta)\nabla u_\rmi)=\nabla\cdot(A_\rme^\theta\nabla u)&in $\mathcal{S}_d$, $t>0$,\\
            u(-\infty,\eta,t)=1,\quad u(\infty,\eta,t)=0&for $\eta\in\mathbb{S}_d^1$, $t>0$\\
            u(\cdot,\cdot,0)=u_0&in $\mathcal{S}_d$.
        \end{dcases*}
    \end{align*}

    As stated above, we adopt the operator splitting method in time discretization.
    Namely, we split the above equations into
    \begin{align}
        \begin{dcases*}
            \frac{\partial v}{\partial t}=\nabla\cdot(A_\rmi^\theta\nabla v_\rmi)&in $\mathcal{S}_d$, $t>0$,\\
            \nabla\cdot((A_\rmi^\theta+A_\rme^\theta)\nabla v_\rmi)=\nabla\cdot(A_\rme^\theta\nabla v)&in $\mathcal{S}_d$, $t>0$,\\
            v(-\infty,\eta,t)=1,\quad v(\infty,\eta,t)=0&for $\eta\in\mathbb{S}_d^1$, $t>0$
        \end{dcases*}
        \label{eq:splitting_AC_split_1}
    \end{align}
    and
    \begin{align}
        \frac{\partial w}{\partial t}=f(w)
        \quad
        \text{in}\ \mathcal{S}_d,\ t>0.
        \label{eq:splitting_AC_split_2}
    \end{align}
    We then adopt the second-order Strang splitting; that is, denoting by $\psi^{\Delta t}$ and $\varphi_f^{\Delta t}$ the time evolution maps, which advance solutions by time $\Delta t$, for \eqref{eq:splitting_AC_split_1} and \eqref{eq:splitting_AC_split_2}, respectively, we compute the time evolution according to \eqref{eq:Strang_splitting}.

    \subsubsection*{Construction of time evolution maps}

    Concerning $\varphi_f^{\Delta t}$, we adopt the second-order explicit Runge--Kutta method.
    Namely, we construct $\varphi_f^{\Delta t}$ by \eqref{eq:second-order_RK}.
    On the other hand, for $\psi^{\Delta t}$, we adopt the trapezoidal rule; that is, we solve the following partial differential equations:
    \begin{align*}
        \begin{dcases*}
            \frac{v-\hat{v}}{\Delta t}=\frac{1}{2}\left[\nabla\cdot(A_\rmi^\theta\nabla(v_\rmi+\hat{v}_\rmi))\right]&in $\mathcal{S}_d$,\\
            \nabla\cdot((A_\rmi^\theta+A_\rme^\theta)\nabla v_\rmi)=\nabla\cdot(A_\rme^\theta\nabla v)&in $\mathcal{S}_d$,\\
            v(-\infty,\eta)=1,\quad v(\infty,\eta)=0&for $\eta\in\mathbb{S}_d^1$,
        \end{dcases*}
    \end{align*}
    where the hat symbol $\hat{\cdot}$ denotes the solution at a previous time step.

    Since all the functions $v$, $v_\rmi$, and $\hat{v}$ are periodic in $\eta$ direction, we can represent them as the Fourier series in the $\eta$ variable with variable coefficients of $\xi$.
    Denote their Fourier coefficients by $\{v_l(\xi)\}_{l\in\mathbb{Z}}$, $\{v_{\rmi,l}(\xi)\}_{l\in\mathbb{Z}}$, and $\{\hat{v}_l(\xi)\}_{l\in\mathbb{Z}}$, respectively.
    We then obtain the following system of ordinary differential equations for each Fourier mode:
    \begin{align*}
        \begin{dcases*}
            \frac{v_l(\xi)-\hat{v}_l(\xi)}{\Delta t}=\frac{a_\rmi^\theta}{2}\left(\frac{\mathrm{d}^2v_{\rmi,l}}{\mathrm{d}\xi^2}(\xi)+\frac{\mathrm{d}^2\hat{v}_{\rmi,l}}{\mathrm{d}\xi^2}(\xi)\right)+b_\rmi^\theta\left(\frac{\mathrm{d}v_{\rmi,l}}{\mathrm{d}\xi}(\xi)+\frac{\mathrm{d}\hat{v}_{\rmi,l}}{\mathrm{d}\xi}(\xi)\right)\rmi w_l\\
            \hspace{180pt}-\frac{c_\rmi^\theta}{2} \left(v_{\rmi,l}(\xi)+\hat{v}_{\rmi,l}(\xi)\right)w_l^2,\quad\xi\in\mathbb{R},\\
            (a_\rmi^\theta+a_\rme^\theta)\frac{\mathrm{d}^2v_{\rmi,l}}{\mathrm{d}\xi^2}(\xi)+2(b_\rmi^\theta+b_\rme^\theta)\frac{\mathrm{d}v_{\rmi,l}}{\mathrm{d}\xi}(\xi)\rmi w_l-(c_\rmi^\theta+c_\rme^\theta)v_{\rmi,l}(\xi)w_l^2\\
            \hspace{120pt}
            =a_\rme^\theta\frac{\mathrm{d}^2v_l}{\mathrm{d}\xi^2}(\xi)+2b_\rme^\theta\frac{\mathrm{d}v_l}{\mathrm{d}\xi}(\xi)\rmi w_l-c_\rme^\theta v_l(\xi)w_l^2,\quad\xi\in\mathbb{R},\\
            v_l(-\infty)=
            \begin{dcases*}
                1&if $l=0$,\\
                0&otherwise,
            \end{dcases*}
            \quad
            v_l(\infty)=0,
        \end{dcases*}
    \end{align*}
    where $w_l=2\pi l/d$ ($l\in\mathbb{Z}$) is a discrete wavenumber.
    To solve the above problem, we need to discretize in $\xi$ direction to approximate derivatives for $\xi$.
    In this subsection, we develop the one-dimensional finite difference method on an unbounded interval.

    We introduce a coordinate transformation $g\colon(-1,1)\rightarrow\mathbb{R}$ by
    \begin{align*}
        g(z)=K\tan\left(\frac{\pi}{2}z\right),
        \quad
        z\in(-1,1),
    \end{align*}
    where $K$ is a positive number.
    We formally extend this function to the one from $[-1,1]$ onto $\overline{\mathbb{R}}=\mathbb{R}\cup\{\pm\infty\}$ by defining $g(-1)=-\infty$ and $g(1)=\infty$ for the sake of convenience.
    We denote the extended function by the same symbol $g$.
    We define a uniform mesh $\{z^j\}_{j=0}^{N_\xi+1}$ and its adjoint $\{\hat{z}^j\}_{j=1}^{N_\xi+1}$ by
    \begin{alignat*}{2}
        z^j
        &=
        -1+j\Delta z,
        &\quad
        &j=0,1,\ldots,N_\xi+1,\\
        \hat{z}^j
        &=
        -1+\left(j-\frac{1}{2}\right)\Delta z,
        &\quad
        &j=1,2,\ldots,N_\xi+1,
    \end{alignat*}
    where $\Delta z=2/(N_\xi+1)$.
    By mapping the uniform mesh $\{z^j\}_{j=0}^{N_\xi+1}$ by the coordinate transformation $g$, we obtain a nonuniform mesh $\{\xi^j\}_{j=0}^{N_\xi+1}$ on $\overline{\mathbb{R}}$ as
    \begin{align*}
        \xi^j=g(z^j),
        \quad
        j=0,1,\ldots,N_\xi+1.
    \end{align*}
    In particular, $\xi^0=-\infty$ and $\xi^{N_\xi+1}=\infty$ hold.
    Based on the chain rule
    \begin{align*}
        \frac{\mathrm{d}\mathsf{F}}{\mathrm{d}\xi}(\xi)
        =
        \frac{1}{g'(z)}\frac{\mathrm{d}}{\mathrm{d}z}(\mathsf{F}(g(z))),
        \quad
        \frac{\mathrm{d}^2\mathsf{F}}{\mathrm{d}\xi^2}(\xi)
        =
        \frac{1}{g'(z)}\frac{\mathrm{d}}{\mathrm{d}z}\left[\frac{1}{g'(z)}\frac{\mathrm{d}}{\mathrm{d}z}(\mathsf{F}(g(z)))\right]
    \end{align*}
    at $\xi=g(z)$ for a function $\mathsf{F}$ defined on $\mathbb{R}$, we approximate the first derivative $(\mathrm{d}\mathsf{F}/\mathrm{d}\xi)(\xi)$ and the second derivative $(\mathrm{d}^2\mathsf{F}/\mathrm{d}\xi^2)(\xi)$ on the nodal points $\xi=\xi^j$ ($j=1,2,\ldots,N_\xi$) by the central finite differences:
    \begin{align*}
        \frac{\mathrm{d}\mathsf{F}}{\mathrm{d}\xi}(\xi^j)
        &\approx
        \frac{1}{2}\left[\frac{1}{g'(\hat{z}^{j+1})}\frac{\mathsf{F}^{j+1}-\mathsf{F}^j}{\Delta z}+\frac{1}{g'(\hat{z}^j)}\frac{\mathsf{F}^j-\mathsf{F}^{j-1}}{\Delta z}\right]
        \eqqcolon
        \overline{\delta}\mathsf{F}^j,\\
        \frac{\mathrm{d}^2\mathsf{F}}{\mathrm{d}\xi^2}(\xi^j)
        &\approx
        \frac{1}{g'(z^j)\Delta z}\left[\frac{1}{g'(\hat{z}^{j+1})}\frac{\mathsf{F}^{j+1}-\mathsf{F}^j}{\Delta z}-\frac{1}{g'(\hat{z}^j)}\frac{\mathsf{F}^j-\mathsf{F}^{j-1}}{\Delta z}\right]
        \eqqcolon
        \overline{\delta}^2\mathsf{F}^j,
    \end{align*}
    where $\mathsf{F}^j\coloneqq\mathsf{F}(\xi^j)=\mathsf{F}(g(z^j))$.
    We then obtain the following linear system: for each $l\in\mathbb{Z}$,
    \begin{align*}
        \begin{dcases}
            \frac{v_l^j-\hat{v}_l^j}{\Delta t}
            =
            \frac{a_\rmi^\theta}{2}\left(\overline{\delta}^2v_{\rmi,l}^j+\overline{\delta}^2\hat{v}_{\rmi,l}^j\right)+b_\rmi^\theta\left(\overline{\delta}v_{\rmi,l}^j+\overline{\delta}\hat{v}_{\rmi,l}^j\right)\rmi w_l-\frac{c_\rmi^\theta}{2}\left(v_{\rmi,l}^j+\hat{v}_{\rmi,l}^j\right)w_l^2,\\
            \hspace{250pt}j=1,2,\ldots,N_\xi,\\
            (a_\rmi^\theta+a_\rme^\theta)\overline{\delta}^2v_{\rmi,l}^j+2(b_\rmi^\theta+b_\rme^\theta)\overline{\delta}v_{\rmi,l}^j\rmi w_l-(c_\rmi^\theta+c_\rme^\theta)v_{\rmi,l}^jw_l^2\\
            \hspace{100pt}
            =
            a_\rme^\theta\overline{\delta}^2v_l^j+2b_\rme^\theta\overline{\delta}v_l^j\rmi w_l-c_\rme^\theta v_l^jw_l^2,
            \quad
            j=1,2,\ldots,N_\xi,\\
            v_l^0=
            \begin{dcases*}
                1&if $l=0$,\\
                0&otherwise,
            \end{dcases*}
            \quad
            v_l^{N_\xi+1}=0.
        \end{dcases}
    \end{align*}
    
    \subsubsection*{Regridding}
    
    To study the long-time behavior of the front, it becomes necessary to re-grid as the front advances.
    First, we estimate the position of the $1/2$-level set by using third-order Lagrange interpolation and set $\xi=0$ there.
    The values of $u$ in the new coordinate system are then defined using the quadratic interpolation.

    For the bidomain FitzHugh--Nagumo equation, we compute $(u^{n+1},v^{n+1})$ from $(u^n,v^n)$ as in \eqref{eq:Strang_splitting_FHN}.

    \subsection{Principal eigenvalues and corresponding eigenfunctions}
    \label{subsec:eigenvalue}

    Although we know that the principal eigenvalue $\lambda_0$ at $l=0$ being equal to $0$ and the asymptotic behavior of principal eigenvalues where $|l|$ is sufficiently small by \cref{thm:2.1}, in order to discuss the stability of the planar front in the bidomain Allen--Cahn equation, we need to compute the principal eigenvalues beyond the range covered by the theorem.
    We solve the eigenvalue problem \eqref{eq:Eig1} by the same strategy in the previous subsection.
    Namely, we represent the eigenfunction $v^\theta$ as the Fourier series in $\eta$ with variable Fourier coefficients of $\xi$ and approximate derivatives concerning $\xi$ by the central finite differences.
    As a result, we obtain the following linear system for the $l$th mode's eigenvalue problem:
    \begin{align*}
        \begin{dcases}
            \lambda_l^\theta v_l^j
            =
            c_\rmf^\theta\overline{\delta}v_{\rmi,l}^j+a_\rmi^\theta\overline{\delta}^2v_{\rmi,l}^j+2b_\rmi^\theta\overline{\delta}v_{\rmi,l}^j\rmi w_l-c_\rmi^\theta v_{\rmi,l}^jw_l^2+f'(u_\rmf^j)v_l^j,
            \quad
            j=1,2,\ldots,N_\xi,\\
            (a_\rmi^\theta+a_\rme^\theta)\overline{\delta}^2v_{\rmi,l}^j+2(b_\rmi^\theta+b_\rme^\theta)\overline{\delta}v_{\rmi,l}^j\rmi w_l-(c_\rmi^\theta+c_\rme^\theta)v_{\rmi,l}^jw_l^2\\
            \hspace{80pt}=
            a_\rme^\theta\overline{\delta}^2v_l+2b_\rme^\theta\overline{\delta}v_{\rmi,l}^j\rmi w_l-(c_\rmi^\theta+c_\rme^\theta)v_{\rmi,l}^jw_l^2,
            \quad
            j=1,2,\ldots,N_\xi,\\
            v_l^0=v_l^{N_\xi+1}=v_{\rmi,l}^0=v_{\rmi,l}^{N_\xi+1}=0,
            \quad
            \sum_{j=1}^{N_\xi}|v_l^j|^2=1.
        \end{dcases}
    \end{align*}
    We know from \eqref{eq:Eig1} the analytical expressions of the principal eigenvalue and corresponding eigenfunction at $l=0$.
    Hence, we solve the above problem using the Newton method while slightly increasing $l$.

    \subsection{Iterative method for asymptotic shape of fronts in the bidomain Allen--Cahn equation}
    \label{subsec:iterative_method}

    We will investigate the existence/nonexistence of zigzag and usual planar fronts in the bidomain Allen--Cahn equation.
    To this end, we develop an iterative algorithm for computing the asymptotic shape of the front in the bidomain Allen--Cahn equation without computing its time evolution.

    Let $u^\theta=u^\theta(\xi,\eta)$ be the front in the bidomain Allen--Cahn equation, and let $c_\xi^\theta$ and $c_\eta^\theta$ be its speed in $\xi$ and $\eta$ directions, respectively.
    Namely, $(u^\theta,c_\xi^\theta,c_\eta^\theta)$ satisfies the following partial differential equations:
    \begin{align}
        \begin{dcases*}
            c_\xi^\theta\frac{\partial u^\theta}{\partial\xi}+c_\eta^\theta\frac{\partial u^\theta}{\partial\eta}+\nabla\cdot(A_\rmi^\theta\nabla u_\rmi^\theta)+f(u^\theta)=0&in $\mathcal{S}_d$,\\
            \nabla\cdot((A_\rmi^\theta+A_\rme^\theta)\nabla u_\rmi^\theta)=\nabla\cdot(A_\rme^\theta\nabla u^\theta)&in $\mathcal{S}_d$,\\
            u^\theta(-\infty,\eta)=1,\quad u^\theta(\infty,\eta)=0&for $\eta\in\mathbb{S}_d^1$.
        \end{dcases*}
        \label{eq:AC_front}
    \end{align}
    To solve this problem, define a constant $f_0$ as
    \begin{align*}
        f_0\coloneqq\frac{f'(0)+f'(1)}{2},
    \end{align*}
    which is negative.
    Then, we rewrite the first equation in \eqref{eq:AC_front} as follows:
    \begin{align*}
        c_\xi^\theta\frac{\partial u^\theta}{\partial\xi}+c_\eta^\theta\frac{\partial u^\theta}{\partial\eta}+\nabla\cdot(A_\rmi^\theta\nabla u_\rmi^\theta)+f_0u^\theta
        =
        -f(u^\theta)+f_0u^\theta.
    \end{align*}
    Performing the Fourier transform in $\eta$ direction, we obtain the following system: for each $\xi\in\mathbb{R}$ and $l\in\mathbb{Z}$,
    \begin{align}
        \begin{dcases*}
            c_\xi^\theta\frac{\mathrm{d}u_l^\theta}{\mathrm{d}\xi}(\xi)+c_\eta^\theta u_l^\theta(\xi)\rmi w_l+a_\rmi^\theta\frac{\mathrm{d}^2u_{\rmi,l}^\theta}{\mathrm{d}\xi^2}(\xi)+2b_\rmi^\theta\frac{\mathrm{d}u_{\rmi,l}^\theta}{\mathrm{d}\xi}(\xi)\rmi w_l\\
            \hspace{50pt}
            -c_\rmi^\theta u_{\rmi,l}^\theta(\xi)w_l^2+f_0u_l^\theta(\xi)=-\mathcal{F}_\eta[f(u^\theta(\xi,\cdot))]_l+f_0u_l^\theta(\xi),\\
            (a_\rmi^\theta+a_\rme^\theta)\frac{\mathrm{d}^2u_{\rmi,l}^\theta}{\mathrm{d}\xi^2}(\xi)+2(b_\rmi^\theta+b_\rme^\theta)\frac{\mathrm{d}u_{\rmi,l}^\theta}{\mathrm{d}\xi}(\xi)\rmi w_l-(c_\rmi^\theta+c_\rme^\theta)u_{\rmi,l}^\theta(\xi)w_l^2\\
            \hspace{90pt}
            =
            a_\rme^\theta\frac{\mathrm{d}^2u_l^\theta}{\mathrm{d}\xi^2}(\xi)+2b_\rme^\theta\frac{\mathrm{d}u_l^\theta}{\mathrm{d}\xi}(\xi)\rmi w_l-c_\rme^\theta u_l^\theta(\xi)w_l^2,\\
            u_l^\theta(-\infty)
            =
            \begin{dcases*}
                1&if $l=0$,\\
                0&otherwise,
            \end{dcases*}
            \quad
            u_l^\theta(\infty)
            =
            0.
        \end{dcases*}
        \label{eq:AC_front_rewrite}
    \end{align}
    To determine the front $u^\theta$ and the speeds $c_\xi^\theta$ and $c_\eta^\theta$, we have to add two more conditions.
    By integrating the first equation in \eqref{eq:AC_front} in $\xi$ and $\eta$ directions, we obtain $c_\xi^\theta$'s explicit representation as
    \begin{align}
        c_\xi^\theta
        =
        \int_{-\infty}^\infty\mathcal{F}_\eta[f(u^\theta(\xi,\cdot))]_0\,\mathrm{d}\xi.
        \label{eq:cxi_representation}
    \end{align}
    We compute $c_\eta^\theta$ as a unique minimizer of the following optimization problem:
    \begin{align}
        c_\eta^\theta
        =
        \argmin_{c_\eta\in\mathbb{R}}\left\|c_\xi^\theta\frac{\partial u^\theta}{\partial\xi}+c_\eta\frac{\partial u^\theta}{\partial\eta}+\nabla\cdot(A_\rmi^\theta\nabla u_\rmi^\theta)+f(u^\theta)\right\|_{L^2}^2,
        \label{eq:ceta_minimization}
    \end{align}
    where we define the $L^2$ norm $\|\mathsf{F}\|_{L^2}$ for a function $\mathsf{F}$ defined on $\mathcal{S}_d$ as
    \begin{align*}
        \|\mathsf{F}\|_{L^2}
        \coloneqq
        \sqrt{(\mathsf{F},\mathsf{F})_{L^2}},
        \quad
        (\mathsf{F},\mathsf{G})_{L^2}
        \coloneqq
        \int_{-\infty}^\infty\int_0^d\mathsf{F}(\xi,\eta)\overline{\mathsf{G}}(\xi,\eta)\,\mathrm{d}\eta\,\mathrm{d}\xi.
    \end{align*}

    Summarizing the above, we compute the front $u^\theta$ and its speeds $c_\xi^\theta$ and $c_\eta^\theta$ corresponding to angle $\theta$ from the ones for angle $\tilde{\theta}$ by the following procedure, where $\tilde{\theta}$ is close to $\theta$.
    \begin{enumerate}
        \item For $n=0$, define $u^{(0)}\coloneqq u^{\tilde{\theta}}$, $c_\xi^{(0)}\coloneqq c_\xi^{\tilde{\theta}}$, and $c_\eta^{(0)}\coloneqq c_\eta^{\tilde{\theta}}$.
        \item For $n\ge1$, we compute $(u^{(n)},c_\xi^{(n)},c_\eta^{(n)})$ from $(u^{(n-1)},c_\xi^{(n-1)},c_\eta^{(n-1)})$ by the following three steps.
            \begin{enumerate}
                \item Compute $u^{(n)}$ as the solution to the following problem in view of \eqref{eq:AC_front_rewrite}:
                    \begin{align*}
                        \begin{dcases}
                            c_\xi^{(n-1)}\frac{\mathrm{d}u_l^{(n)}}{\mathrm{d}\xi}(\xi)+c_\eta^{(n-1)}u_l^{(n)}(\xi)\rmi w_l+a_\rmi^\theta\frac{\mathrm{d}^2u_{\rmi,l}^{(n)}}{\mathrm{d}\xi^2}(\xi)+2b_\rmi^\theta\frac{\mathrm{d}u_{\rmi,l}^{(n)}}{\mathrm{d}\xi}(\xi)\rmi w_l\\
                            \hspace{20pt}
                            -c_\rmi^\theta u_{\rmi,l}^{(n)}(\xi)w_l^2+f_0u_l^{(n)}(\xi)
                            =
                            -\mathcal{F}_\eta[f(u^{(n-1)}(\xi,\cdot))]_l+f_0u_l^{(n-1)}(\xi),\\
                            (a_\rmi^\theta+a_\rme^\theta)\frac{\mathrm{d}^2u_{\rmi,l}^{(n)}}{\mathrm{d}\xi^2}(\xi)+2(b_\rmi^\theta+b_\rme^\theta)\frac{\mathrm{d}u_{\rmi,l}^{(n)}}{\mathrm{d}\xi}(\xi)\rmi w_l-(c_\rmi^\theta+c_\rme^\theta)u_{\rmi,l}^{(n)}(\xi)w_l^2\\
                            \hspace{50pt}
                            =
                            a_\rme^\theta\frac{\mathrm{d}^2u_l^{(n)}}{\mathrm{d}\xi^2}(\xi)+2b_\rme^\theta\frac{\mathrm{d}u_l^{(n)}}{\mathrm{d}\xi}(\xi)\rmi w_l-c_\rme^\theta u_l^{(n)}(\xi)w_l^2,\\
                            u_l^{(n)}(-\infty)=
                            \begin{dcases*}
                                1&if $l=0$,\\
                                0&otherwise,
                            \end{dcases*}
                            \quad
                            u_l^{(n)}(\infty)=0.
                        \end{dcases}
                    \end{align*}
                \item Compute $c_\xi^{(n)}$ by \eqref{eq:cxi_representation}; that is,
                    \begin{align*}
                        c_\xi^{(n)}
                        =
                        \int_{-\infty}^\infty\mathcal{F}_\eta[f(u^{(n)}(\xi,\cdot))]_0\,\mathrm{d}\xi.
                    \end{align*}
                \item Compute $c_\eta^{(n)}$ as a unique minimizer of \eqref{eq:ceta_minimization}; that is,
                    \begin{align*}
                        c_\eta^{(n)}
                        \coloneqq
                        \argmin_{c_\eta\in\mathbb{R}}
                        \left\|c_\xi^{(n)}\frac{\partial u^{(n)}}{\partial\xi}+c_\eta\frac{\partial u^{(n)}}{\partial\eta}+\nabla\cdot(A_\rmi^\theta\nabla u_\rmi^{(n)})+f(u^{(n)})\right\|_{L^2}^2.
                    \end{align*}
            \end{enumerate}
        \item Define $u^\theta$, $c_\xi^\theta$, and $c_\eta^\theta$ as limits of $u^{(n)}$, $c_\xi^{(n)}$, and $c_\eta^{(n)}$, respectively.
        Namely,
        \begin{align*}
            u^\theta=\lim_{n\to\infty}u^{(n)},
            \quad
            c_\xi^\theta=\lim_{n\to\infty}c_\xi^{(n)},
            \quad
            c_\eta^{\theta}=\lim_{n\to\infty}c_\eta^{(n)}.
        \end{align*}
    \end{enumerate}
    In actual computation, we approximate derivatives for $\xi$ by central finite differences developed in \cref{subsec:planar_fronts_pulses} and integration for $\xi$ by the trapezoidal rule.
    
    \subsection{Accuracy of numerical scheme}\label{sect:accuracy}
    
    \begin{figure}[tb]
        \begin{minipage}{.5\hsize}
            \centering
            \includegraphics[width=\hsize]{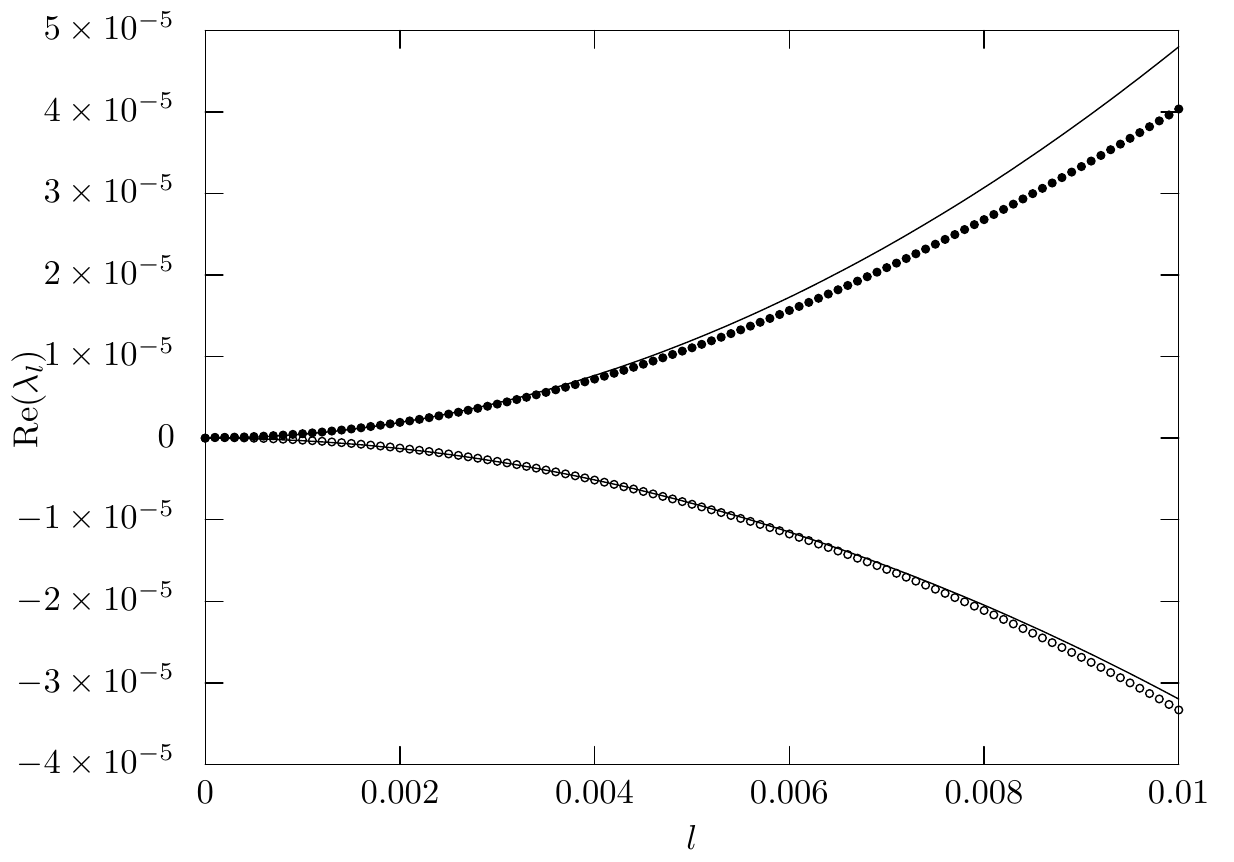}
        \end{minipage}%
        \begin{minipage}{.5\hsize}
            \centering
            \includegraphics[width=\hsize]{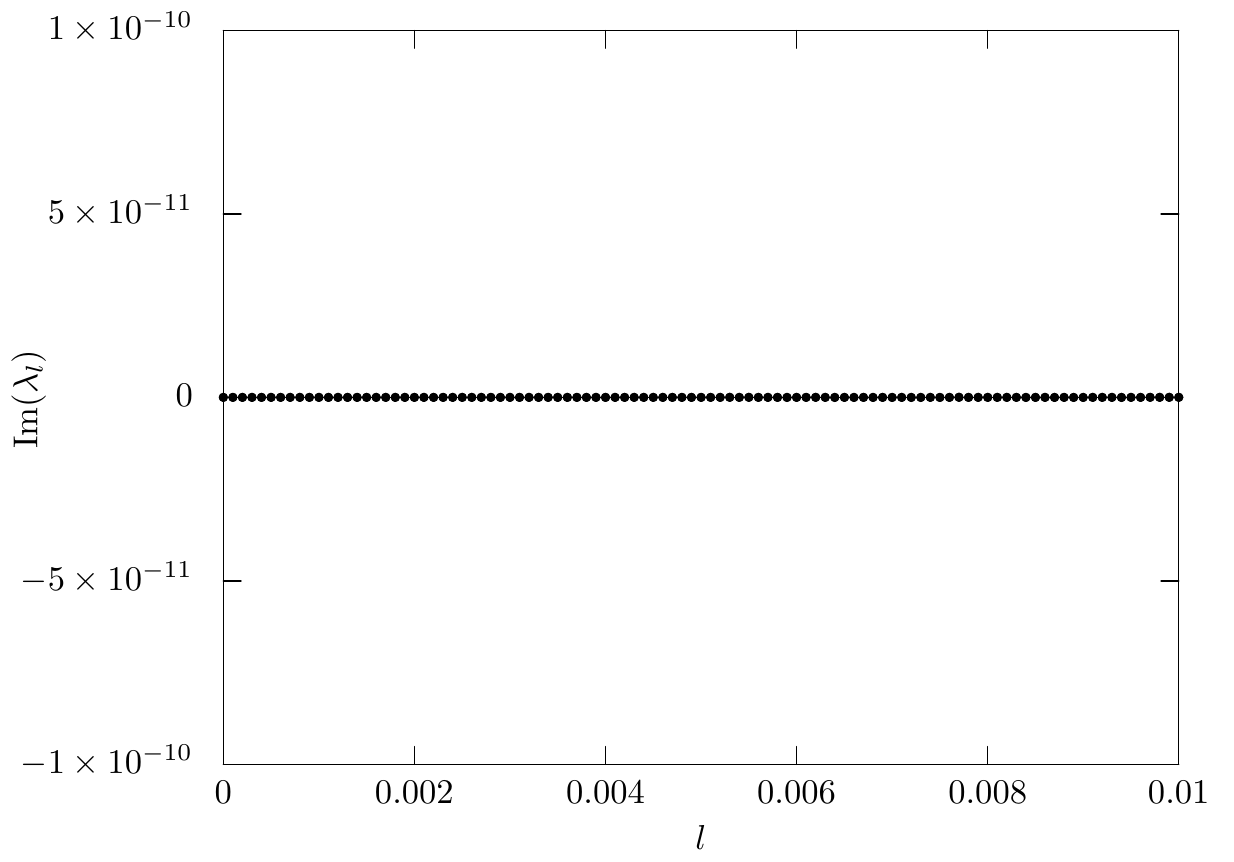}
        \end{minipage}%

        \begin{minipage}{.5\hsize}
            \centering
            \includegraphics[width=\hsize]{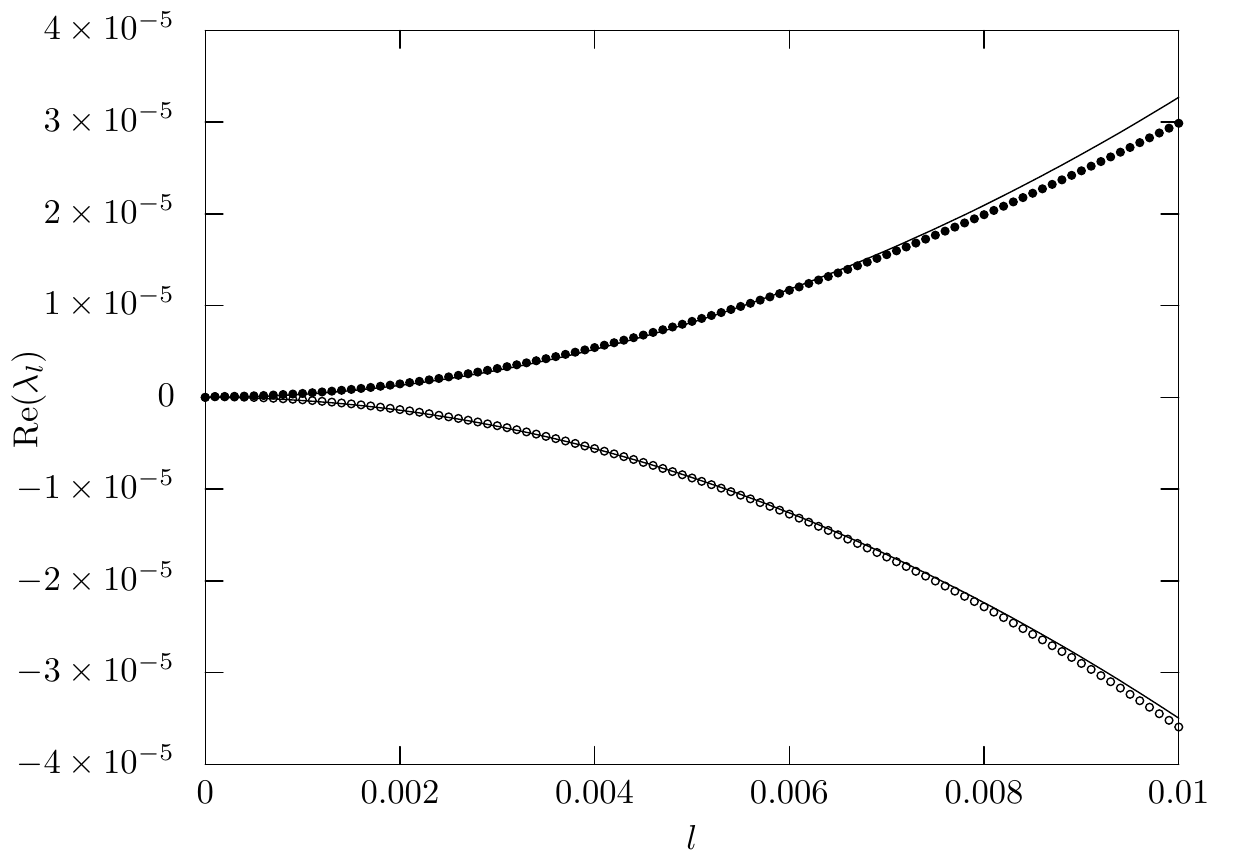}
        \end{minipage}%
        \begin{minipage}{.5\hsize}
            \centering
            \includegraphics[width=\hsize]{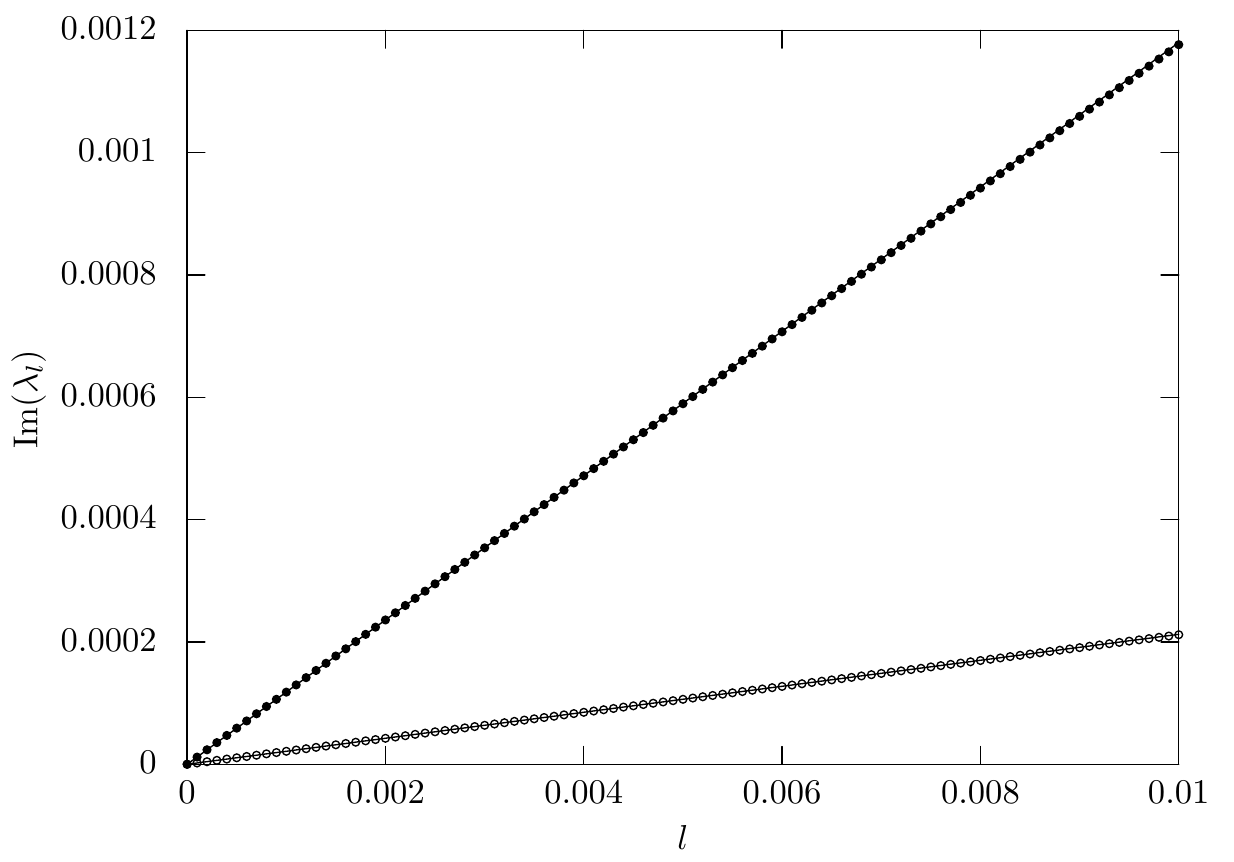}
        \end{minipage}%
        \caption{Asymptotic behavior of principal eigenvalues and comparison with the numerical results.
        The first row shows the result for $\theta=\pi/4$ and the second row for $\theta=\pi/5$, while the first column shows the real parts of the principal eigenvalues and the second row their imaginary parts.
        Dots and circles represent the numerical results.
        The solid lines drawn near them correspond to the principal eigenvalues' asymptotic behavior shown in Theorem \ref{thm:2.1}.}
        \label{fig:comparison_eigenvalue}
    \end{figure}%
    \begin{figure}[tb]
        \begin{minipage}{.4\hsize}
            \centering
            \includegraphics[width=\hsize]{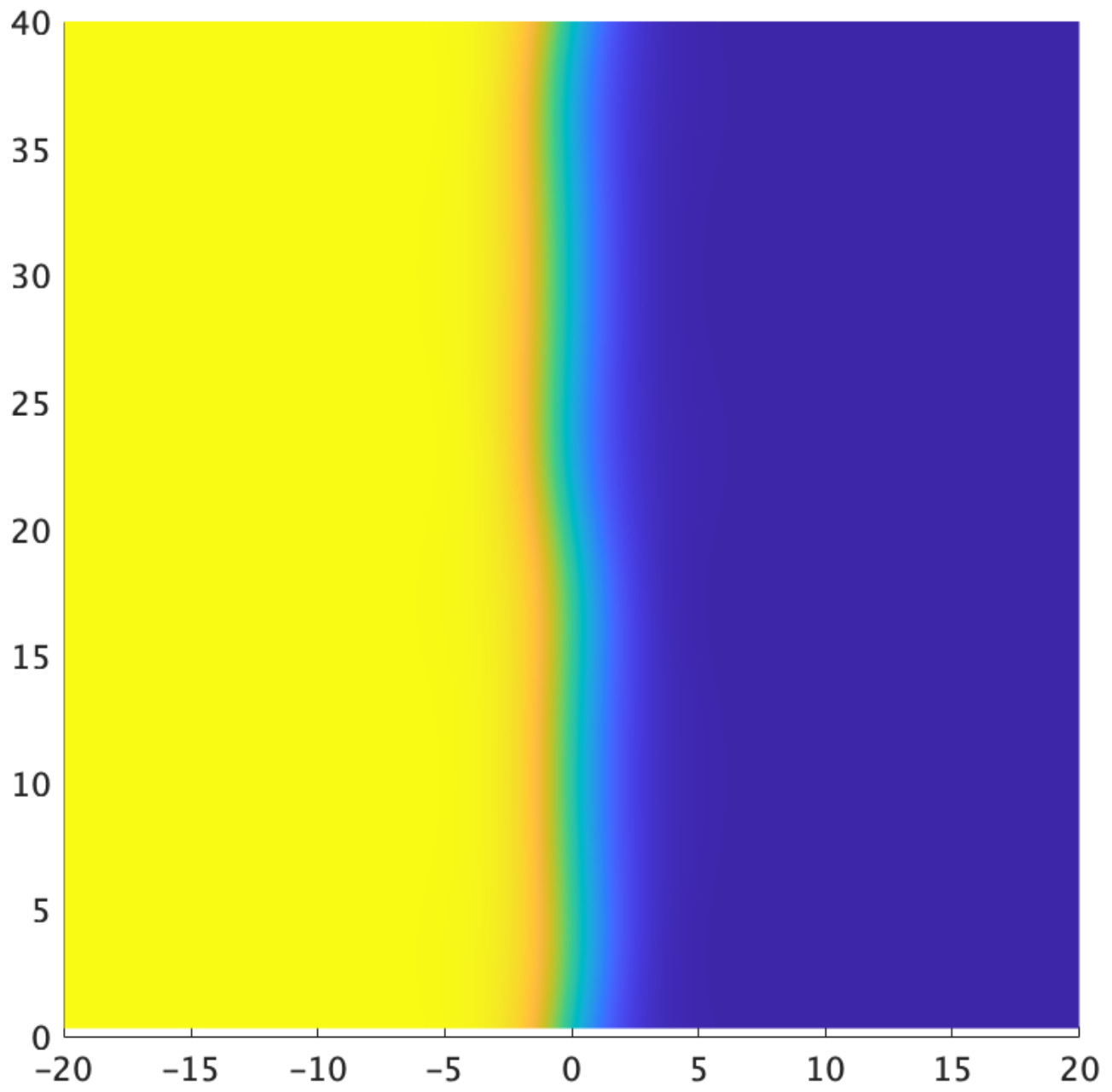}
            
            (a) initial value
        \end{minipage}%
        \begin{minipage}{.6\hsize}
            \centering
            \includegraphics[width=\hsize]{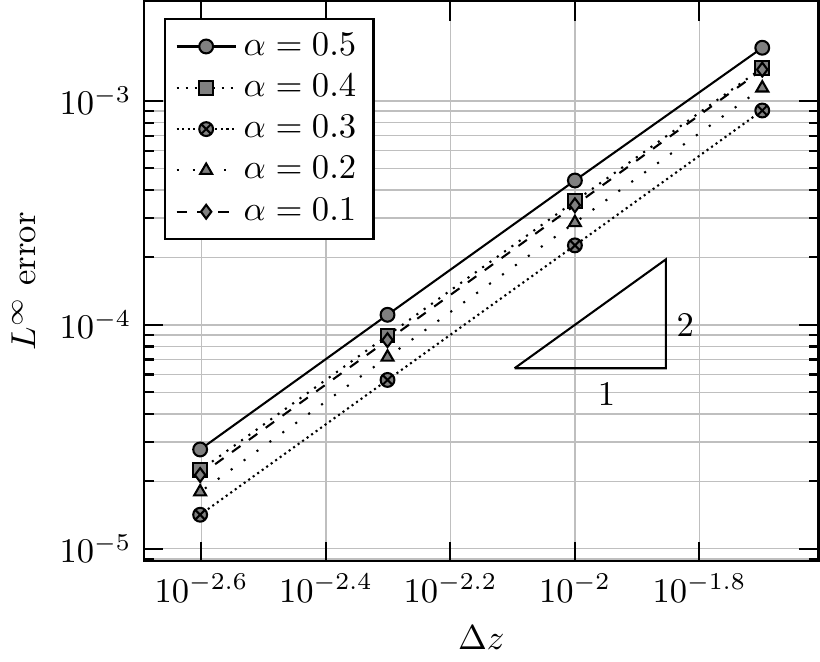}
            
            (b) $L^\infty$ error
        \end{minipage}%
        \caption{Numerical investigation of accuracy of numerical scheme for traveling fronts.}
        \label{fig:Linf_error}
    \end{figure}
    
    We first check whether the numerical scheme developed in \cref{subsec:eigenvalue} works correctly by comparing it with the principal eigenvalues' asymptotic behavior in \cref{thm:2.1}.
    We can observe from \cref{fig:comparison_eigenvalue} that when $l$ is small, the numerical results for the real and imaginary parts of the principal eigenvalues are in good agreement with the asymptotic behavior shown in \cref{thm:2.1}, and in this sense, we can conclude that our numerical method works well.
    
    We also check the accuracy of the scheme developed in \cref{subsec:planar_fronts_pulses} numerically.
    Since we employ the operator splitting method of Strang type for time discretization, we expect temporal accuracy to be second order.
    Concerning the spatial discretization, we use central finite-difference approximations in the $\xi$ direction and spectral discretizations in the $\eta$ direction, so we also expect second order accuracy. We emphasize, however, that our numerical scheme is not conventional in that we are performing a coordinate transformation to map an infinite domain to a finite one, and that we are inserting a regridding/interpolation operation at each time step. It is thus of importance to numerically verify the expected accuracy of our numerical scheme.
    Using the function shown in \cref{fig:Linf_error} (a) as the initial value, we change $N_\xi$ to $99$, $199$, $399$, $799$, $1599$, i.e., $\Delta z$ to $1/50$, $1/100$, $1/200$, $1/400$, $1/800$, and $\Delta t=\Delta z$, and compare the $L^\infty$ error of the solution at $t=20$.
    Here, the $L^\infty$ error $L^\infty(\Delta z)$ for mesh size $\Delta z$ is estimated by 
    \begin{align*}
        L^\infty(\Delta z)
        \coloneqq
        \|u_{\Delta z}-u_{\Delta z/2}\|_{L^\infty(\mathcal{S}_d)},
    \end{align*}
    where $u_{\Delta z}$ denotes numerical solution with spatial mesh size $\Delta z$.
    \Cref{fig:Linf_error} (b) depicts results of numerical experiments.
    We vary the parameter $\alpha$, which controls the speed of fronts, and observe that the accuracy of numerical scheme is indeed second order.

    \section{Bidomain Allen--Cahn equation}
    \label{sec:AC}

    This section investigates the asymptotic behavior and stability of fronts in the bidomain Allen--Cahn equation.

    \subsection{Stability of planar fronts}
    \label{subsec:stability_front}

    \begin{figure}[tb]
        \centering
        \includegraphics[width=.5\hsize]{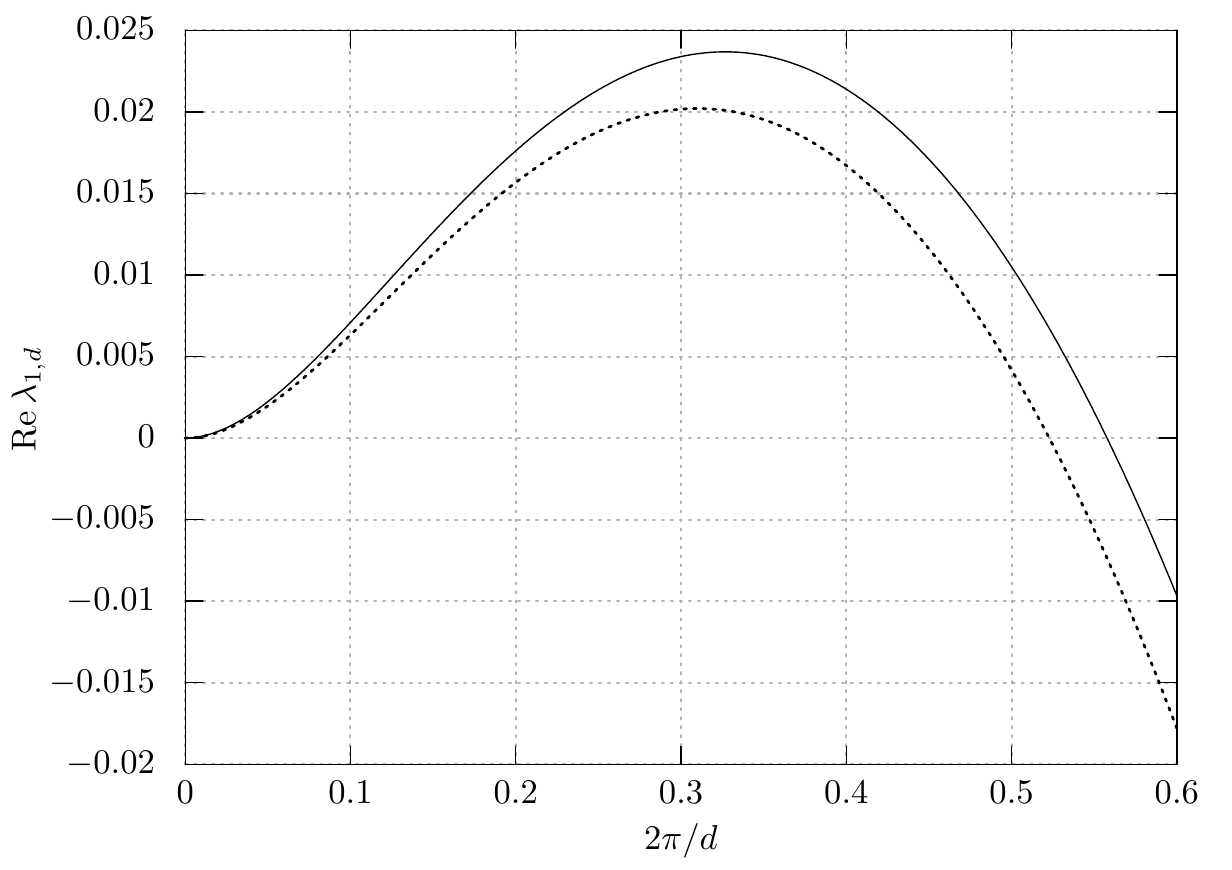}

        \caption{Graphs of $\Re\lambda_{1,d}$ as functions of $d$.
        The solid line shows the result for $\theta=\pi/4$, and the broken line does the result for $\theta=\pi/5$.
        We choose the other parameters as $a=0.9$ and $b=0$.}
        \label{fig:lambda}
    \end{figure}%
    \begin{figure}[tb]
        \begin{minipage}{.25\hsize}
            \centering
            \includegraphics[width=\hsize]{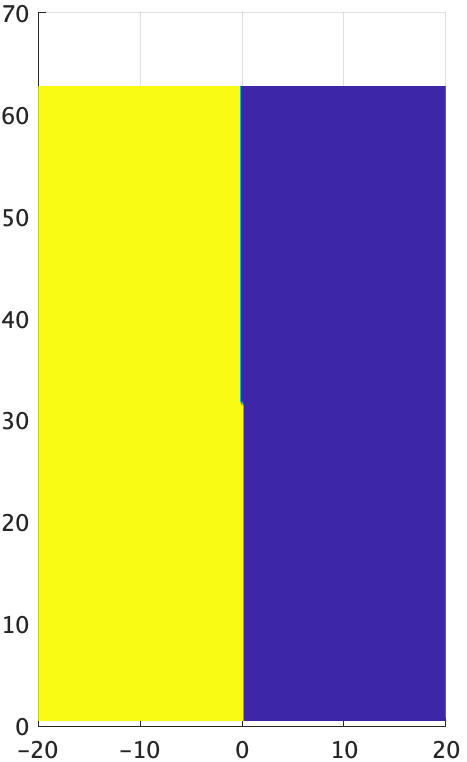}

            $t=0$
        \end{minipage}%
        \begin{minipage}{.25\hsize}
            \centering
            \includegraphics[width=\hsize]{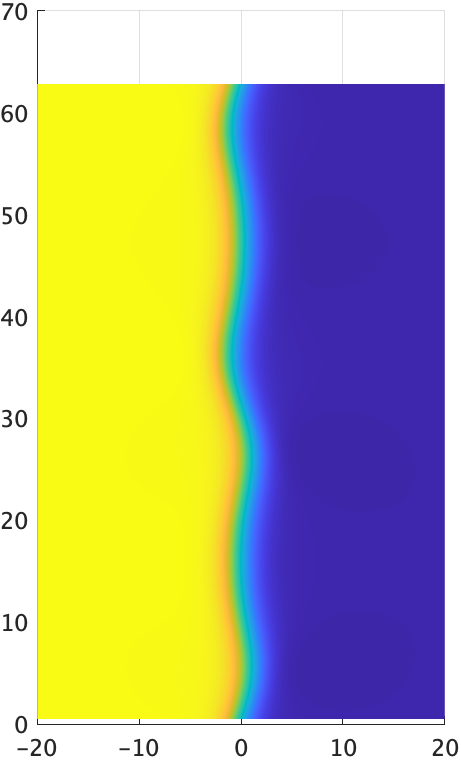}

            $t=100$
        \end{minipage}%
        \begin{minipage}{.25\hsize}
            \centering
            \includegraphics[width=\hsize]{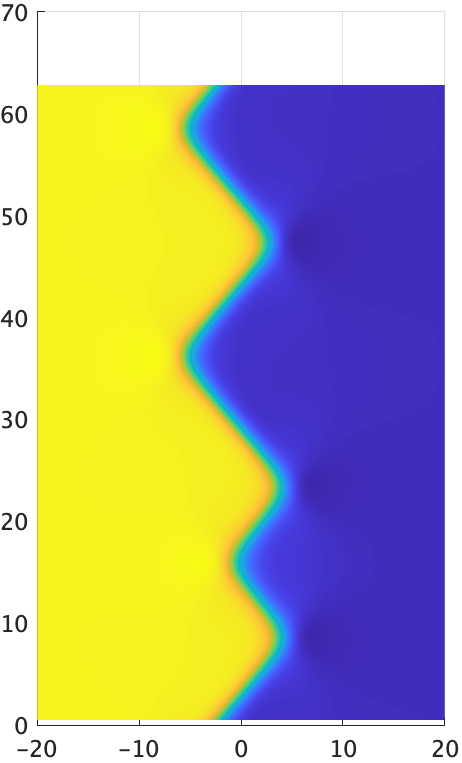}

            $t=500$
        \end{minipage}%
        \begin{minipage}{.25\hsize}
            \centering
            \includegraphics[width=\hsize]{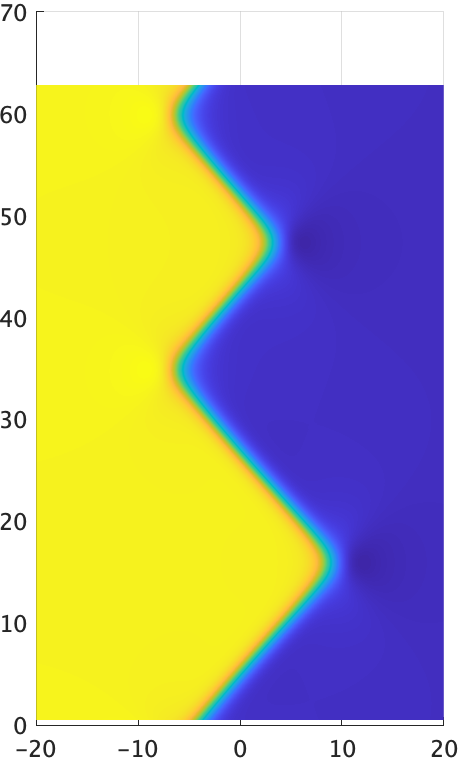}

            $t=1000$
        \end{minipage}%
        \caption{Numerical computations of the bidomain Allen--Cahn equation when $d=2\pi/0.1$, $\theta=\pi/4$, $a=0.9$, $b=0$, and $\alpha=0.4$.}
        \label{fig:pi4_01}
    \end{figure}%
    \begin{figure}[tb]
        \begin{minipage}{.25\hsize}
            \centering
            \includegraphics[width=\hsize]{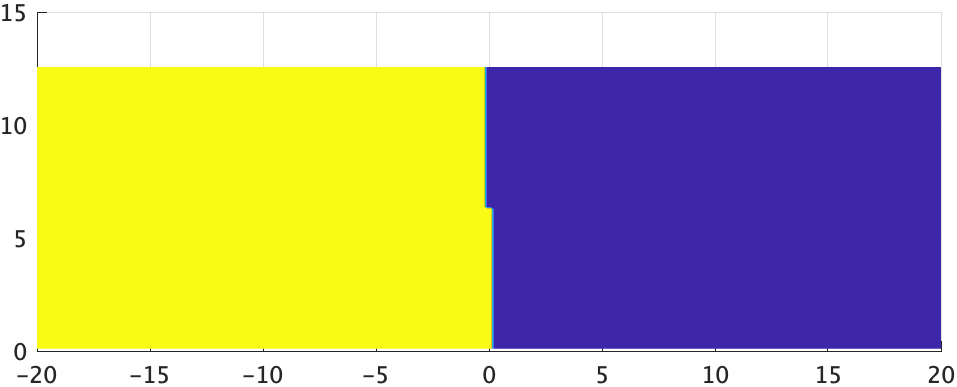}

            $t=0$
        \end{minipage}%
        \begin{minipage}{.25\hsize}
            \centering
            \includegraphics[width=\hsize]{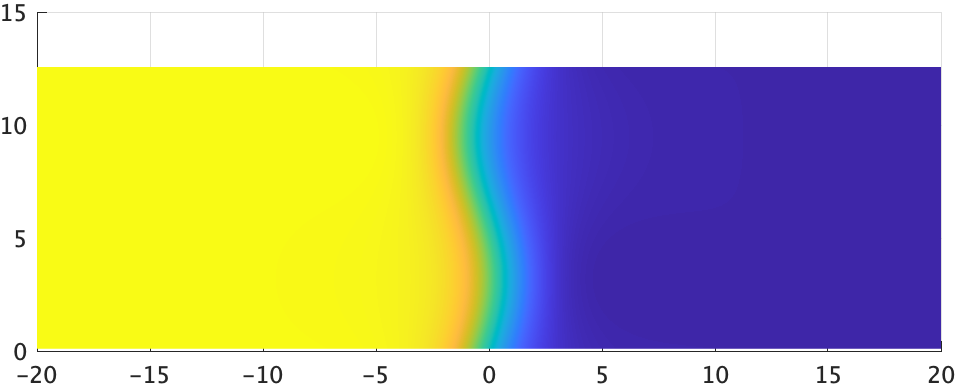}

            $t=100$
        \end{minipage}%
        \begin{minipage}{.25\hsize}
            \centering
            \includegraphics[width=\hsize]{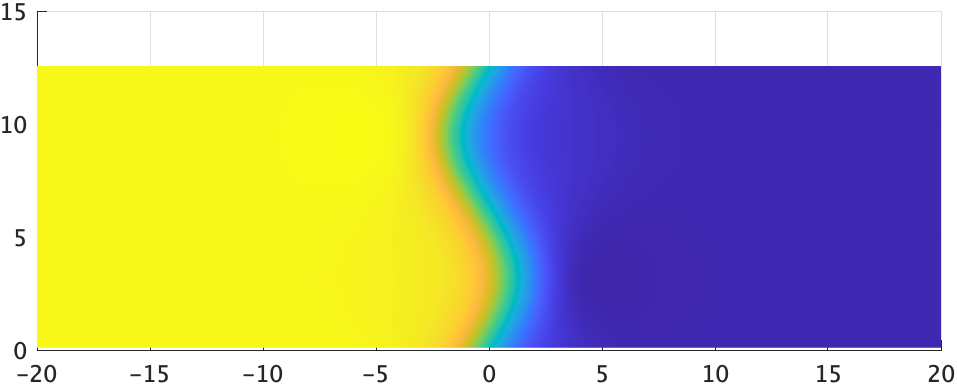}

            $t=500$
        \end{minipage}%
        \begin{minipage}{.25\hsize}
            \centering
            \includegraphics[width=\hsize]{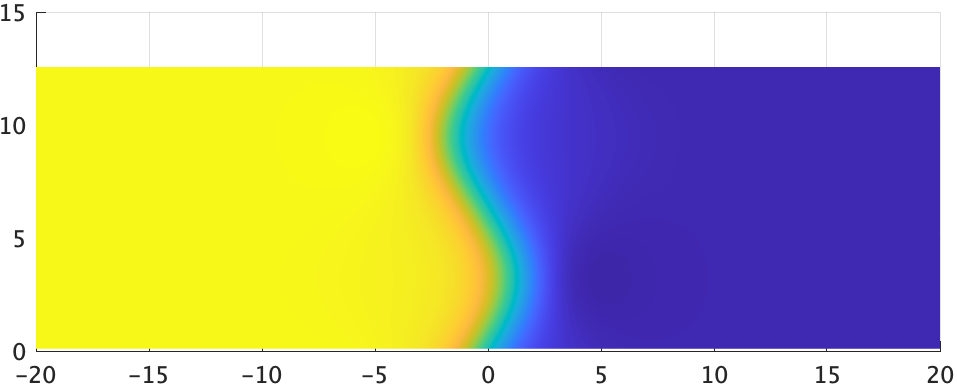}

            $t=1000$
        \end{minipage}%
        \caption{Numerical computations of the bidomain Allen--Cahn equation when $d=2\pi/0.5$, $\theta=\pi/4$, $a=0.9$, $b=0$, and $\alpha=0.4$.}
        \label{fig:pi4_05}
    \end{figure}%
    \begin{figure}[tb]
        \begin{minipage}{.25\hsize}
            \centering
            \includegraphics[width=\hsize]{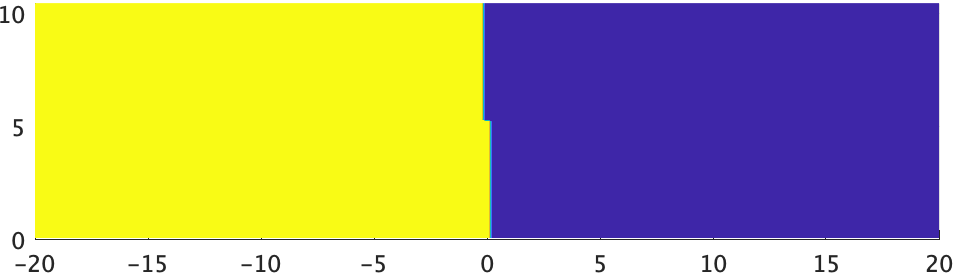}

            $t=0$
        \end{minipage}%
        \begin{minipage}{.25\hsize}
            \centering
            \includegraphics[width=\hsize]{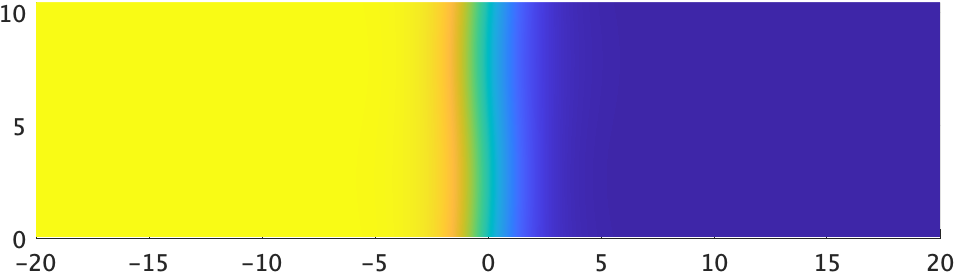}

            $t=100$
        \end{minipage}%
        \begin{minipage}{.25\hsize}
            \centering
            \includegraphics[width=\hsize]{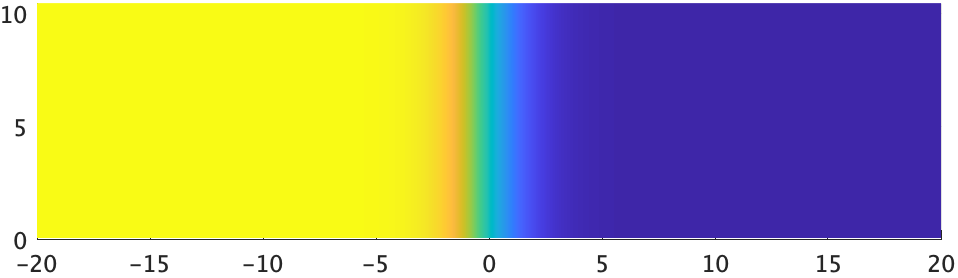}

            $t=500$
        \end{minipage}%
        \begin{minipage}{.25\hsize}
            \centering
            \includegraphics[width=\hsize]{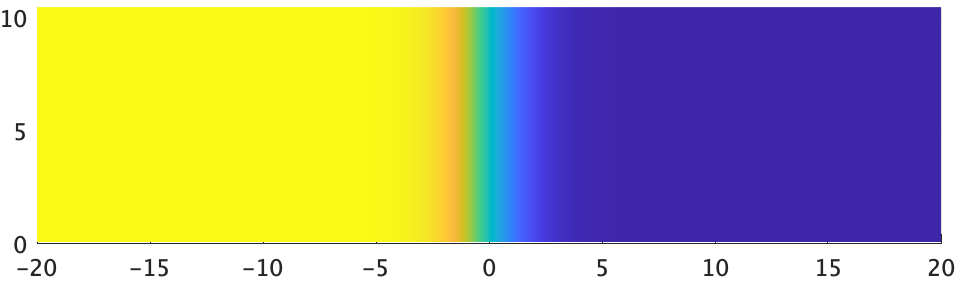}

            $t=1000$
        \end{minipage}%
        \caption{Numerical computations of the bidomain Allen--Cahn equation when $d=2\pi/0.6$, $\theta=\pi/4$, $a=0.9$, $b=0$, and $\alpha=0.4$.}
        \label{fig:pi4_06}
    \end{figure}%
    \begin{figure}[tb]
        \begin{minipage}{.25\hsize}
            \centering
            \includegraphics[width=\hsize]{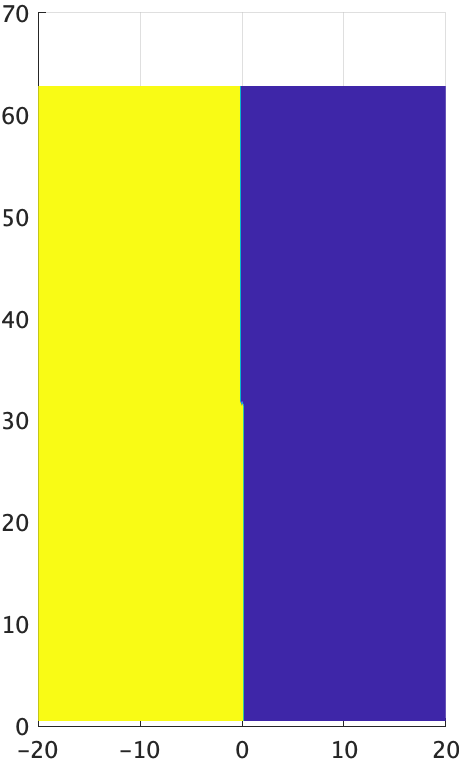}
            
            $t=0$
        \end{minipage}%
        \begin{minipage}{.25\hsize}
            \centering
            \includegraphics[width=\hsize]{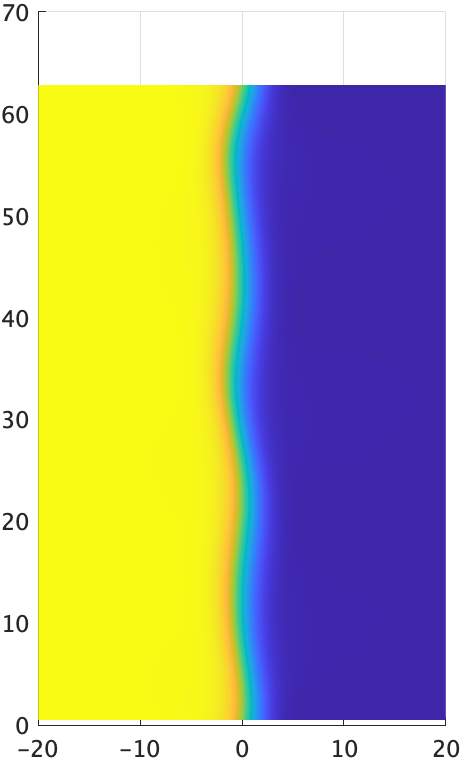}

            $t=200$
        \end{minipage}%
        \begin{minipage}{.25\hsize}
            \centering
            \includegraphics[width=\hsize]{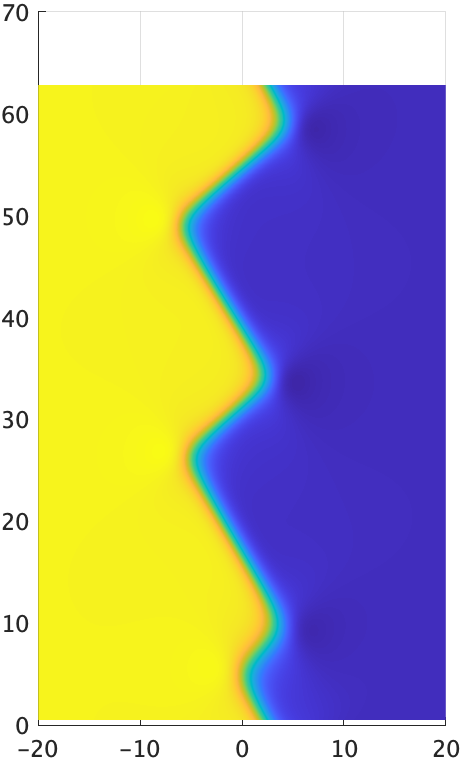}

            $t=600$
        \end{minipage}%
        \begin{minipage}{.25\hsize}
            \centering
            \includegraphics[width=\hsize]{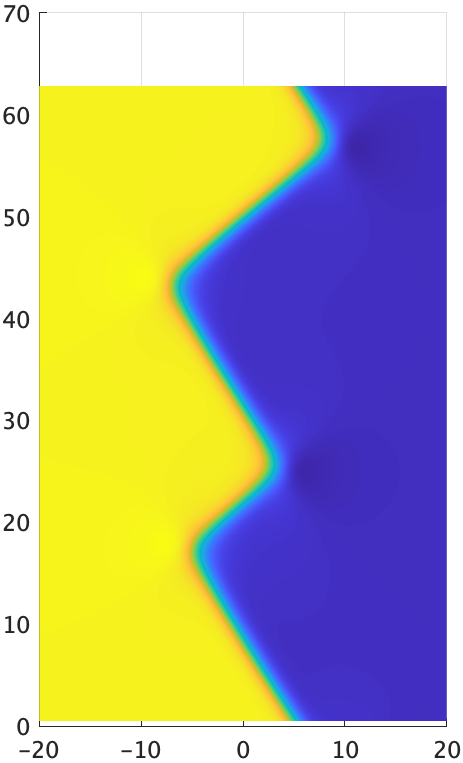}

            $t=1000$
        \end{minipage}%
        \caption{Numerical computations of the bidomain Allen--Cahn equation when $d=2\pi/0.1$, $\theta=\pi/5$, $a=0.9$, $b=0$, and $\alpha=0.4$.}
        \label{fig:pi5_01}
    \end{figure}%
    \begin{figure}[tb]
        \begin{minipage}{.25\hsize}
            \centering
            \includegraphics[width=\hsize]{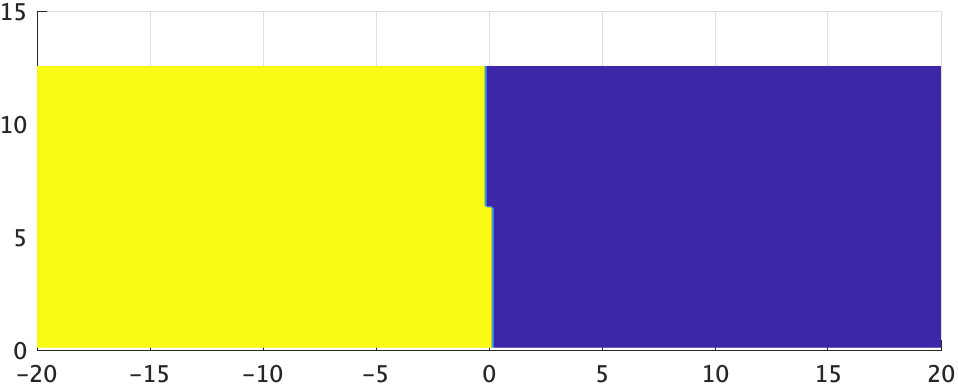}
            
            $t=0$
        \end{minipage}%
        \begin{minipage}{.25\hsize}
            \centering
            \includegraphics[width=\hsize]{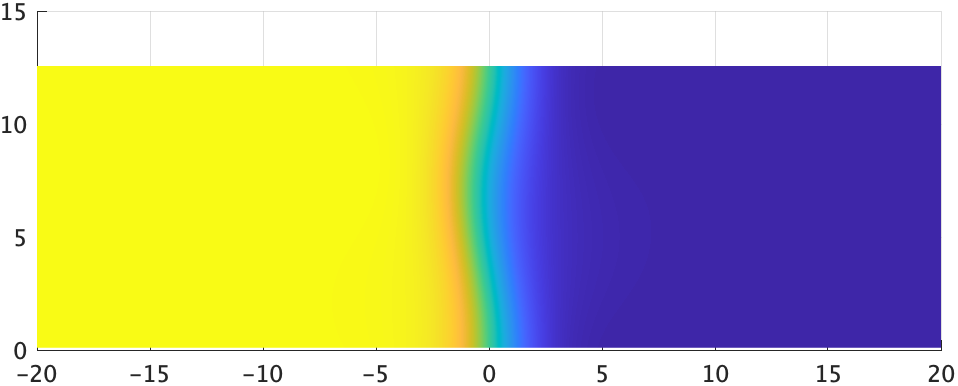}

            $t=200$
        \end{minipage}%
        \begin{minipage}{.25\hsize}
            \centering
            \includegraphics[width=\hsize]{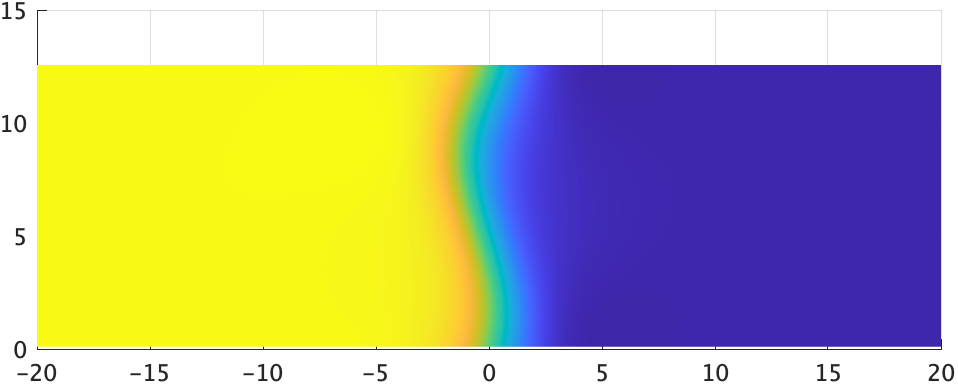}

            $t=600$
        \end{minipage}%
        \begin{minipage}{.25\hsize}
            \centering
            \includegraphics[width=\hsize]{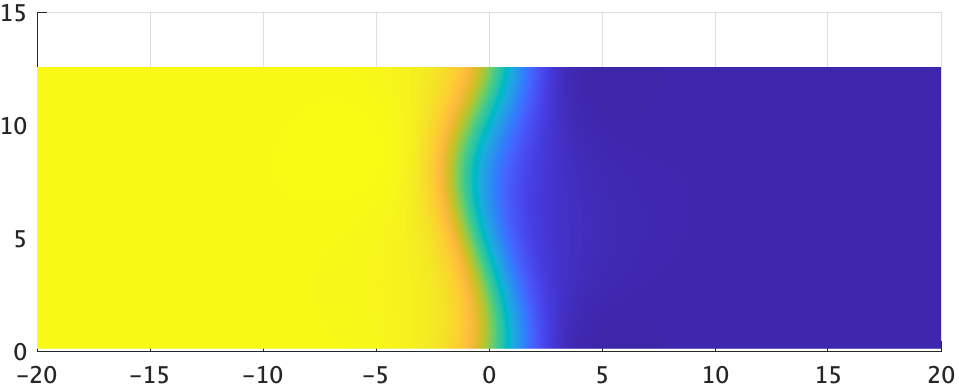}

            $t=1000$
        \end{minipage}%
        \caption{Numerical computations of the bidomain Allen--Cahn equation when $d=2\pi/0.5$, $\theta=\pi/5$, $a=0.9$, $b=0$, and $\alpha=0.4$.}
        \label{fig:pi5_05}
    \end{figure}%
    \begin{figure}[tb]
        \begin{minipage}{.25\hsize}
            \centering
            \includegraphics[width=\hsize]{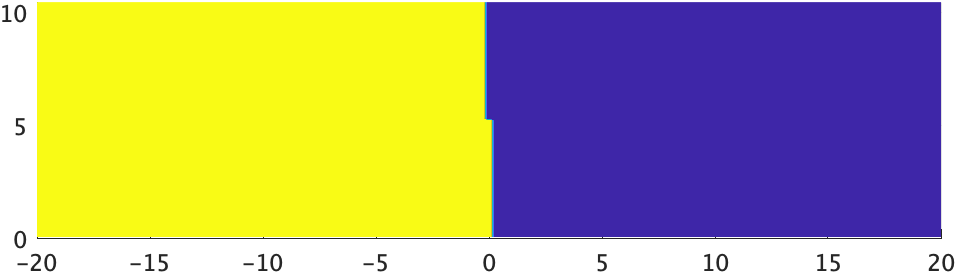}
            
            $t=0$
        \end{minipage}%
        \begin{minipage}{.25\hsize}
            \centering
            \includegraphics[width=\hsize]{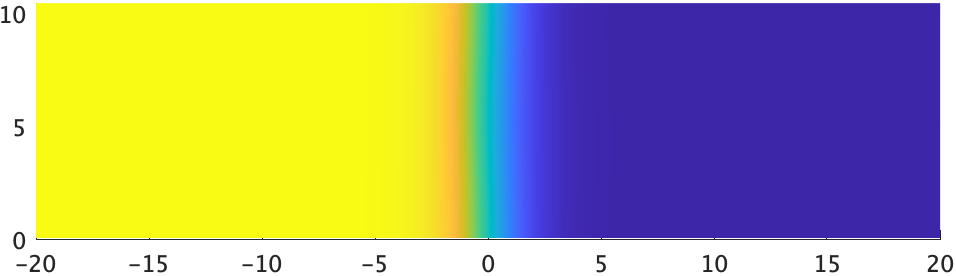}

            $t=200$
        \end{minipage}%
        \begin{minipage}{.25\hsize}
            \centering
            \includegraphics[width=\hsize]{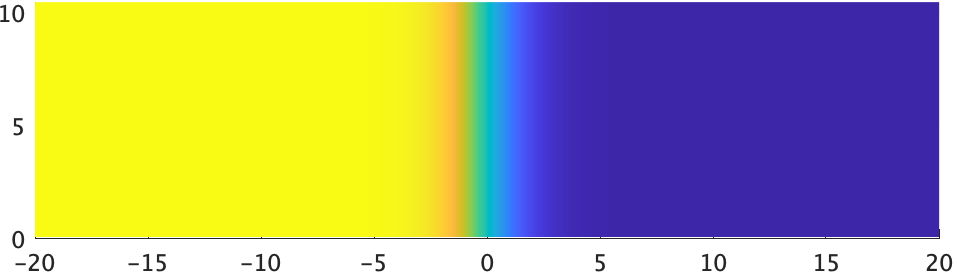}

            $t=600$
        \end{minipage}%
        \begin{minipage}{.25\hsize}
            \centering
            \includegraphics[width=\hsize]{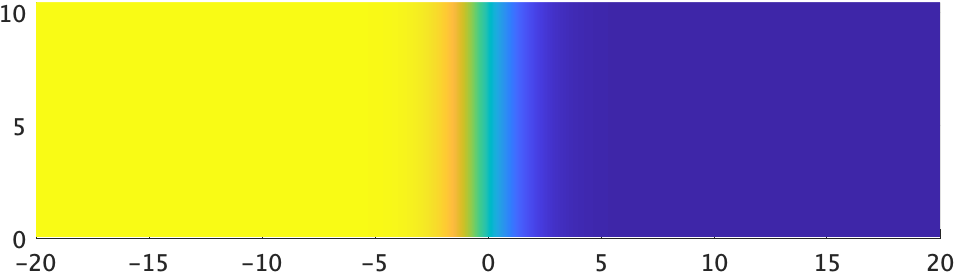}

            $t=1000$
        \end{minipage}%
        \caption{Numerical computations of the bidomain Allen--Cahn equation when $d=2\pi/0.6$, $\theta=\pi/5$, $a=0.9$, $b=0$, and $\alpha=0.4$.}
        \label{fig:pi5_06}
    \end{figure}%

    We study how the principal eigenvalue and the strip region's width affect planar fronts' stability and observe that Hopf bifurcation occurs when the planar front is unstable.
    To this end, we employ the numerical method developed in \cref{subsec:eigenvalue} for computing principal eigenvalues and the one in \cref{subsec:planar_fronts_pulses} for computing the time evolution of a slightly perturbed planar front.

    We now investigate the planar fronts' stability of the bidomain Allen--Cahn equation on the strip region $\mathcal{S}_d$.
    Define the discrete wavenumber as $w_{l,d}\coloneqq2\pi l/d$ ($l\in\mathbb{Z}$, $d>0$) and denote the principal eigenvalue of $\mathcal{P}_{w_{l,d}}$ as $\lambda_{l,d}$.
    First, let us numerically investigate how the sign of the real part of the first mode's principal eigenvalue $\lambda_{1,d}$ changes as $d$ varies.
    As shown in \cref{fig:lambda}, when $2\pi/d$ is small (i.e., $d$ is large), $\Re\lambda_{1,d}$ is positive, and its value is gradually increasing.
    However, the increase eventually stops and starts to decrease, and after a specific value of $d$, it becomes negative.
    This observation indicates that when $d$ is large, the planar front is unstable and zigzag front appears, but as $d$ decreases, all the unstable modes disappear, and the planar front becomes stable.
    We select the parameters as $a=0.9$, $b=0$, and $\alpha=0.4$.
    \Cref{fig:pi4_01,fig:pi4_05,fig:pi4_06} consider the case of $\theta=\pi/4$, and \cref{fig:pi5_01,fig:pi5_05,fig:pi5_06} consider the case of $\theta=\pi/5$.
    We vary the width $d$ of the strip region $\mathcal{S}_d$ as $2\pi/0.1$, $2\pi/0.5$, and $2\pi/0.6$.
    We set the initial values to be completely flat with a slight perturbation added.
    Comparing these numerical results with those in \cref{fig:lambda}, we can say that the above scenario is, to some extent, correct.
    These angles correspond to points where the Frank plot is nonconvex, and \cref{thm:stability} implies that planar fronts in these directions are unstable.
    After destabilization, the planar front becomes a zigzag front, which eventually forms one peak with coarsening. 
    Such a zigzag front does not exist in the monodomain Allen-Cahn equation~\cite[subsection 6.1]{MM16}, a significant characteristic of the bidomain Allen-Cahn equation.
    Moreover, as we also observe in \cref{subsec:coexistence} in detail, the zigzag front's appearance is caused by the Hopf bifurcation.

    \subsection{Zigzag fronts}

    As we observed in the previous subsection, zigzag fronts appear when planar fronts are unstable. 
    In this subsection, we study their asymptotic behavior.

    \subsubsection{Speed of the zigzag front}

    \begin{figure}
        \centering
        \includegraphics[width=.3\hsize]{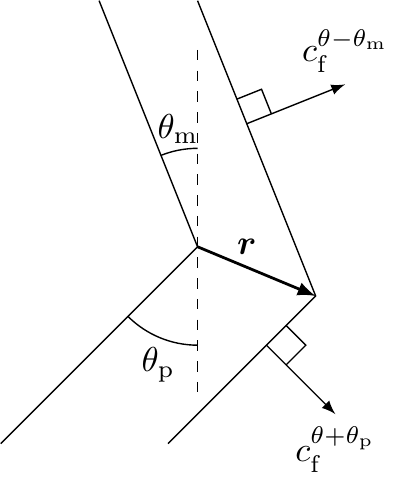}

        \caption{A schematic of the coarsened zigzag front}
        \label{fig:coarsened_zigzag_front}
    \end{figure}%
    \begin{figure}[tb]
        \begin{minipage}{.33\hsize}
            \centering
            \includegraphics[width=\hsize]{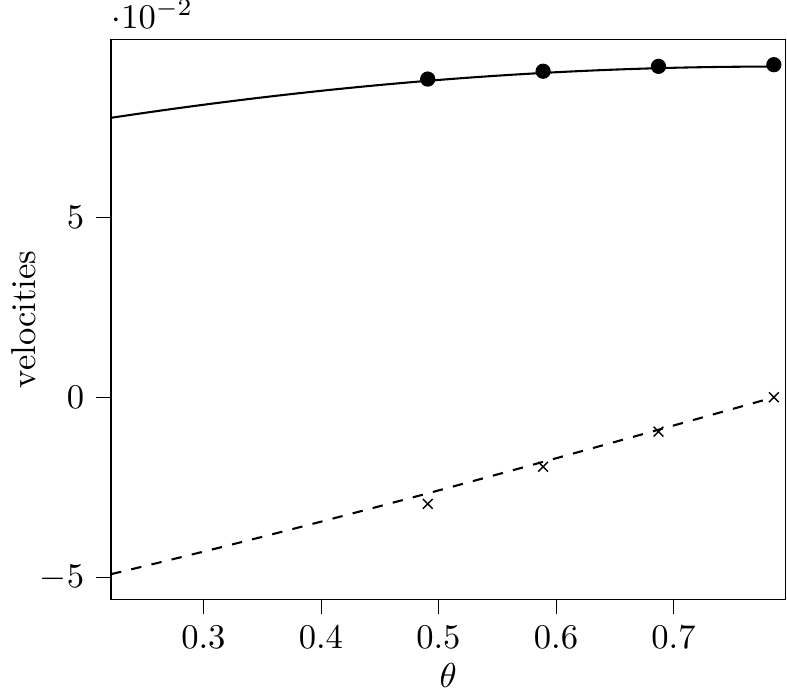}

            (a) $a=0.7$
        \end{minipage}%
        \begin{minipage}{.33\hsize}
            \centering
            \includegraphics[width=\hsize]{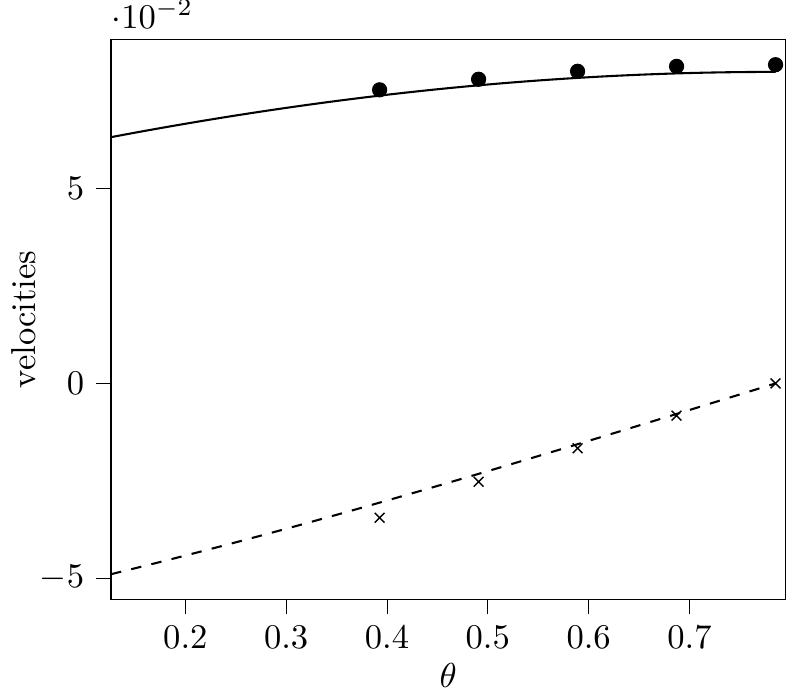}

            (b) $a=0.8$
        \end{minipage}%
        \begin{minipage}{.33\hsize}
            \centering
            \includegraphics[width=\hsize]{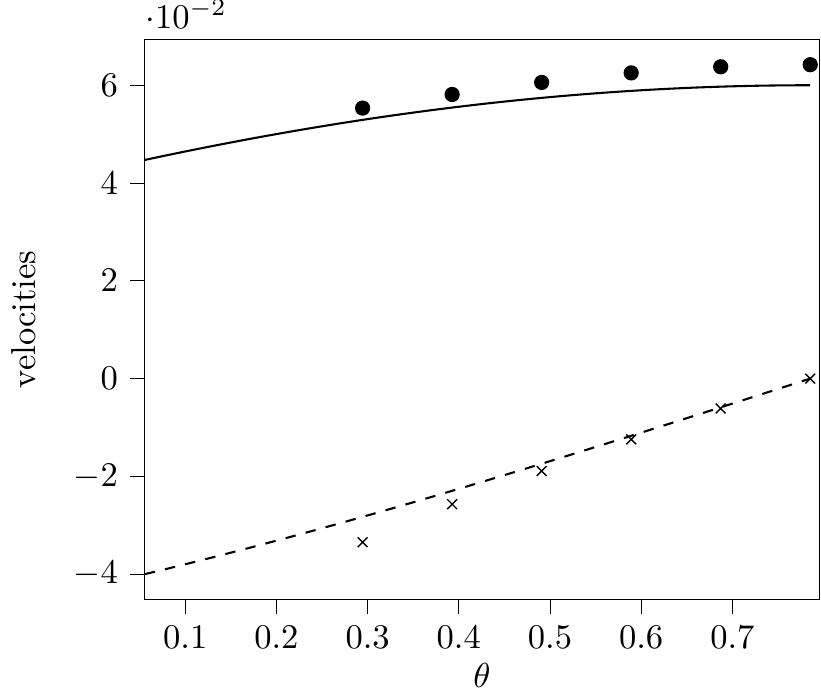}

            (c) $a=0.9$
        \end{minipage}%
        \caption{Comparison of the theoretical speed of the zigzag front with the numerical results.
        Solid lines show the theoretical speeds in $\xi$ direction, broken lines show theoretical speeds in $\eta$ direction, dots show the numerical results of speed in $\xi$ direction, and crosses show the numerical results of speed in $\eta$ direction.}
        \label{fig:speed_comparison}
    \end{figure}

    As we saw in the previous subsection, when the planar front is unstable, it eventually forms a single peak while coarsening. 
    Considering this front in the situation as $d$ tends to infinity, we can compute the zigzag front speed analytically by elementary geometric arguments. 
    Therefore, we below perform the calculation to confirm that the numerical experiments agree with them and provide a benchmark for the numerical method's correctness developed in the \cref{subsec:iterative_method}.

    Let us consider the coarsened ideal zigzag front, as depicted in \cref{fig:coarsened_zigzag_front}.
    As time passes by $\Delta t$, the planar front in the direction of $\bm{n}^{\theta+\theta_{\mathrm{p}}}$ travels distance $c_\rmf^{\theta+\theta_{\mathrm{p}}}\Delta t$, and the one in $\bm{n}^{\theta-\theta_{\mathrm{m}}}$ does the distance $c_\rmf^{\theta-\theta_{\mathrm{m}}}\Delta t$.
    From these facts, we can write down the vector $\bm{r}$ in \cref{fig:coarsened_zigzag_front} concretely as follows:
    \begin{align*}
        \bm{r}=
        \frac{\Delta t}{\sin(\theta_{\mathrm{m}}+\theta_{\mathrm{p}})}
        \begin{pmatrix}
            c_\rmf^{\theta+\theta_{\mathrm{p}}}\sin\theta_{\mathrm{m}}+c_\rmf^{\theta-\theta_{\mathrm{m}}}\sin\theta_{\mathrm{p}}\\
            -c_\rmf^{\theta+\theta_{\mathrm{p}}}\cos\theta_{\mathrm{m}}+c_\rmf^{\theta-\theta_{\mathrm{m}}}\cos\theta_{\mathrm{p}}
        \end{pmatrix}.
    \end{align*}
    Therefore, by dividing the vector $\bm{r}$ by $\Delta t$ and taking the limit of $\Delta t\to0$, we obtain the velocity vector $\bm{v}$ as
    \begin{align*}
        \bm{v}
        =
        \frac{1}{\sin(\theta_{\mathrm{m}}+\theta_{\mathrm{p}})}
        \begin{pmatrix}
            c_\rmf^{\theta+\theta_{\mathrm{p}}}\sin\theta_{\mathrm{m}}+c_\rmf^{\theta-\theta_{\mathrm{m}}}\sin\theta_{\mathrm{p}}\\
            -c_\rmf^{\theta+\theta_{\mathrm{p}}}\cos\theta_{\mathrm{m}}+c_\rmf^{\theta-\theta_{\mathrm{m}}}\cos\theta_{\mathrm{p}}
        \end{pmatrix}.
    \end{align*}
    If $b=0$, then, as will be confirmed in \cref{subsubsec:zigzag_Frank}, $\theta+\theta_{\mathrm{p}}$ and $\theta-\theta_{\mathrm{m}}$ approach the angles of contact between the Frank diagram and its convex hull, where the planar front's speed is equal:
    \begin{align*}
        \tilde{c}
        \coloneqq
        c_\rmf^{\theta+\theta_{\mathrm{p}}}
        =
        c_\rmf^{\theta-\theta_{\mathrm{m}}}.
    \end{align*}
    Hence, in this case, the velocity vector $\bm{v}$ is rewritten in a more straightforward form:
    \begin{align}
        \bm{v}
        =
        \frac{\tilde{c}}{\sin(\theta_{\mathrm{m}}+\theta_{\mathrm{p}})}
        \begin{pmatrix}
            \sin\theta_{\mathrm{m}}+\sin\theta_{\mathrm{p}}\\
            -\cos\theta_{\mathrm{m}}+\cos\theta_{\mathrm{p}}
        \end{pmatrix}.
        \label{eq:front_velocity}
    \end{align}
    The first element corresponds to the speed in $\xi$ direction, and the second one to the speed in $\eta$ direction.

    \Cref{fig:speed_comparison} compares the theoretical value of the front velocity of equation \eqref{eq:front_velocity} with the one calculated from the numerical results for $a=0.7$, $a=0.8$, and $a=0.9$.
    We set the width $d$ of the strip region $\mathcal{S}_d$ to $100$.
    When $a$ is not so far from $1/2$ (e.g., $a=0.7$, $0.8$), the numerical results agree well with the theoretical values.
    When $a=0.9$, some disparities exist, but the numerical results well understand the fronts' qualitative properties.

    \subsubsection{Correspondence between angles of zigzag fronts and the Frank diagram}
    \label{subsubsec:zigzag_Frank}

    \begin{figure}[tb]
        \centering
        \includegraphics[width=.6\hsize]{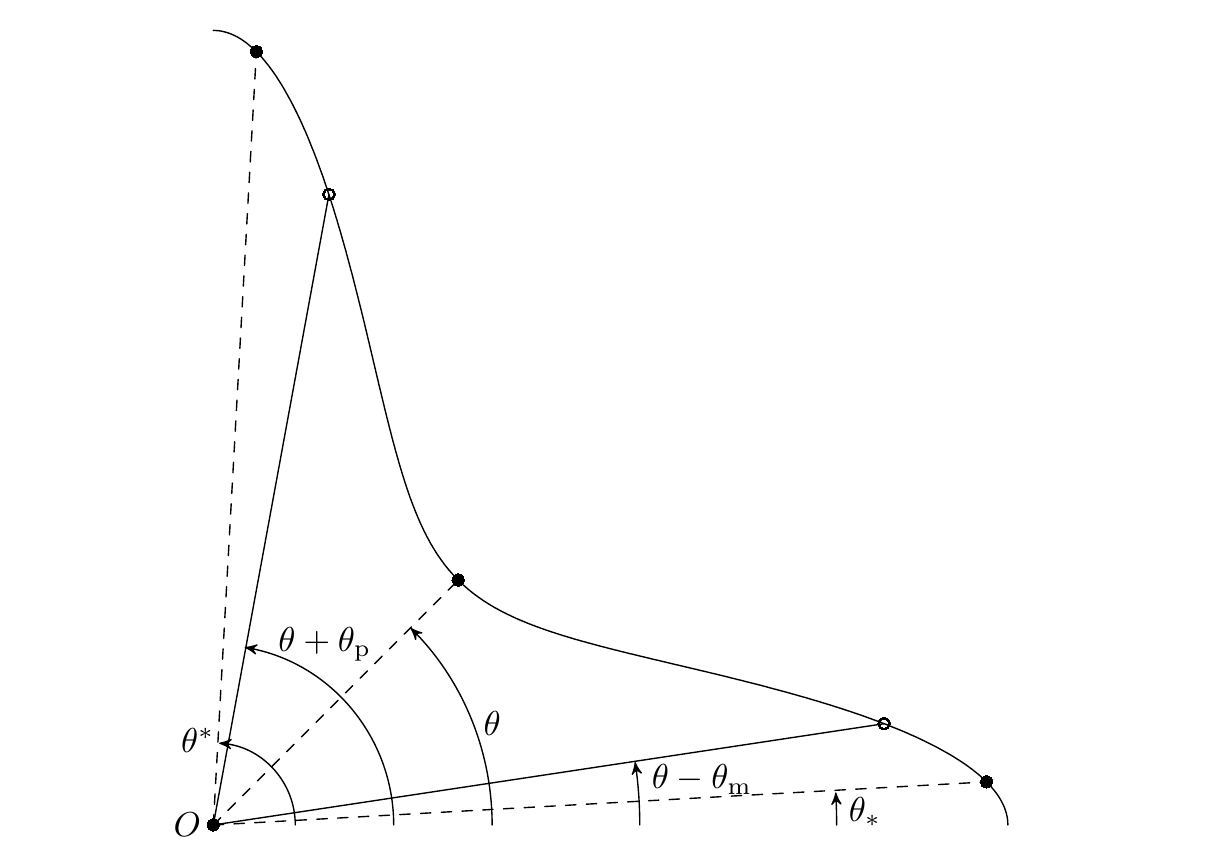}
        \caption{Correspondence of angles in the Frank plot.}
        \label{fig:frank_part}
    \end{figure}%
    \begin{figure}[tb]
        \centering
        \includegraphics[width=.7\hsize]{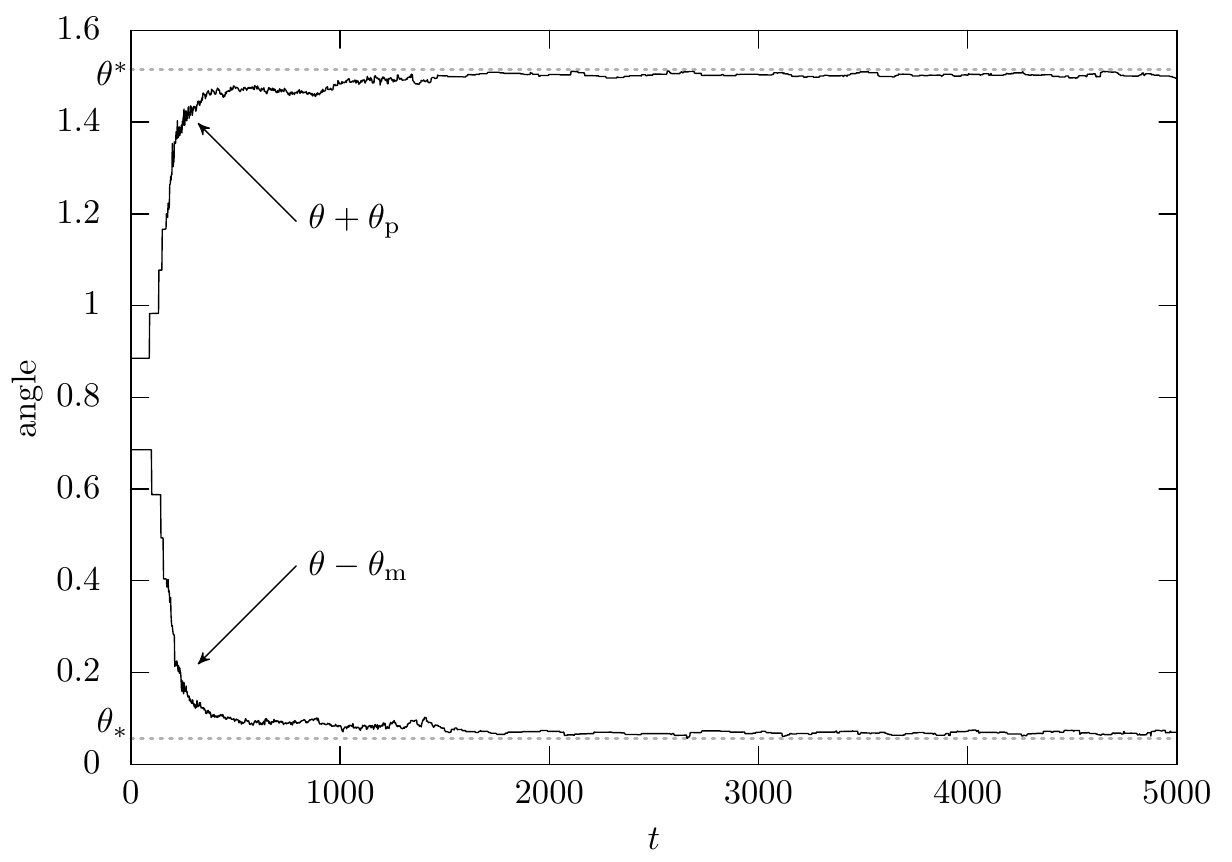}
        \caption{Graphs of $\theta-\theta_{\mathrm{m}}$ and $\theta+\theta_{\mathrm{p}}$.}
        \label{fig:asymptotic_angle}
    \end{figure}

    In \cref{subsec:stability_front}, we investigated a relationship between the principal eigenvalue's real part and the planar front's instability.
    In this section, we focus on the zigzag front itself.
    If we look at the figure for $t=1000$ in \cref{fig:pi4_01}, we can observe several corners (peaks) on the zigzag interface.
    In particular, looking at the time evolution, there is some law at work in the mechanism of corners.
    In this section, we study the asymptotic behavior of zigzag fronts in terms of the peak's angles and the Frank diagram.

    As we checked in \cref{subsec:stability_front}, the situation of \cref{thm:stability} is consistent with \cref{fig:pi4_01,fig:pi4_05,fig:pi4_06,fig:pi5_01,fig:pi5_05,fig:pi5_06}, where $\theta=\pi/4$ and $\theta=\pi/5$ correspond to places where the Frank plot is nonconvex (see also \cref{fig:Frank}).
    Planar fronts are indeed unstable, and zigzag fronts appear, but the zigzag fronts themselves appear to have some stability; that is, they propagate over a long period while maintaining their shape.
    Combining this observation with \cref{thm:stability}, we reach the following conjecture.

    \begin{conj}
        The angles of the peak of the zigzag front caused by destabilization asymptotically approach the angles of the contact points between the Frank diagram and its convex hull.
    \end{conj}

    Let us verify this conjecture numerically.
    Define angles $\theta_{\mathrm{m}}$ and $\theta_{\mathrm{p}}$ as in \cref{fig:coarsened_zigzag_front}.
    Namely, the angles of the zigzag front forming the peak are $\theta-\theta_{\mathrm{m}}$ and $\theta+\theta_{\mathrm{p}}$.
    Focus on \cref{fig:frank_part} and assume that the direction of the zigzag front is $\theta$ and that the point $P_\theta=(\cos\theta,\sin\theta)^\top/K(\theta)$ on the Frank plot corresponding to the angle $\theta$ is in the region where the Frank plot is nonconvex.
    Let us denote the angles of the two closest contacts from point $P_\theta$ on the Frank plot as $\theta_*$, $\theta^*$ from the smaller one.
    Then, we expect that
    \begin{align*}
        \theta-\theta_{\mathrm{m}}(t)\rightarrow\theta_*,
        \quad
        \theta+\theta_{\mathrm{p}}(t)\rightarrow\theta^*
    \end{align*}
    hold as $t\to\infty$.
    \Cref{fig:asymptotic_angle} shows the results of one numerical experiment, which verifies that our predictions are correct.

    \subsection{Coexistence of zigzag and planar fronts}
    \label{subsec:coexistence}

    \begin{figure}[tb]
        \centering
        \includegraphics[width=\hsize]{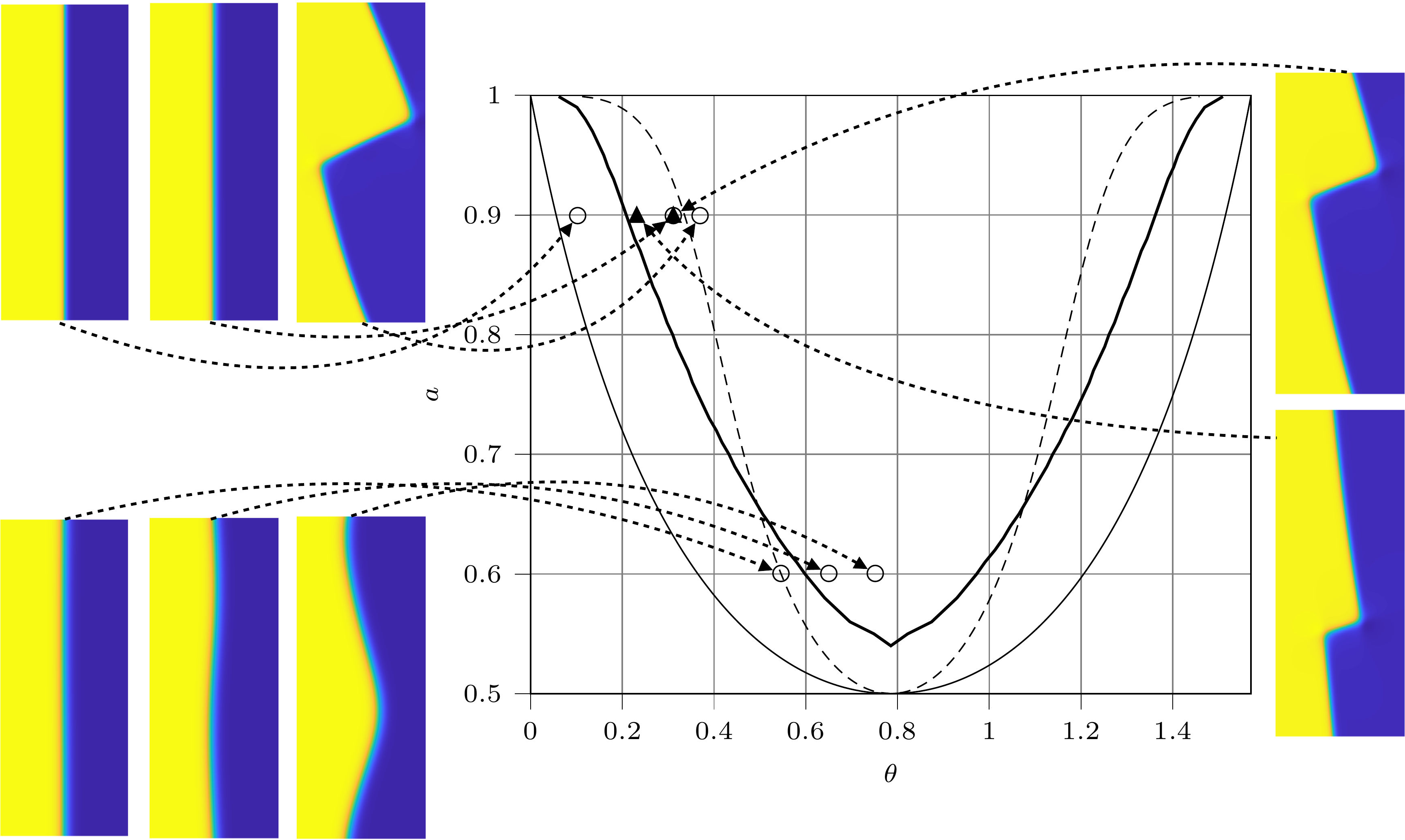}

        \caption{Phase diagram concerning the existence of the zigzag front and the planar front.
        We choose parameters as $b=0$, $d=100$, and $\alpha=0.4$.}
        \label{fig:phase_diagram}
    \end{figure}

    In the monodomain Allen--Cahn equation, the planar front exists in all directions, while the zigzag front in the bidomain Allen--Cahn equation may degenerate into the planar front, as we can see in \cref{fig:pi4_06}, for instance.
    Moreover, the planar front's stability is described by the sign of the Frank plot's curvature as in \cref{thm:stability}.
    However, as shown in \cref{fig:Frank}, there is a portion inside the convex hull of the Frank diagram such that the curvature of the Frank plot is positive.
    We can then conjecture that a coexistence region of the zigzag front and the planar front may exist.
    We investigate this conjecture numerically in this section.

    We search the boundary of existence/nonexistence by varying two parameters $a$ and $\theta$ using the numerical method developed in \cref{subsec:iterative_method}.
    \Cref{fig:phase_diagram} depicts the result.
    The thin solid line represents the angles of contact between the Frank diagram and its convex hull; that is, points above this line are inside the convex hull.
    The broken line expresses the position where the Frank plot's curvature is equal to $0$; above this line, the curvature is negative, while it is positive below this line.
    The thick solid lines represent the boundary of existence/nonexistence of the zigzag front.
    More precisely, the zigzag front exists above this line.

    Furthermore, we can observe bifurcation phenomena.
    For a ``large'' $a$, the subcritical Hopf bifurcation occurs as $\theta$ being a bifurcation parameter.
    Indeed, when numerical computations are performed from the initial values shown in the leftmost figure of \cref{fig:pi4_01} at $a=0.9$, the zigzag front does not occur in the region to the left of the broken line due to the stability of the planar front, but slightly beyond the broken line to the right, the planar front becomes unstable and a zigzag front with a large peak occurs.
    At $a=0.9$, the thick line passes to the left of the dashed line, suggesting that planar fronts and zigzag fronts coexist between these lines.
    In fact, if we use the zigzag front as the initial value and gradually decreases the parameter $\theta$ , the zigzag front does not disappear even if the value of $\theta$ crosses the broken line to the left, but will continue to exist until it reaches the bold line.
    On the other hand, for a ``small'' $a$, the bold line passes to the right of the broken line, suggesting that a supercritical Hopf bifurcation occurs.
    Indeed, the numerical results for $a=0.6$ do not show the same hysteresis as observed for $a=0.9$.
    
    Taking these facts together, it can be inferred that the type of Hopf bifurcation is switched near the intersection of the bold and broken lines, and a double Hopf bifurcation, which is a bifurcation of codimension $2$, is expected to occur.

    \subsection{Coarsening}

    \begin{figure}[tb]
        \centering
        \includegraphics[width=.8\hsize]{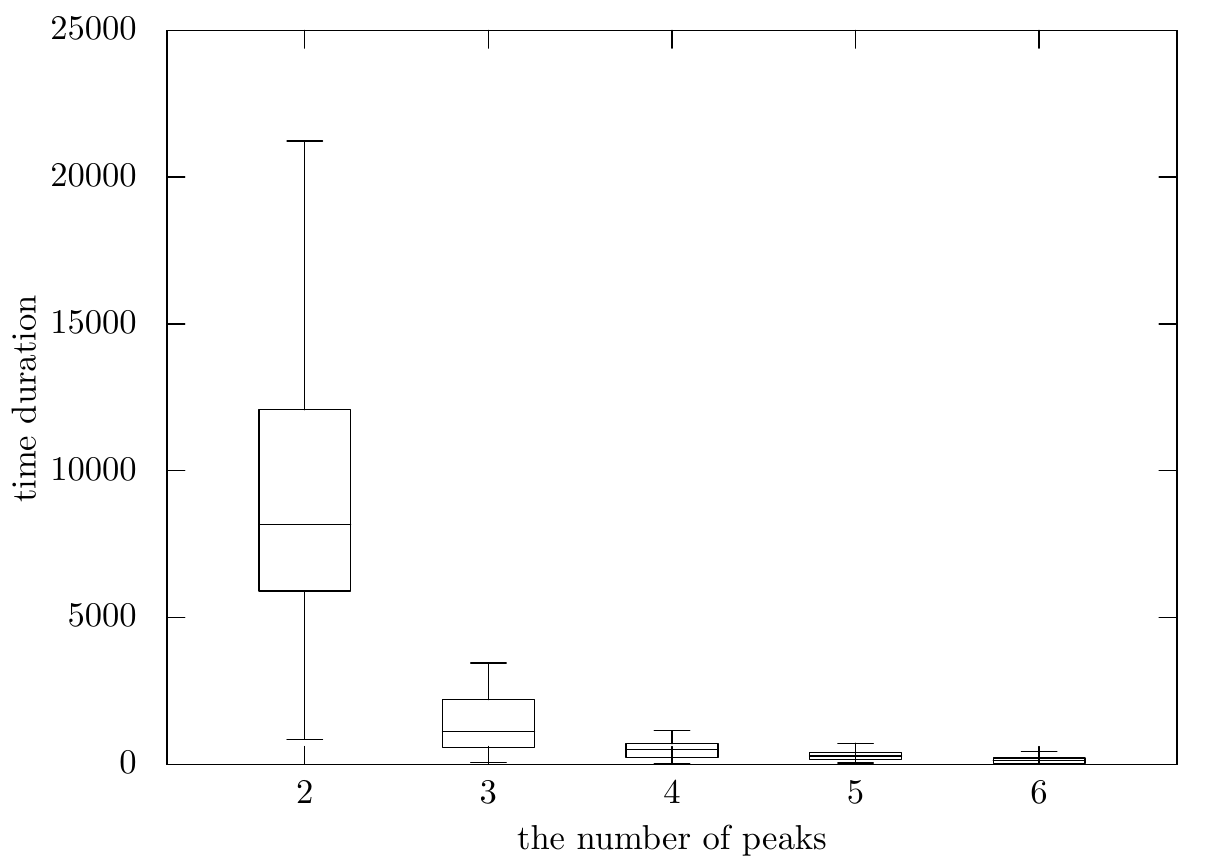}

        \caption{The box plot on the time evolution of the number of peaks in the bidomain Allen--Cahn equation.}
        \label{fig:box_plot}
    \end{figure}

    We, so far, investigated the stability, the asymptotic behavior, and the bifurcation phenomena.
    In this section, we briefly deal with coarsening phenomena for zigzag fronts.

    For example, in \cref{fig:pi4_01}, the zigzag front that appears in the bidomain Allen--Cahn equation undergoes coarsening with time evolution; that is, the number of peaks gradually decreases, common in all cases.
    Therefore, we present the statistical data of the time evolution of the number of peaks.
    \Cref{fig:box_plot} shows the box plot on the time evolution of the number of peaks in the bidomain Allen--Cahn equation, where we choose parameters so that $d=100$, $\theta=\pi/4$, $a=0.9$, $b=0$, and $\alpha=0.4$.
    We made this box plot by performing several numerical experiments of the bidomain Allen--Cahn equation with randomly perturbed initial data.

    We can interpret the result of this numerical experiment as follows.
    First, near the initial time, many peaks are generated from the random initial data, but they coarsen one after another in a relatively short time.
    After that, the number of peaks decreases to two, but the state with two peaks maintains its shape for a relatively long time; that is, the two-peak state retains stability in a sense.

    \subsection{Spreading fronts}

    \begin{figure}[tb]
        \begin{minipage}{.4\hsize}
            \centering
            \includegraphics[width=\hsize]{frank-crop.pdf}
        \end{minipage}
        \qquad
        \begin{minipage}{.4\hsize}
            \centering
            \includegraphics[width=\hsize]{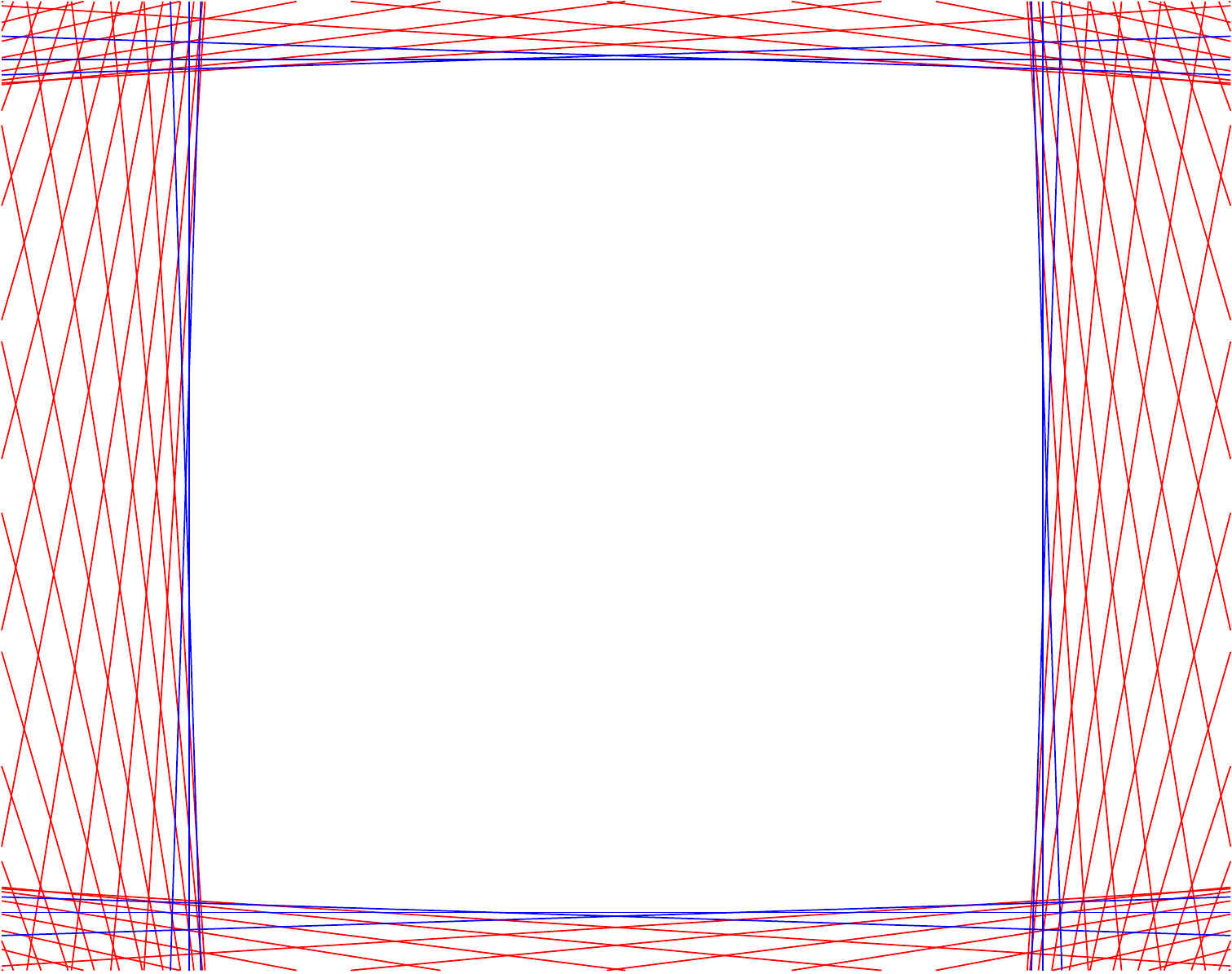}
        \end{minipage}
        \caption{The Frank diagram (left) and the Wulff shape (right) when $a=0.9$ and $b=0$.}
        \label{fig:Frank_Wulff}
    \end{figure}%
    \begin{figure}[tb]
        \begin{minipage}{.33\hsize}
            \centering
            \includegraphics[width=\hsize]{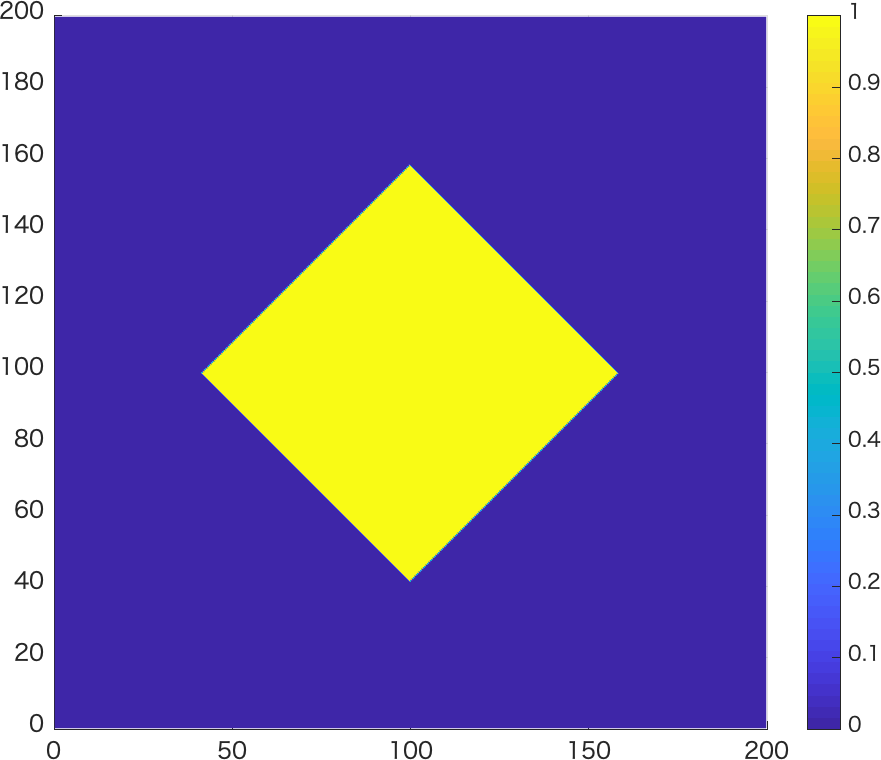}

            $t=0$
        \end{minipage}%
        \begin{minipage}{.33\hsize}
            \centering
            \includegraphics[width=\hsize]{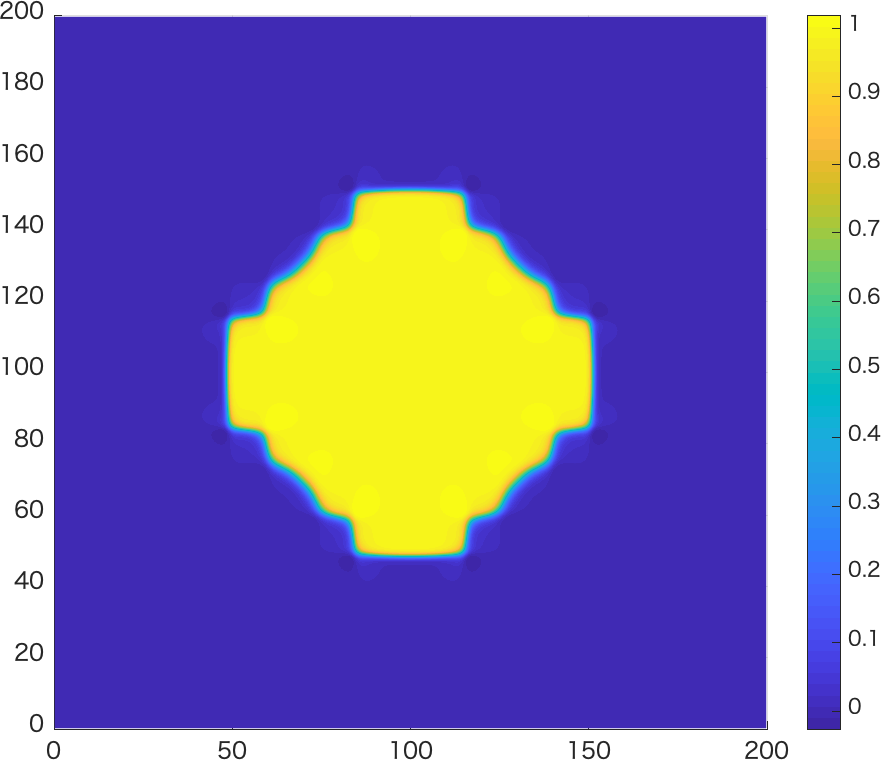}

            $t=200$
        \end{minipage}%
        \begin{minipage}{.33\hsize}
            \centering
            \includegraphics[width=\hsize]{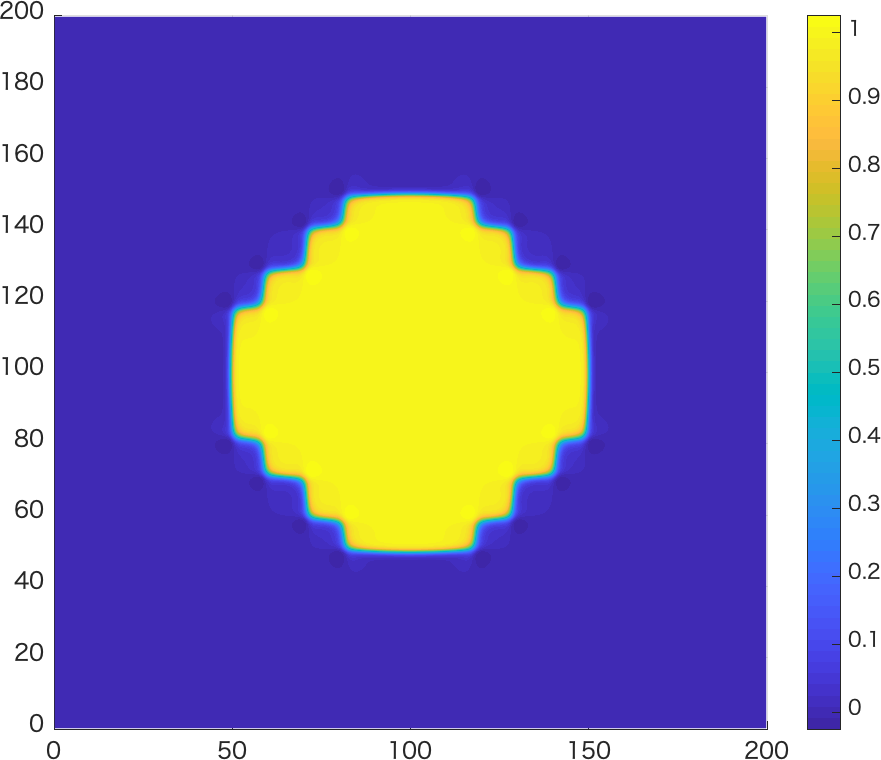}

            $t=400$
        \end{minipage}%

        \begin{minipage}{.33\hsize}
            \centering
            \includegraphics[width=\hsize]{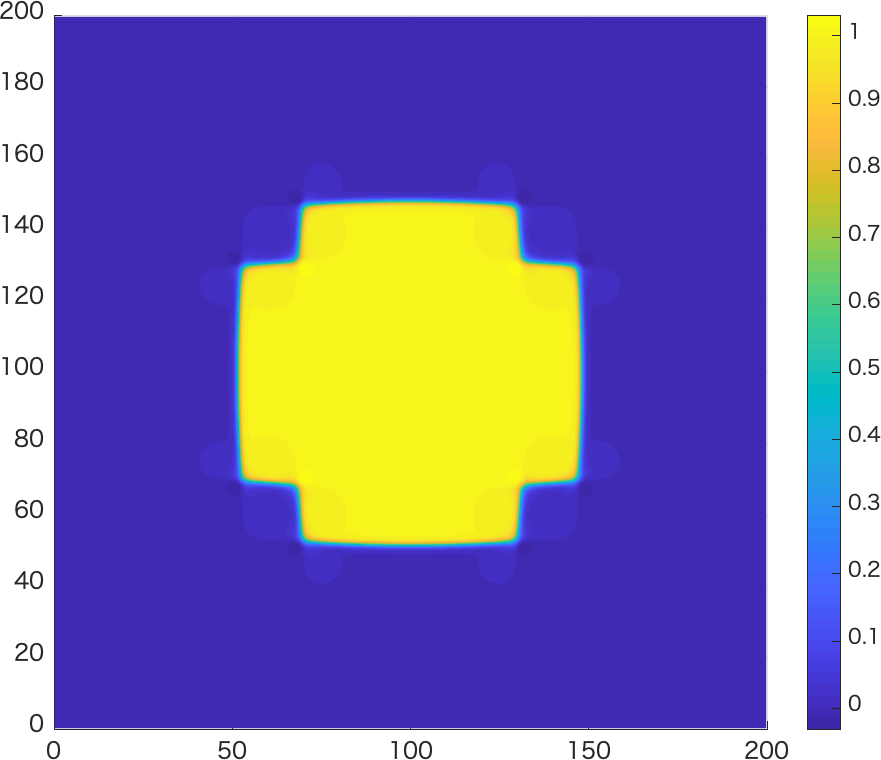}

            $t=800$
        \end{minipage}%
        \begin{minipage}{.33\hsize}
            \centering
            \includegraphics[width=\hsize]{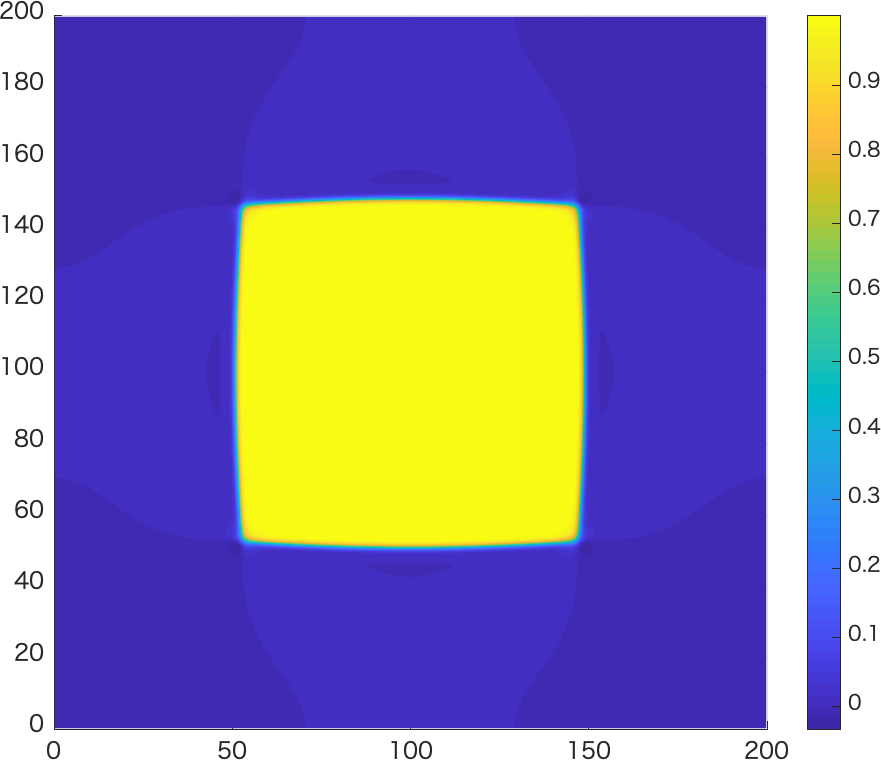}

            $t=2000$
        \end{minipage}%
        \begin{minipage}{.33\hsize}
            \centering
            \includegraphics[width=\hsize]{wulff-crop.pdf}

            the Wulff shape
        \end{minipage}%
        \caption{The spreading front and the Wulff shape when $a=0.9$, $b=0$, and $\alpha=0.49$.}
        \label{fig:spreading_front}
    \end{figure}

    So far, we have investigated the qualitative properties of planar and zigzag fronts in the strip region.
    However, it is also essential to study how the front spreads from the initial value with compact support.
    In particular, it is mathematically interesting to investigate the asymptotic shape of spreading fronts.
    The Wulff shape is a tool to elucidate a part of this.

    The \textit{Wulff shape} $\mathscr{W}$ is defined by
    \begin{align}
        \mathscr{W}
        =
        \bigcap_{0\le\theta<2\pi}\left\{(x,y)^\top\in\mathbb{R}^2\mid x\cos\theta+y\sin\theta\le K(\theta)\right\}.
        \label{eq:Wulff}
    \end{align}
    \Cref{fig:Frank_Wulff} shows the Wulff shape and the corresponding Frank diagram when $a=0.9$ and $b=0$.
    There is a kind of duality between the Wulff shape and the Frank diagram, which we below show their relation roughly.

    Denote the Frank diagram's convex hull by $\hat{\mathscr{F}}$ and denote a point on the Frank plot $\mathscr{F}$ by $P_\theta=(\cos\theta,\sin\theta)^\top/K(\theta)$.
    Moreover, we define a set $\mathcal{S}$ by
    \begin{align*}
        \mathcal{S}
        \coloneqq
        \{\theta\in[0,2\pi)\mid P_\theta\in\partial\hat{\mathscr{F}}\}.
    \end{align*}
    Then, we can prove that the Wulff shape also has the following expression:
    \begin{align*}
        \mathscr{W}
        =
        \bigcap_{\theta\in\mathcal{S}}\left\{(x,y)^\top\in\mathbb{R}^2\mid x\cos\theta+y\sin\theta\le K(\theta)\right\}.
    \end{align*}
    Based on the above expression of the Wulff shape, Mori and Matano proposed the following conjecture in \cite{MM16}.

    \begin{conj}
        Suppose that the normalized speed $c_*$ is positive, and that the planar front in the direction of $\bm{n}^\theta$ is spectrally stable for all $\theta\in\mathcal{S}$.
        Then, the asymptotic shape of the spreading front of the bidomain Allen--Cahn equation is described by the Wulff shape \eqref{eq:Wulff}.
    \end{conj}

    In the anisotropic Allen--Cahn equation, it is well-known that the Wulff shape gives the asymptotic shape of the spreading front \cite{EPS96,ES96,BCP97,AGHMS10}.
    In this sense, the above prediction is not so far off the mark.
    Furthermore, in \cite{MM16}, what can happen is predicted if the above conjecture's stability condition is not satisfied.
    A probable scenario is the following.
    The spreading front approaches the Wulff shape at large length scales.
    However, at finer length scales, the front is modulated by oscillatory waves of the small amplitude corresponding to instabilities in the intermediate wavelengths.

    Let us show one result of numerical experiments.
    \Cref{fig:spreading_front} shows the numerical results of the spreading front for $a=0.9$, $b=0$, and $\alpha=0.49$, and the corresponding Wulff shape.
    As time goes by, we can observe that the spreading front gradually converges to the Wulff shape while coarsening and oscillation, which agrees with the above scenario.
    In this numerical computation, we adopt the numerical scheme developed in \cref{subsec:spreading} and solved the bidomain Allen--Cahn equation in a rectangular region with periodic boundary conditions in both $\xi$ and $\eta$ directions.
    
    \section{Bidomain FitzHugh--Nagumo equation}
    \label{sec:FHN}

    \begin{figure}[tb]
        \begin{minipage}{.25\hsize}
            \centering
            \includegraphics[width=\hsize]{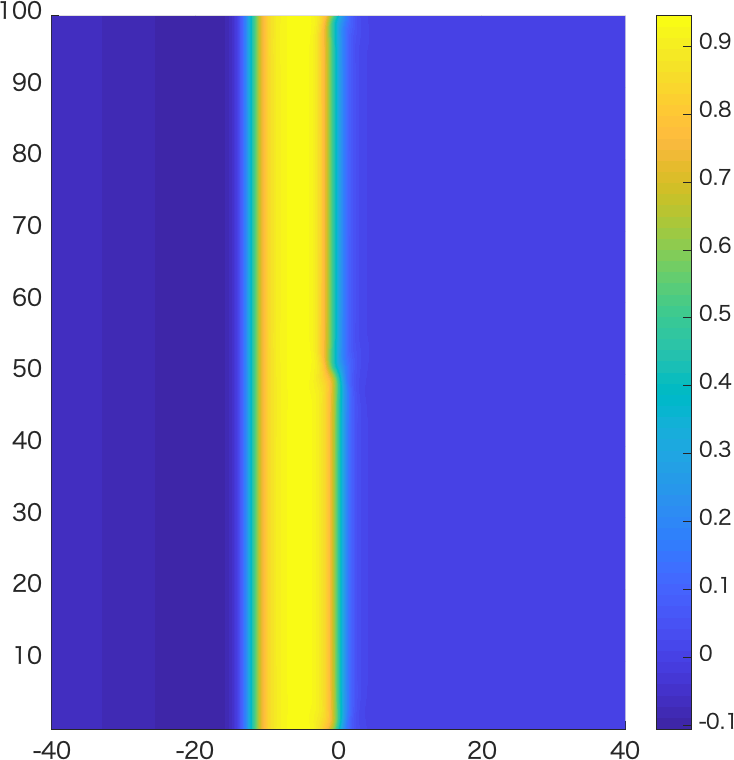}

            $t=0$
        \end{minipage}%
        \begin{minipage}{.25\hsize}
            \centering

            \includegraphics[width=\hsize]{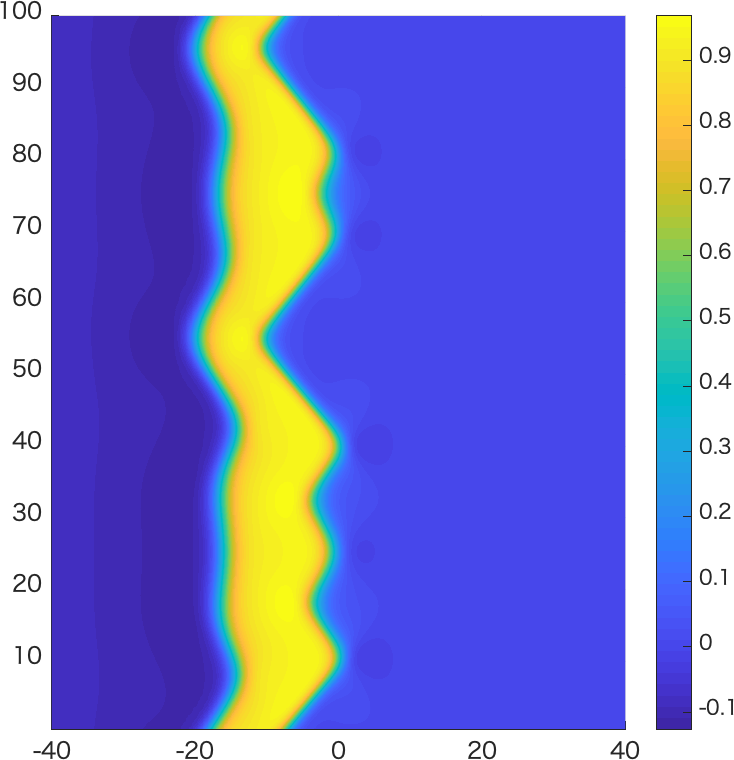}
            $t=200$
        \end{minipage}%
        \begin{minipage}{.25\hsize}
            \centering
            \includegraphics[width=\hsize]{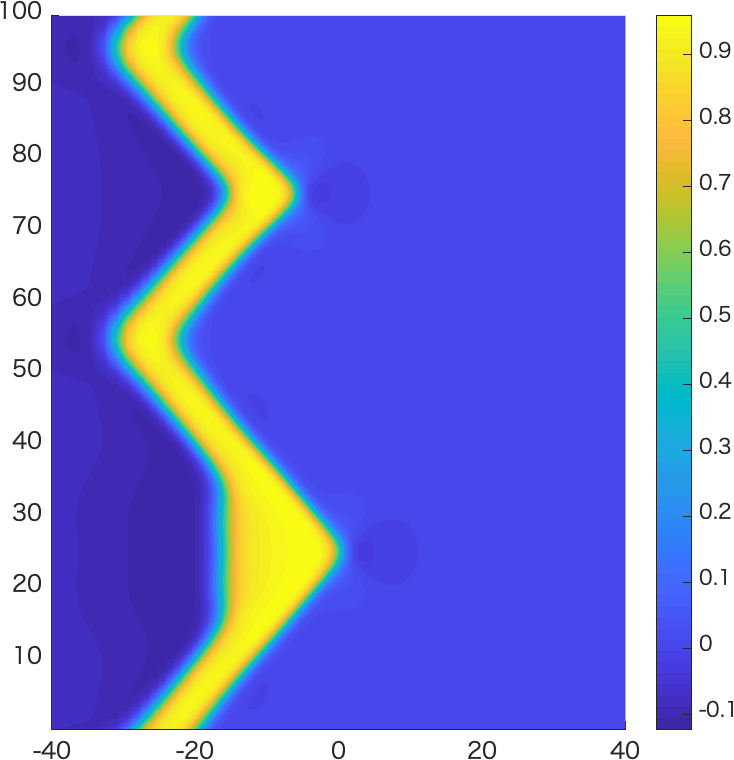}

            $t=400$
        \end{minipage}%
        \begin{minipage}{.25\hsize}
            \centering
            \includegraphics[width=\hsize]{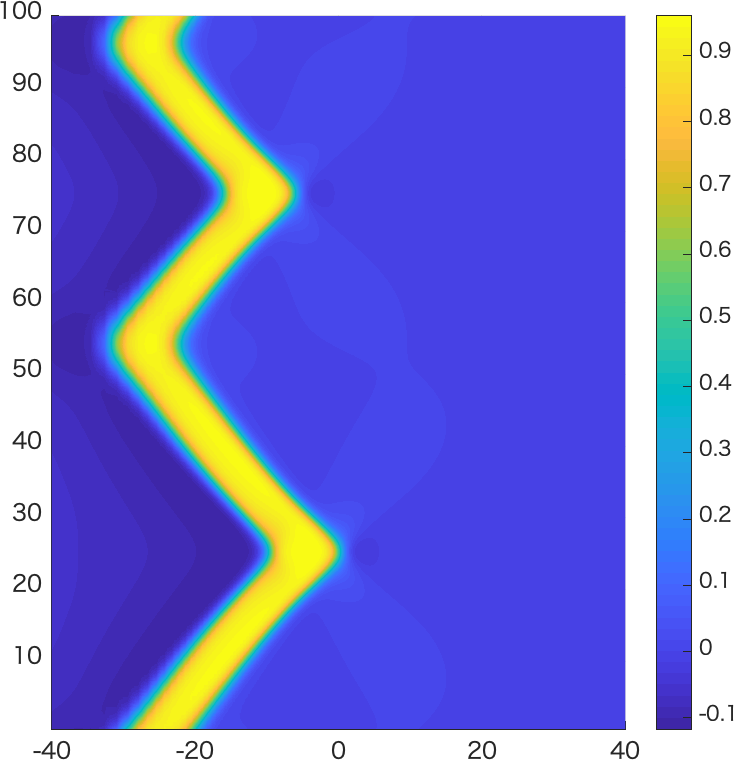}

            $t=600$
        \end{minipage}%

        \begin{minipage}{.25\hsize}
            \centering
            \includegraphics[width=\hsize]{alpha033_00001.png}

            $t=0$
        \end{minipage}%
        \begin{minipage}{.25\hsize}
            \centering
            \includegraphics[width=\hsize]{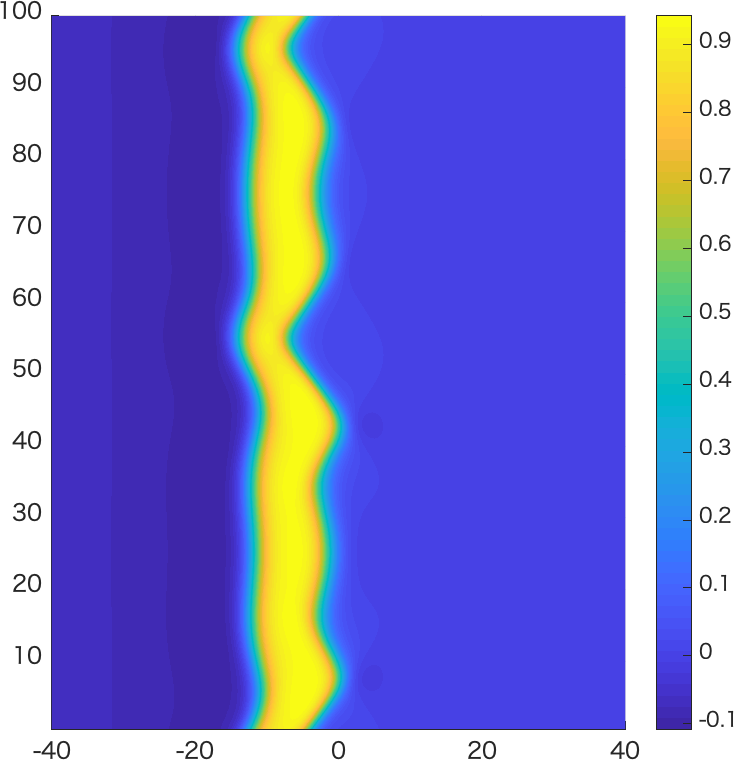}

            $t=120$
        \end{minipage}%
        \begin{minipage}{.25\hsize}
            \centering
            \includegraphics[width=\hsize]{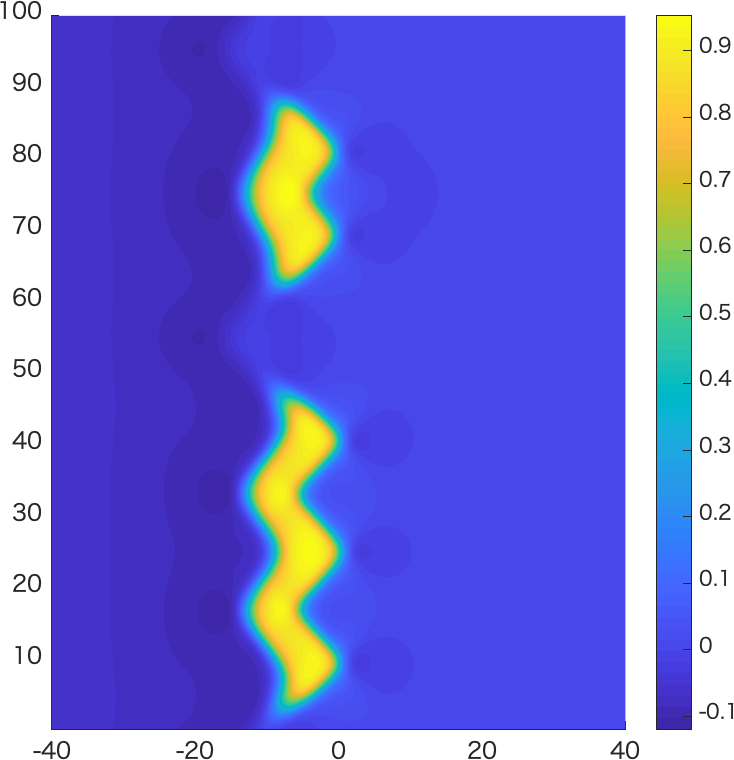}
            
            $t=200$
        \end{minipage}%
        \begin{minipage}{.25\hsize}
            \centering
            \includegraphics[width=\hsize]{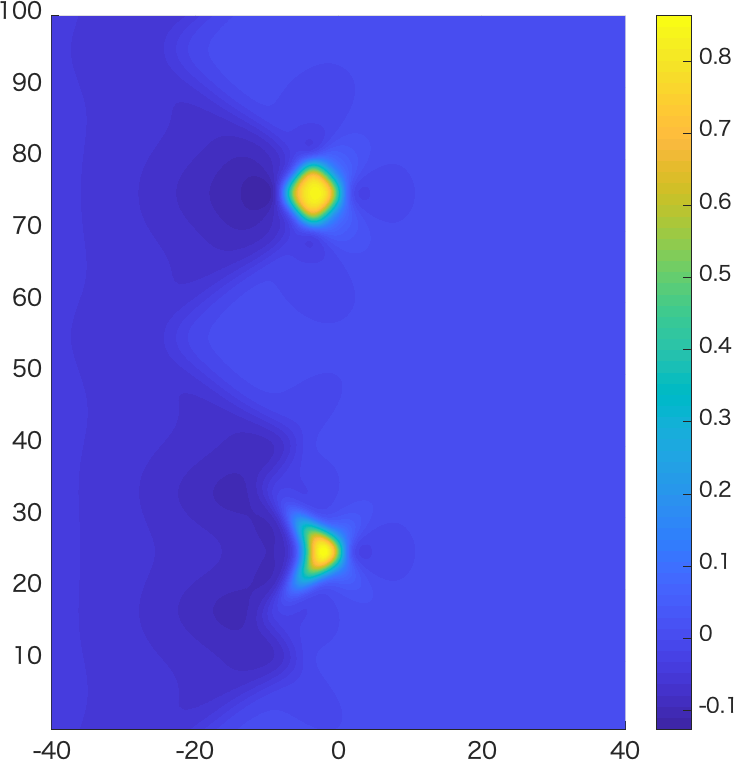}

            $t=300$
        \end{minipage}
        \caption{Time evolution in the bidomain FitzHugh--Nagumo equation.
        The first and the second rows are results for $\alpha=0.3$ and $\alpha=0.33$, respectively.
        We define other parameters as $a=0.9$, $b=0$, $\theta=\pi/4$, $\epsilon=10^{-3}$, and $\gamma=3$.}
        \label{fig:traveling_pulse}
    \end{figure}

    \subsection{Existence and stability of the planar pulse}
    \label{subsec:pulse}

\begin{figure}[tb]
        \centering
        \includegraphics[width=.6\hsize]{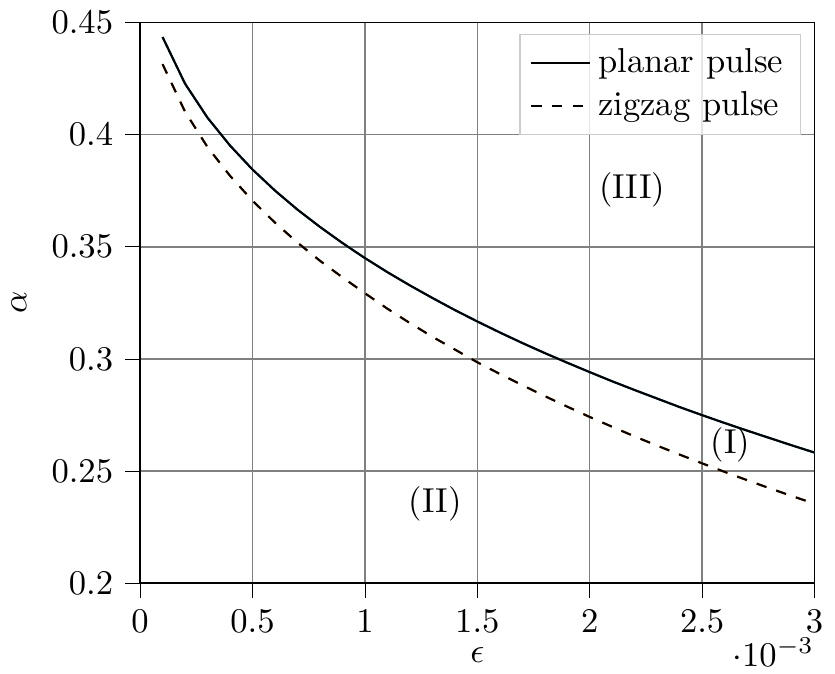}

        \caption{The existence region of planar pulse (below the solid line) and the one of zigzag pulse (below the broken line).}
        \label{fig:existence_pulse}
    \end{figure}

    In the bidomain Allen--Cahn equation, there always exists a planar front, whereas, in the bidomain FitzHugh--Nagumo equation, this is not the case.
    In the bidomain FitzHugh--Nagumo equation, depending on $(\epsilon,\alpha)$, the pulse may exist or vanish.
    As can be seen in \cref{fig:existence_pulse}, planar pulses exist stably in the regions (I) and (II) below the solid line, but not in the region (III) above.
    Also, the smaller $\alpha$ is, the larger the $\epsilon$-region where the planar pulse exists.
    
    First of all, let us ensure that if a planar pulse exists in the monodomain FitzHugh--Hagumo equation, it must also exist in the bidomain FitzHugh--Nagumo equation.
    Let $(u_{\mathrm{p}}^\theta,v_{\mathrm{p}}^\theta)$ be a planar pulse in the bidomain FitzHugh--Nagumo equation propagating in the direction of $\bm{n}^\theta$, and $c_{\mathrm{p}}^\theta$ be its velocity.
    The planar pulse satisfies the homogeneous Dirichlet boundary conditions at infinity:
    \begin{align*}
    	\lim_{|\xi|\to\infty}u_{\mathrm{p}}^\theta(\xi)
	=
	\lim_{|\xi|\to\infty}v_{\mathrm{p}}^\theta(\xi)
	=
	0.
    \end{align*}
    The planar pulse is a solution to the following boundary value problem:
    \begin{align*}
    	\begin{dcases}
		c_{\mathrm{p}}^\theta\frac{\mathrm{d}u_{\mathrm{p}}^\theta}{\mathrm{d}\xi}+Q(\bm{n}^\theta)\frac{\mathrm{d}^2u_{\mathrm{p}}^\theta}{\mathrm{d}\xi^2}+f(u_{\mathrm{p}}^\theta,v_{\mathrm{p}}^\theta)=0,\\
		c_{\mathrm{p}}^\theta\frac{\mathrm{d}v_{\mathrm{p}}^\theta}{\mathrm{d}\xi}+g(u_{\mathrm{p}}^\theta,v_{\mathrm{p}}^\theta)=0,\\
		\lim_{|\xi|\to\infty}u_{\mathrm{p}}^\theta(\xi)=\lim_{|\xi|\to\infty}v_{\mathrm{p}}^\theta(\xi)=0.
	\end{dcases}
    \end{align*}
    Let $(u_{\mathrm{p}}^*,v_{\mathrm{p}}^*)$ be the normalized planar pulse, and $c_{\mathrm{p}}^*$ be its speed; that is, $(u_{\mathrm{p}}^*,v_{\mathrm{p}}^*,c_{\mathrm{p}}^*)$ is a solution to the following boundary value problem:
    \begin{align*}
	\begin{dcases}
		c_{\mathrm{p}}^*\frac{\mathrm{d}u_{\mathrm{p}}^*}{\mathrm{d}\xi}+\frac{\mathrm{d}^2u_{\mathrm{p}}^*}{\mathrm{d}\xi}+f(u_{\mathrm{p}}^*,v_{\mathrm{p}}^*)=0,\\
		c_{\mathrm{p}}^*\frac{\mathrm{d}v_{\mathrm{p}}^*}{\mathrm{d}\xi}+g(u_{\mathrm{p}}^*,v_{\mathrm{p}}^*)=0,\\
		\lim_{|\xi|\to\infty}u_{\mathrm{p}}^*(\xi)=\lim_{|\xi|\to\infty}v_{\mathrm{p}}^*(\xi)=0.
	\end{dcases}
    \end{align*}
    Since there is no diffusion term in the equation for $v$ in the bidomain FitzHugh--Nagumo equation, we can construct the planar pulse by performing the same scaling
    \begin{align*}
    	u_{\mathrm{p}}^\theta(\xi)=u_{\mathrm{p}}^*\left(\xi/\sqrt{Q(\bm{n}^\theta)}\right),
	\quad
	v_{\mathrm{p}}^\theta(\xi)=v_{\mathrm{p}}^*\left(\xi/\sqrt{Q(\bm{n}^\theta)}\right),
	\quad
	c_{\mathrm{p}}^\theta=\sqrt{Q(\bm{n}^\theta)}c_{\mathrm{p}}^*
    \end{align*}
    as in the case of the bidomain Allen--Cahn equation.
    In other words, in the region of $(\epsilon,\alpha)$ where the planar pulse exists in the monodomain FitzHugh--Nagumo equation, the planar pulse always exists in the bidomain FitzHugh--Nagumo equation.
    
    Since we have found no difference between the monodomain and bidomain FitzHugh--Nagumo equations for planar pulses, our subsequent interest is in the region of $(\epsilon,\alpha)$ where zigzag pulses exist.
    For example, the first row of \cref{fig:traveling_pulse} shows the result when $\alpha=0.3$.
    When we add a small perturbation to the planar pulse, it becomes a zigzag pulse with coarsening and this state exists for a long time.
    However, in the second row, when $\alpha=0.33$, although coarsening occurs, the pulse is torn off before long and eventually disappears.
    These observations motivate us to examine the region of $(\epsilon,\alpha)$ where the zigzag pulse stably exists numerically.

    \begin{figure}
        \begin{minipage}{.25\hsize}
            \centering
            \includegraphics[width=\hsize]{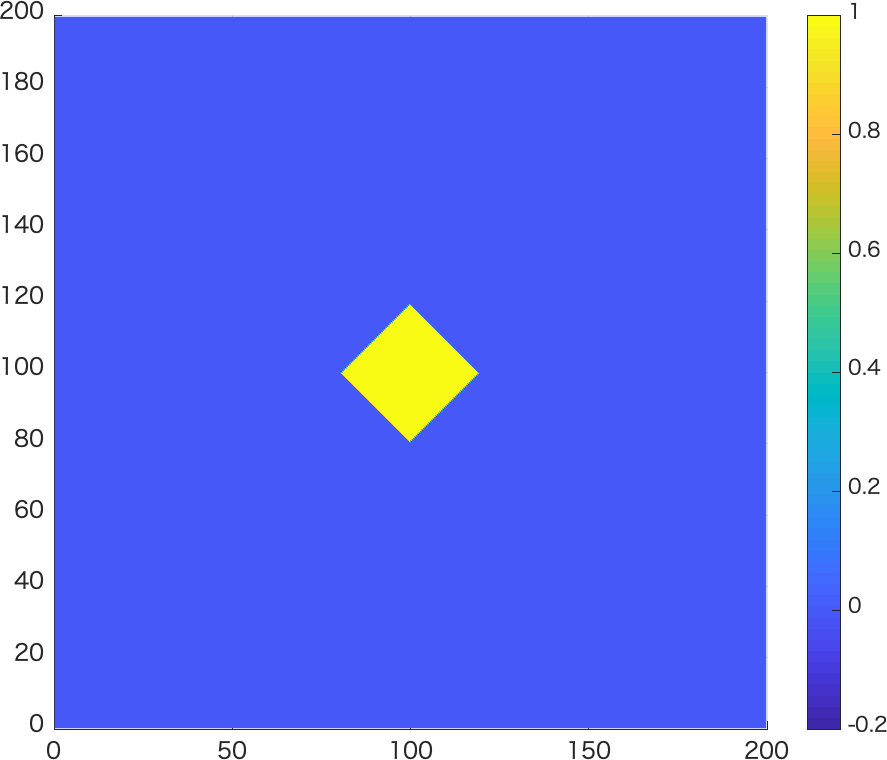}

            $t=0$
        \end{minipage}%
        \begin{minipage}{.25\hsize}
            \centering
            \includegraphics[width=\hsize]{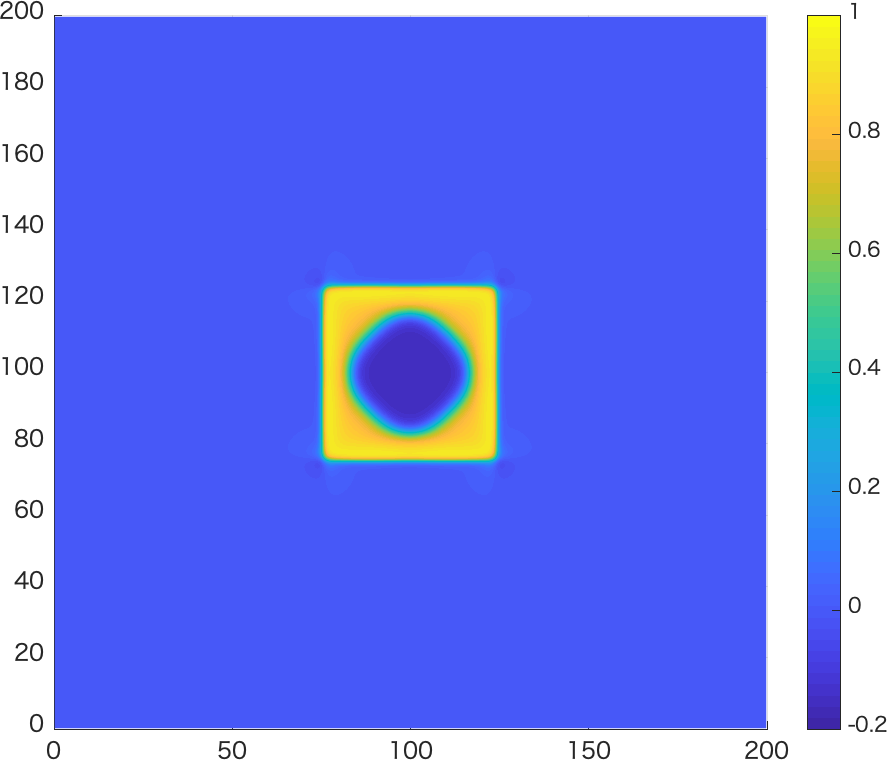}

            $t=280$
        \end{minipage}%
        \begin{minipage}{.25\hsize}
            \centering
            \includegraphics[width=\hsize]{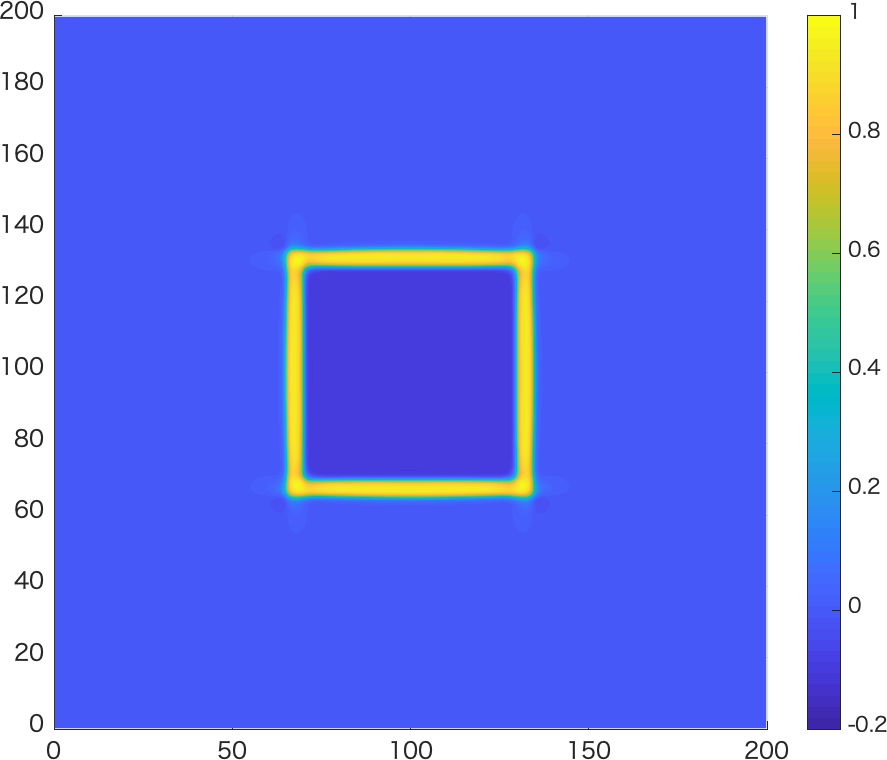}

            $t=600$
        \end{minipage}%
        \begin{minipage}{.25\hsize}
            \centering
            \includegraphics[width=\hsize]{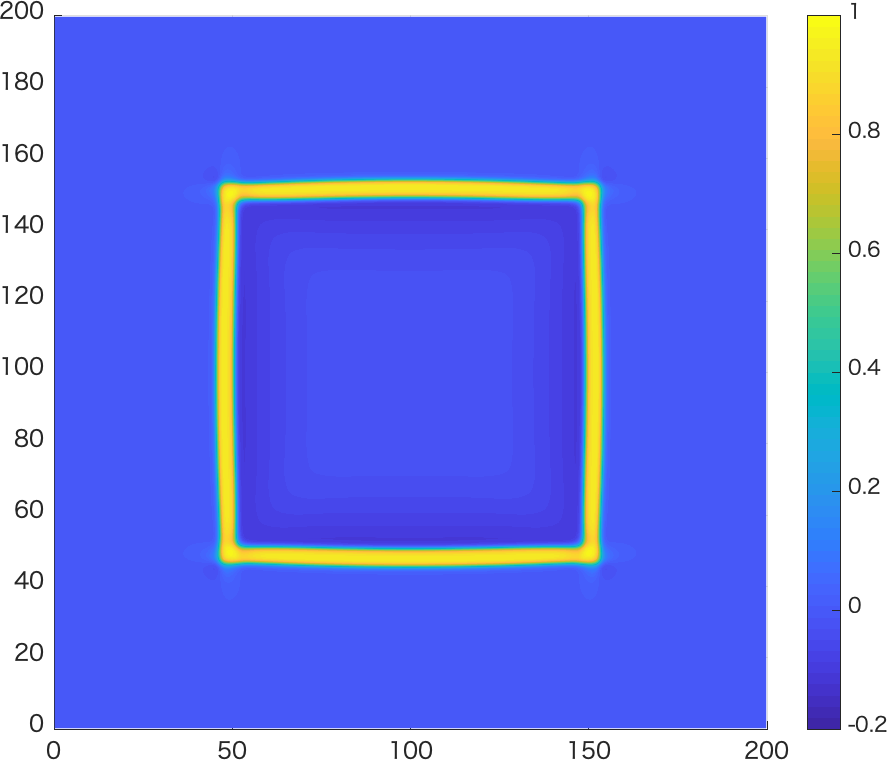}

            $t=1200$
        \end{minipage}%

        \begin{minipage}{.25\hsize}
            \centering
            \includegraphics[width=\hsize]{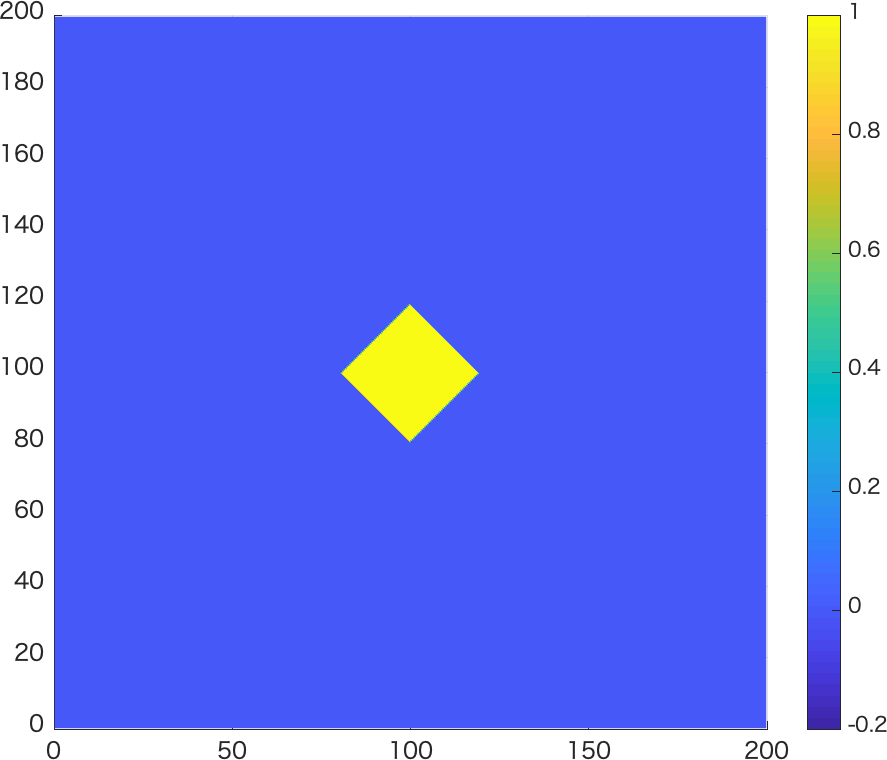}

            $t=0$
        \end{minipage}%
        \begin{minipage}{.25\hsize}
            \centering
            \includegraphics[width=\hsize]{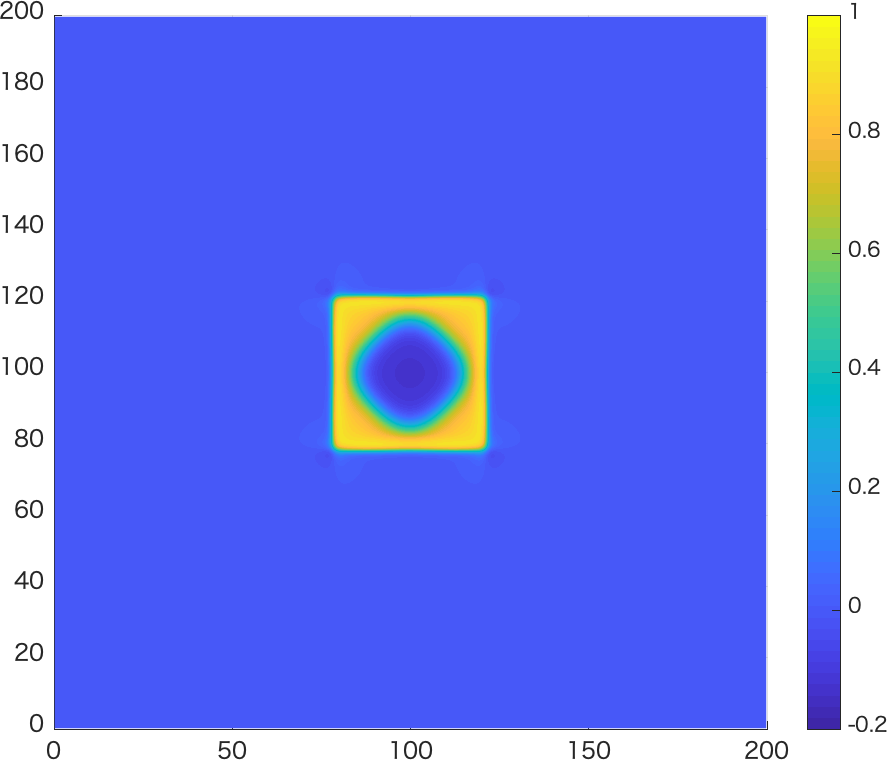}

            $t=260$
        \end{minipage}%
        \begin{minipage}{.25\hsize}
            \centering
            \includegraphics[width=\hsize]{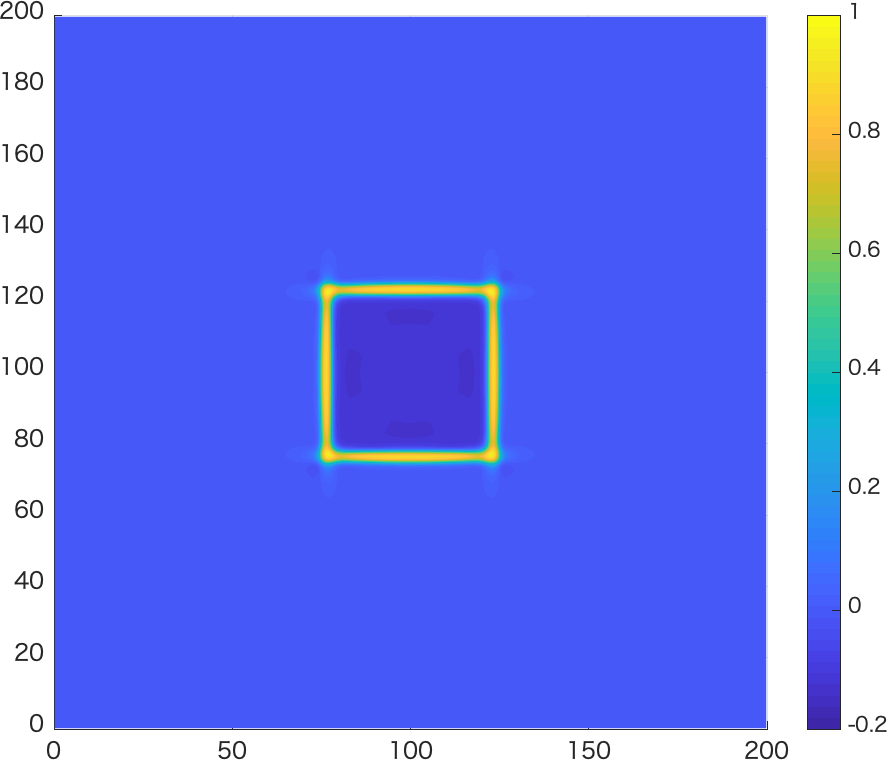}

            $t=400$
        \end{minipage}%
        \begin{minipage}{.25\hsize}
            \centering
            \includegraphics[width=\hsize]{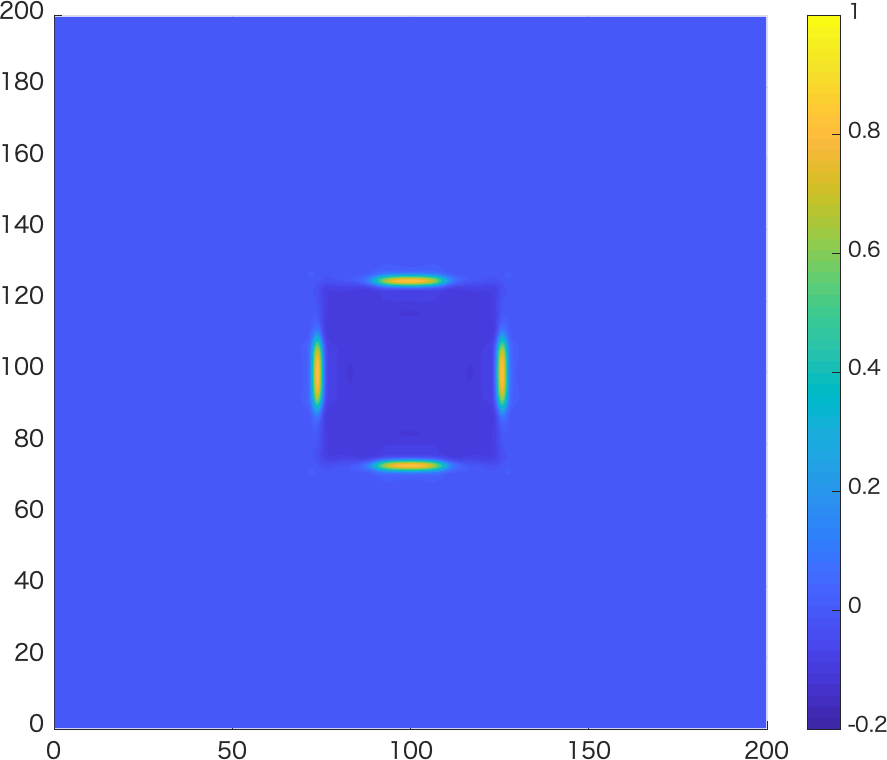}

            $t=500$
        \end{minipage}%
        \caption{Spreading pulse in the bidomain FitzHugh--Nagumo equation.
        The first row represents results for $\alpha=0.32$ and the second row for $\alpha=0.34$.
        We define other parameters as $a=0.9$, $b=0$, $\epsilon=10^{-3}$, and $\gamma=3$.}
        \label{fig:spreading_pulse}
    \end{figure}

    \Cref{fig:existence_pulse} shows the results of the numerical experiments, where we fix the parameters $a=0.9$, $b=0$, $\theta=\pi/4$, and $\gamma=3$ and vary $\alpha$ and $\epsilon$ as free parameters. We indicate the parametric regions $(\alpha,\epsilon)$ in which the planar pulses exist and/or zigzag pulses are stable. 
    The solid line indicates the boundary of existence/nonexistence of planar pulses, and the planar pulse exists below the solid line, i.e., in the regions (I) and (II). Note that this solid line is determined completely by the 1D FitzHugh-Nagumo model. As we saw above, planar pulse solutions of the bidomain FitzHugh-Nagumo model and the pulse solutions of the 1D FitzHugh-Nagumo equations are identical. The broken line shows whether the zigzag pulse is stable in the bidomain FitzHugh--Nagumo equation, and it is below the broken line, i.e., in the region (II).
    The above means that there are no stable planar or zigzag pulses in region (I). 
    
    We note that $\alpha$ controls the excitability of the FitzHugh-Nagumo system; the system is more excitable when $\alpha$ is close to $0$ and less excitable as $\alpha$ approaches $1/2$. The parameter $\epsilon$ controls the separation of time scales of the slow and fast variables. \Cref{fig:existence_pulse} indicates that the bidomain FitzHugh Nagumo pulse is prone to propagation failure when the system is less excitable and the separation of time scales is less pronounced.

    \subsection{Correspondence between spreading pulses and the Wulff shape}

    In the bidomain Allen--Cahn equation, we observed that the Wulff shape describes the spreading front's asymptotic shape.
    We see if the same consideration holds for the bidomain FitzHugh--Nagumo equation.
    We show our numerical results in \cref{fig:spreading_pulse}.
    We fix parameters as $a=0.9$, $b=0$, $\epsilon=10^{-3}$, and $\gamma=3$.
    There exists in this setting a threshold value $\alpha_*$ around $0.33$ on existence/nonexistence of zigzag pulse as in \cref{fig:existence_pulse}.
    Reflecting this observation, we can see that the spreading pulse exists and converges to the Wulff shape for $\alpha=0.32$ (the first row in \cref{fig:spreading_pulse}) and that it disappears for $\alpha=0.34$ (the second row in \cref{fig:spreading_pulse}).
    
    \section{Discussion}\label{sec:discussion}

    In this paper, we performed a detailed computational study of the asymptotic behavior of planar fronts and pulses of the bidomain Allen-Cahn and bidomain FitzHugh-Nagumo equation. For this purpose, we developed a numerical scheme that simulates the propagation of fronts and pulses on an infinite two-dimensional strip domain. We confirm that planar fronts of the bidomain Allen-Cahn equation are unstable when the Frank diagram is not convex in the direction of propagation. The destabilized planar fronts generically approach a zigzag rotating front whose shape and speed can be explained by the geometry of the Frank diagram. Destablization of the planar front thus takes place through a Hopf bifurcation. We have shown that this Hopf bifurcation can be either supercritical or subcritical depending on the parameter regime. For the bidomain FitzHugh-Nagumo equation, we have shown that the pulse solution does not necessarily develop into a zigzag rotating pulse. The pulse solution can entirely disappear especially when the parameters are close to the boundary of the parametric region in which planar pulses exist.

    Theoretical studies on arrhythmogenesis using FitzHugh-Nagumo and similar excitable systems have focused on the monodomain case, in which case the bidomain equation reduces to a reaction diffusion system (see \eqref{eq:monodomain_FHN}). In this paper, we have demonstrated that replacing the Laplacian in the monodomain model with the bidomain operator leads to qualitatively different asymptotic behaviors that could play an important role in arrythmogenesis. It is well-known that regions of cardiac ischemia or infarction induce cardiac arrhythmias, and the authors of \cite{CPS14} have suggested that the Frank diagram can become non-convex in regions of cardiac ischemia or healed infarcts. Our results suggest that planar fronts, passing through such a region, can become unstable and deform into a zigzag front. Furthermore, the zigzag front may maintain its zigzag shape even if it passes into a healthy region of the heart, as is suggested from our demonstration of a subcritical Hopf bifurcation. We have also demonstrated that zigzag pulses may completely disintegrate when passing through regions where the Frank diagram is non-convex and the underlying tissue is less excitable (see discussion at end of \cref{subsec:pulse}), leading to propagation failure. We note that regions of ischemia are precisely those locations for which the tissue is less excitable. These scenarios thus point to novel modes in which cardiac propagation can fail due solely to the bidomain nature of cardiac tissue. 
    
    \bibliographystyle{plain}
    \bibliography{2021_MMNS_arXiv}
\end{document}